# L'inauguration des séries trigonométriques dans la *Théorie analytique de la chaleur* de Fourier et dans la controverse des cordes vibrantes

**Alain Herreman**

Résumé – Cet article applique les notions de texte et d'énoncé inaugural à la *Théorie analytique de la chaleur* de Fourier et à la controverse des cordes vibrantes. Il établit que seul le livre de Fourier est un texte inaugural. Il propose en outre une analyse sémiotique des modes d'expression de la généralité qui permet dans chaque cas de clarifier les conditions d'introduction de la possibilité de représenter une fonction arbitraire par des séries trigonométriques.

Abstract (The inauguration of trigonometrical series in Fourier's *Théorie analytique de la chaleur* and in the vibrating strings controversy)
This paper applies the notions of inaugural text and statement to Fourier's *Théorie analytique de la chaleur* and to the vibrating strings controversy. It establishes that only Fourier's book is an inaugural text. It presents also a semiotic analysis of the expression modes of generality which allows to clarify in both cases the conditions of introducing of the possibility to represent an arbitrary function by trigonometrical series.

# I - Introduction[1]

La possibilité de représenter une fonction quelconque par une série trigonométrique a été envisagée au moins deux fois. Une première fois vers 1750 au cours de la controverse sur les cordes vibrantes qui opposa Daniel Bernoulli, D'Alembert, Euler et ensuite Lagrange. Une deuxième fois, par Fourier, plus de cinquante ans plus tard, dans le cadre de sa théorie de la chaleur. Nous voudrions ici étudier ces deux moments en les considérant dans la perspective d'une analyse des conditions d'expression de la généralité en mathématiques.
Les énoncés mathématiques (les théorèmes, leurs démonstrations, les définitions, les exemples, les figures, etc.) sont des énoncés la plupart du temps généraux. On peut étudier les moyens par lesquels leur généralité est exprimée ainsi que la manière dont ces moyens ont été historiquement introduits. Les séries trigonométriques en sont un exemple : elles offrent une représentation sinon de toutes les fonctions, tout du moins d'une vaste totalité d'entre elles, grâce à laquelle il est possible d'énoncer des théorèmes généraux. On sait aussi l'importance de cette représentation dans le développement de la théorie des fonctions tout au long du 19e siècle. Par ailleurs, la représentation des fonctions par des séries trigonométriques a été introduite en tant que représentation *conforme* des fonctions, c'est-à-dire pouvant représenter *toutes* les fonctions, tout ce qu'elle représente étant une fonction, toutes les propriétés des fonctions pouvant être obtenues à partir d'elle, et toutes ses propriétés étant elles-mêmes des propriétés des fonctions. C'est là une façon à la fois particulière et remarquable d'introduire une nouvelle représentation. On peut, pour les distinguer, appeler *textes inauguraux* les *textes* qui introduisent une telle représentation et qui s'attachent à en soutenir la conformité. On peut aussi distinguer les *énoncés* qui affirment qu'une représentation est conforme à ce qu'elle représente et les appeler des *énoncés inauguraux*. Nous avons montré que *La Géométrie* de Descartes était un texte inaugural qui comprenait plusieurs énoncés inauguraux [Herreman 2012]. Les notions d'énoncés et de textes inauguraux ont ainsi déjà été définies, présentées et discutées à cette occasion. Nous rappellerons ici ces définitions et nous montrerons que *La théorie analytique de la chaleur* de Fourier est aussi un texte inaugural. Nous montrerons en revanche que les textes de la controverse sur les cordes vibrantes ne forment pas un texte inaugural bien qu'ils comprennent l'énoncé inaugural correspondant. Nous esquisserons aussi une analyse des conditions d'expression de la généralité dans ces deux ensembles de textes. Nous préciserons d'abord le rôle des séries trigonométriques dans la généralité de la *Théorie analytique de la chaleur* de Fourier. Nous déterminerons ensuite inversement les conditions de possibilité sémiotiques de l'introduction de la représentation d'une fonction arbitraire par une série trigonométrique dans le texte de Fourier et lors de la controverse des cordes vibrantes. Nous pourrons ainsi rendre compte, suivant ce point de vue, de l'introduction de la possibilité de représenter une fonction arbitraire par des séries trigonométriques dans les deux cas.

La *Théorie analytique de la chaleur* et la controverse des cordes vibrantes ont déjà fait l'objet d'importants travaux de la part d'historiens des mathématiques et de la physique spécialistes des 18e et 19e siècles[2]. Notre propos est ici seulement de les considérer suivant une perspective différente : celle des moyens d'expression de la généralité et des conditions de leur introduction. Nous avons aussi essayé de nous en tenir à ce seul point de vue qui ne saurait bien sûr prétendre se substituer à aucun autre. Il conduit d'ailleurs souvent à des observations qui ont déjà été faites (le cas échéant par les auteurs des textes considérés). Nous avons essayé de l'indiquer toutes les fois où cela nous est apparu et paru justifié. L'intérêt dans ce cas est d'y être arrivé par une autre voie. Il nous a semblé qu'il pouvait aussi dans quelques cas conduire à nuancer d'autres analyses. Nous avons là aussi essayé de l'indiquer.

# II - L'inauguration des séries trigonométriques dans la *Théorie analytique de la chaleur* de Joseph Fourier

## 1 - Introduction : explication du cheminement

Nous allons rappeler dans cette partie la définition d'un énoncé et d'un texte inaugural et quelques-uns de leurs enjeux à partir de l'exemple de la représentation d'une fonction arbitraire par une série trigonométrique. Nous donnerons ensuite brièvement les éléments de l'exposé de Fourier dans la *Théorie analytique de la chaleur* nécessaires à la suite. Nous commencerons notre analyse en déterminant la contribution de la généralité de la représentation par des séries trigonométriques à la généralité de la théorie de la chaleur. Nous le ferons en précisant le statut de l'équation générale de la chaleur et celui des équations particulières données pour les divers cas de solide. Nous déterminerons ensuite les conditions d'introduction dans ce texte de la possibilité de représenter une fonction arbitraire par des séries trigonométriques. Nous nous intéresserons ensuite à l'inauguration proprement dite en mettant surtout en évidence la relation entre certaines caractéristiques de l'exposé de Fourier et les problèmes sémiotiques inhérents à une telle inauguration.

---

2 Pour une présentation et des commentaires du développement de la théorie de la chaleur de Fourier voir Grattan-Guinness & Ravetz [1972], Herivel [1975], Charbonneau [1976], Friedman [1977], Grattan-Guinness [1990, vol. 2, chp. 9 ; 2000 ; 2005], Dhombres & Robert [1998, 443-620]. Parmi les travaux traitant de la controverses des cordes vibrantes on peut citer : Riemann [1854] ; Burkhardt [1908, 10-47 trad. fr. in Jouve 2007, II, 21-51] ; Langer [1947] ; Truesdell [1960] ; Ravetz [1961] ; Grattan-Guinness [1970, 2-21] ; Youschkevitch [1976] ; Lützen [1982, 15-24] ; Darrigol [2007] ; Jouve [2007, I, 23-25 ; 2008]. L'édition annotée des mémoires de D'Alembert sur le sujet se trouve dans Jouve [2007]. Sur les problèmes vibratoires avant cette controverse, voir Cannon & Dostrovsky [1981]. Sur D'Alembert et les équations aux différences partielles, voir Guilbaud & Jouve [2009], Demidov [1982, 1989].

## 2 - Enoncés et textes inauguraux

### a) Un énoncé inaugural

Dans l'avant dernier article du chapitre III de la *Théorie analytique de la chaleur*, Fourier conclut les développements qui précèdent par l'énoncé suivant :

> «si l'on propose une fonction $f\,x$ dont la valeur est représentée, dans un intervalle déterminé, depuis x=0 jusqu'à x=X, par l'ordonnée d'une ligne courbe tracée arbitrairement, on pourra toujours développer cette fonction en une série qui ne contiendra que les sinus ou les cosinus, ou les sinus et cosinus des arcs multiples, ou les seuls cosinus des multiples impairs. » art. 235

Fourier soutient dans cet énoncé qu'il est toujours possible de développer en série trigonométrique une fonction dont la valeur est représentée sur un intervalle donné par une courbe définie arbitrairement[3]. On peut distinguer dans cet énoncé fondamental un ensemble remarquable de caractéristiques. C'est un énoncé relatif aux fonctions mais plus précisément un énoncé *général* sur les fonctions : il se rapporte aux fonctions dans leur *totalité*[4]. L'affirmation centrale selon laquelle « *on pourra toujours développer cette fonction en une série qui ne contiendra que les sinus ou les cosinus etc.* » implique *deux totalités* : d'une part, la totalité des fonctions « *dont la valeur est représentée, dans un intervalle déterminé (...) par l'ordonnée d'une ligne courbe tracée arbitrairement* », d'autre part, celle des « *séries trigonométriques* ». Ces deux totalités existent bien toutes les deux. Chacune peut être considérée sans que l'autre ne le soit. Si l'on considère maintenant les *moyens d'expression* de ces totalités, c'est-à-dire les moyens d'expression ou les représentations qui font de chaque totalité *une* totalité, il y en a aussi deux : d'une part l'expression géométrique d'une courbe sur un intervalle (soumise à la condition de définir une fonction), d'autre part l'expression analytique d'une série trigonométrique. Ces deux totalités et leur mode d'expression n'ont pas le même statut. Les courbes (même soumises à la condition de correspondre à une fonction[5]) ont, au moment où Fourier écrit, un caractère *pré-établi*. Cela ne veut pas dire qu'elle serait bien définie, mais que Fourier ne saurait prétendre introduire lui-même cette totalité en tant que telle : c'est une totalité qu'il considère de fait comme donnée : *il y a* des courbes. Leur extension est notamment fixée par les phénomènes naturels qu'elles représentent[6]. C'est une totalité dont *il dispose*. Il n'en est pas tout à fait de même des séries trigonométriques. Il est vrai néanmoins qu'elles ont déjà été considérées lors de la controverse sur les cordes vibrantes. Mais justement, elles l'ont été dans une controverse ni close, ni même tout à fait éteinte si l'on pense à Lagrange parti prenant de cette controverse et toujours interlocuteur de Fourier. Au moment où Fourier écrit les séries trigonométriques *n*'ont *pas* le statut d'une totalité reçue. Et surtout, elles ne l'ont pas dans le rapport considéré ici : à savoir représenter n'importe quelle fonction. C'est bien l'enjeu de cet énoncé que de soutenir ce rapport. La correspondance considérée entre les fonctions et les séries trigonométriques a un sens très fort : non seulement les séries trigonométriques permettent de représenter *toutes* les courbes, mais il n'y a pas de série trigonométrique qui ne soit pas la représentation d'une courbe. Ces expressions, nouvelles au regard des courbes qu'elles représentent (et peu importe ici que l'on en attribue la nouveauté à Daniel Bernoulli), n'excèdent d'aucune manière ce qu'elles représentent : la représentation n'introduit pas de courbes « imaginaires », « fausses », « pathologiques » ou « formelles ». Au-delà de cette correspondance entre objets ou expressions, il y a aussi une correspondance entre leurs propriétés : pour Fourier les séries trigonométriques restituent *toutes* les *propriétés* aussi bien des courbes que des phénomènes de propagation thermique qu'elles représentent[7]. Les deux totalités en présence sont aussi de nature différente : elles sont hétérogènes. On peut simplement repérer cette hétérogénéité en remarquant que les expressions des deux totalités ne se combinent pas : on peut les *juxtaposer*, par exemple en écrivant une série trigonométrique sous sa représentation graphique, mais si l'on insère dans une série trigonométrique un morceau de courbe on obtient une expression mixte qui n'appartient à aucune des deux totalités[8]. Cette hétérogénéité, jointe aux caractéristiques précédentes, rend en fait impossible la démonstration de cet énoncé[9].

---

[3] Nous nous limitons ici à la représentation par des séries trigonométriques. Pour la représentation par la formule intégrale, introduite par Fourier comme extension de la précédente, voir Annaratone [1997].

[4] Le terme de « totalité », moins marqué mathématiquement, est préféré à celui d' « ensemble ». Il réduit le risque d'introduire subrepticement une vision ensembliste, et avec elle une uniformisation anachronique des totalités considérées.

[5] Cette condition sera sous-entendue dans la suite.

On peut ainsi distinguer cinq caractéristiques de cet énoncé, auxquelles il est pratique d'associer un nom (indiqué entre parenthèses), que l'on peut reformuler brièvement comme suit :

---

6 Il pourra ainsi écarter les fonctions prenant une valeur infinie en un point : « *il est impossible qu'aucune question naturelle conduise à supposer que la fonction fx devient infinie lorsqu'on donne à x une valeur singulière comprise entre des limites données.* » art. 417. Plus loin : « *Si l'on conçoit qu'entre les limites de l'intégration certaines valeurs de* α *deviennent infinies, la construction indique dans quel sens la proposition générale doit être entendue. Mais nous ne considérons point ici les cas de cette nature, parce qu'ils n'appartiennent point aux questions physiques.* » art. 423.
7 Cette caractéristique ne peut, en l'occurrence, être directement dégagée de l'énoncé. Elle est néanmoins avérée et nous l'indiquons ici parce qu'elle étend la correspondance sur les expressions qui vient d'être considérée.
8 Cette observation ne sert qu'à faire reconnaître qu'*il y a* de l'hétérogénéité, ce qui suffit ici à notre propos. L'hétérogénéité a un caractère historique multiple qui interdit d'en donner une caractérisation générale pertinente : suivant les systèmes d'expressions utilisés, les nombres et les lettres ont pu être ou non hétérogènes mais ensuite constituer les expressions d'un même système d'expression. C'est là, bien sûr, un aspect important de l'histoire de l'algèbre.
9 Il n'entre ici aucune conception *figée* de la démonstration. Il s'agit seulement d'indiquer les problèmes sémiotiques que rencontre *nécessairement* la démonstration d'un tel énoncé.

1. l'énoncé met en jeu deux totalités (dualisme) ;
2. l'une des totalités est tenue pour pré-établie (réalisme[10]) ;
3. l'autre totalité n'est pas tenue pour pré-établie sous la forme considérée ou dans le rapport considéré à la totalité pré-établie (inauguration) ;
4. l'énoncé affirme que la deuxième totalité, celle qui n'est pas pré-établie, est une représentation conforme de la première (conformité) ;
5. la démonstration de la conformité des deux totalités est impossible (incommensurabilité) ;

Ces cinq conditions peuvent servir à caractériser un *énoncé inaugural*[11]. L'énoncé dans lequel Fourier soutient la possibilité de représenter toutes les courbes par des séries trigonométriques est ainsi un énoncé inaugural, de même que les diverses formulations qui en ont été données dans ses mémoires précédents[12].

### *b) De l'énoncé inaugural au texte inaugural*

Pour soutenir cet énoncé Fourier doit établir la possibilité de représenter par une série trigonométrique *toute* courbe représentant une fonction sur un intervalle. Il faut pour cela considérer *toutes* les courbes et les *parcourir* d'une manière ou d'une autre. Or, *la possibilité de parcourir une totalité dépend directement du système d'expression par lequel elle est exprimée*. Grâce à l'expression « générale » d'une série trigonométrique, cela est remarquablement facile pour les séries trigonométriques[13]. Il est dès lors très facile en disposant d'une telle expression de faire un raisonnement valant pour *toutes* les séries trigonométriques. Il n'y a rien de semblable pour les courbes au moment où Fourier écrit. La démonstration d'un énoncé portant sur toutes les courbes doit *nécessairement* comprendre un procédé ou un dispositif pour résoudre ce *problème sémiotique*[14]. Certaines différences entre les deux systèmes d'expression qui font l'intérêt de la représentation trigonométrique sont des problèmes au moment de l'inauguration. Nous avons vu Descartes recourir dans une situation semblable à la possibilité, introduite pour la circonstance, de représenter une courbe par l'instrument qui permettait de la tracer. En s'appuyant sur le type spécifique d'instruments qu'il retenait, il reportait sur eux la possibilité de soutenir que les relations entre les points tracées par leur moyen ont un rapport à une même droite exprimée par une équation polynomiale. Pour les problèmes de géométrie, il tirait parti des caractéristiques remarquables, et en partie fausses..., offertes par le problème de Pappus. Les problèmes rencontrés dépendent étroitement des systèmes d'expressions en présence ; ceux rencontrés pour établir la possibilité d'associer une courbe à une série trigonométrique sont eux-mêmes différents de ceux rencontrés dans le sens inverse. Les systèmes d'expression en présence rendent parfois possible cette démonstration (il s'agit alors d'un *théorème de représentation*). Les énoncés inauguraux sont ceux pour lesquels cela est au contraire impossible. Mais c'est aussi cette impossibilité qui fait leur importance historique : l'intérêt du système d'expression inauguré est d'introduire une représentation tenue à la fois pour conforme mais néanmoins irréductible à la précédente. Les problèmes que posent la démonstration de ces énoncés sont précisément ceux dont l'usage de la représentation inaugurée nous affranchira ensuite. Celui qui soutient un énoncé inaugural est *nécessairement* confronté à ces problèmes. Il devra dès lors nécessairement trouver des expédients pour pallier l'impossibilité de parcourir la totalité des expressions du système d'expressions pré-établi et leurs propriétés ainsi que les différences entres les deux systèmes d'expressions en présence (dualité). Ces problèmes déterminent un type particulier de texte que nous appelons un *texte inaugural*.

---

10 Quand nous parlons de « réalisme » aussi bien à propos de *La Géométrie* de Descartes que de la *Théorie analytique de la chaleur* de Fourier nous ne sommes pas en train de confondre les philosophies, explicites ou implicites, de leurs auteurs. Nous désignons simplement *le fait* sémiotique que leurs textes présupposent des totalités pré-établies : celui de Descartes présuppose la totalité pré-établie des problèmes de géométrie et celle des instruments qu'il décrit et celui de Fourier la totalité pré-établie des courbes définies sur un intervalle. Il ne s'agit pas de confondre la nature de l'existence des problèmes de géométrie et celle des instruments dans *La Géométrie*. On ne saurait *a fortiori* les confondre avec celle des courbes définies sur un intervalle dans la *Théorie analytique de la chaleur*. Il ne s'agit donc pas de le faire et nous ne le faisons pas. En parlant dans tous ces cas de « réalisme », nous ne faisons pas de la philosophie, nous désignons un fait sémiotique.
11 Une discussion plus précise de ces conditions est donnée dans Herreman [2012, 71-79].
12 On peut citer comme autres énoncés : « *Il n'y a ainsi aucune fonction fx, ou partie de fonction, que l'on ne puisse exprimer en une suite trigonométrique.* » Fourier 1822, art. 503. « *Il résulte de mes recherches sur cet objet que les fonctions arbitraires même discontinues peuvent toujours être représentées par les développements en sinus ou cosinus d'arcs multiples, et que les intégrales qui contiennent ces développements sont précisément aussi générales que celles où entrent les fonctions arbitraires d'arcs multiples.* » Fourier 1807/1972, 183 ; « *L'exposition détaillée des résultats de notre analyse ne peut laisser aucun doute sur le véritable sens dans lequel ils doivent être pris ; les développements de sinus et de cosinus multiples ont évidemment toute la généralité que comporte les fonctions arbitraires.* » Fourier 1807/1972, 113 ; « *une fonction quelconque peut toujours être développée en séries de sinus ou de cosinus d'arcs multiples* » Fourier 1807/1972, 115 ; « *cette remarque est essentielle, en ce qu'elle conduit à connaître comment les fonctions entièrement arbitraires peuvent aussi être développées en séries de sinus d'arcs multiples.* » Fourier 1807, voir aussi Fourier 1822, p. 234. "*Les théorèmes dont il s'agit donnent le moyen de réduire une fonction arbitraire et même discontinue en séries trigonométriques. Cette nouvelle extension de la théorie des suites était nécessaire pour autoriser l'emploi des solutions particulières des équations aux différences partielles.*» Fourier, "Notes jointes à l'extrait du mémoire sur la chaleur" in Hérivel 1980, 59.

Ces problèmes peuvent faire l'objet de remarques explicites. Fourier en fait peu. Il trouve tout de même « *remarquable que l'on puisse exprimer par des séries convergentes (...) des lignes et des surfaces qui ne sont pas assujetties à une loi continue* » (art. 230). Ces problèmes se manifestent surtout par le dispositif argumentatif qu'il développe pour soutenir son énoncé. Il avertit ainsi à plusieurs reprises son lecteur, et à chaque fois dans des passages relatifs à cette conformité, qu'il est amené à faire des développements inhabituels[15]. Ils sont surtout

---
13 Cela suppose de *disposer* de la possibilité de parcourir *tous* les nombres, ce qui renvoie en l'occurrence à une autre inauguration : cette inauguration en présuppose une autre (et en fait, plusieurs autres).

14 Le parti pris de cette analyse est de montrer la part de ces problèmes sémiotiques dans le texte de Fourier. Il va de soi que des problèmes de bien d'autres natures ont aussi leur part dans ce texte !

15 « *Nous aurions regardé comme inutile les développements contenus dans les articles précédents, si nous n'avions point à exposer une théorie entièrement nouvelle, dont il est nécessaire de fixer les principes.* » (art.197, p. 198). Citons cet autre passage dans lequel il annonce devoir déroger au mode d'exposition habituel des mémoires scientifiques en promettant de revenir à une plus grande concision dans ses publications ultérieures : « *Nous avons démontré dans cet Ouvrage tous les principes de la Théorie de la chaleur, et résolu toutes les questions fondamentales. On aurait pu les exposer sous une forme plus concise, omettre les questions simples, et présenter d'abord les* indirectement abordés dans son « Discours préliminaire », sur lequel nous reviendrons, au travers des considérations philosophiques plus générales sur le rapport de l'Analyse mathématique aux phénomènes naturels.

## 3 - La *Théorie analytique de la chaleur* de Joseph Fourier : un texte inaugural

### *a) L'inauguration dans le texte*

Nous allons maintenant présenter les éléments de l'exposé de Fourier utiles à la suite de notre analyse. Nous commencerons par rappeler la structure d'ensemble du texte pour ensuite présenter plus en détail le chapitre III qui comprend l'inauguration des séries trigonométriques.

La *Théorie analytique de la chaleur* est composée de neuf chapitres précédés d'un « Discours préliminaire »[16]. L'introduction qui suit ce « Discours préliminaire » en est le premier chapitre. Fourier y définit les grandeurs physiques considérées (température, capacité spécifique, *C*, conductibilité extérieure, *h*, et intérieure, *K*) et leur mesure, le modèle et les principes adoptés pour décrire la diffusion de la chaleur à l'intérieur et à la surface d'un corps et analysée au niveau des molécules. L'étude du flux de chaleur à l'intérieur d'un solide homogène compris entre deux plans infinis dont la température est fixée lui permet d'obtenir l'expression du flux d'abord à une dimension unique ( $F = -K \dfrac{dv}{dz}$ , section IV, arts. 68-72), puis à trois dimensions à l'intérieur d'un parallélipipède ( $-K\dfrac{dv}{dx}, -K\dfrac{dv}{dy}, -K\dfrac{dv}{dz}$ , section VII, arts. 85-99[17]). Le chapitre II est consacré à déterminer l'équation aux dérivées partielles de la chaleur à l'intérieur d'un solide. Il l'établit les équations de la chaleur, à l'intérieur et à la surface, pour différents solides : une armille, une sphère solide, un cylindre solide, un parallélipipède de longueur infinie, un cube (sections I-V).

---

*conséquences les plus générales ; mais on a voulu montrer l'origine même de la Théorie et ses progrès successifs. Lorsque cette connaissance est acquise, et que les principes sont entièrement fixés, il est préférable d'employer immédiatement les méthodes analytiques les plus étendues, comme nous l'avons fait dans les recherches ultérieures. C'est aussi la marche que nous suivrons désormais dans les Mémoires qui seront joints à cet Ouvrage, et qui en forment en quelque sorte le complément, et par là nous aurons concilié, autant qu'il peut dépendre de nous, le développement nécessaire des principes avec la précision qui convient aux applications de l'Analyse. »* Fourier 1822, Discours préliminaire, Oeuvres I, xxiv-xxv.

16 Charbonneau [1976] et Dhombres & Robert [1998, chp. VIII] contiennent une présentation et un commentaire continus croisés avec le manuscrit de 1807. Charbonneau [1976] propose une analyse génétique limitée aux aspects mathématiques mais prenant en compte l'ensemble des manuscrits non publiés. Dhombres & Robert [1998] proposent une analyse physico-mathématique, biographique et épistémologique plus complète. Grattan-Guinness & Ravetz [1972] comparent le manuscrit de 1807, qu'ils éditent, avec les exposés ultérieurs de la théorie. Fourier [1822, 602-637] donne lui-même une présentation détaillée de son cheminement dans la « Table des matières contenues dans cet ouvrage » située à la fin du livre.
17 La numérotation des articles passe de 90 à 84, puis de 87 à 95, et deux articles portent ensuite le numéro 100. La numérotation donnée dans la « Table des matières » ne correspond pas toujours à celle du texte.

Fourier précise qu'il préfère recourir aux équations particulières qu'il vient toutes d'établir, et qu'il résoudra ensuite chacune à leur tour, plutôt qu'à l'*équation générale* pour « *les corps solides d'une figure quelconque* » :

> « *On pourrait former les équations générales qui représentent le mouvement de la chaleur dans les corps solides d'une figure quelconque, et les appliquer aux cas particuliers. Mais cette méthode entraîne quelquefois des calculs assez compliqués que l'on peut facilement éviter. Il y a plusieurs de ces questions qu'il est préférable de traiter d'une manière spéciale, en exprimant les conditions qui leur sont propres ; nous allons suivre cette marche et examiner séparément les questions que l'on a énoncées dans la première section de l'introduction ; nous nous bornerons d'abord à former les équations différentielles, et nous en donnerons les intégrales dans les chapitres suivants.* » art. 101

Il établit néanmoins ensuite l'équation générale de la chaleur à *l'intérieur* d'un solide *quelconque* en considérant un parallélipipède rectangle dont chacune des six faces a une température déterminée constante, pour lequel il retrouve les expressions du flux obtenues au chapitre précédent ( $-K\frac{dv}{dx}, -K\frac{dv}{dy}, -K\frac{dv}{dz}$ , Section VI, art. 138), et obtient ainsi l'équation générale $\frac{dv}{dt} = \frac{K}{C.D}(\frac{d^2v}{dx^2} + \frac{d^2v}{dy^2} + \frac{d^2v}{dz^2})$ (art.142)[18]. Puis il détermine l'équation générale pour la *surface*[19] ( $m\frac{dv}{dx} + n\frac{dv}{dy} + p\frac{dv}{dz} + \frac{h}{k}vq = 0$ où $m, n, p, q$ sont des fonctions connues des coordonnées des points de la surface, art. 147). Il montre ensuite comment les équations pour *l'intérieur* d'un cylindre et d'une sphère se déduisent de l'équation générale (arts. 155-6).

Après ces deux premiers chapitres, toute la suite du livre, chapitres III à IX, est consacrée à la détermination de la distribution de la chaleur dans les divers solides étudiés en résolvant le système d'équations différentielles propre à chacun, le neuvième et dernier chapitre considérant une masse solide homogène dont toutes les dimensions sont infinies (et donc sans bord).

Le chapitre III, le premier de cette série de chapitres, étudie la chaleur dans une plaque rectangulaire semi-infinie. C'est au cours de ce chapitre que Fourier soutient et introduit son énoncé inaugural sur la représentation des fonctions par des séries trigonométriques. Nous en ferons donc une présentation plus précise.

---

[18] « *L'équation à laquelle on est parvenu dans la question précédente, représente le mouvement de la chaleur dans l'intérieur de tous les solides. Quelle que soit en effet la forme du corps, il est manifeste qu'en le décomposant en molécules prismatiques, on obtiendra ce même résultat. On pourrait donc se borner à démontrer ainsi l'équation de la propagation de la chaleur. Mais afin de rendre plus complète l'exposition des principes, et pour que l'on trouve rassemblés dans un petit nombre d'articles consécutifs les théorèmes qui servent à établir l'équation générale de la propagation dans l'intérieur des solides, et celles qui se rapportent à l'état de la surface, nous procéderons, dans les deux sections suivantes, à la recherche de ces équations, indépendamment de toute question particulière, et sans recourir aux propositions élémentaires que nous avons expliquées dans l'introduction.* » art. 131

[19] « *Si le solide a une forme déterminée, et si la chaleur primitive se dissipe successivement dans l'air atmosphérique entretenu à une température constante, il faut ajouter à l'équation générale et à celle qui représente l'état initial, une troisième condition relative à l'état de la surface. Nous allons examiner dans les articles suivants, la nature de l'équation qui exprime cette dernière condition.* » art. 146.

La lame semi-infinie étudiée a sa température à la base maintenue à 1, celle à ses deux côtés à 0 et sa température initiale à l'intérieur est supposée nulle (voir figure). Il s'agit de déterminer la distribution de la chaleur à l'état permanent.

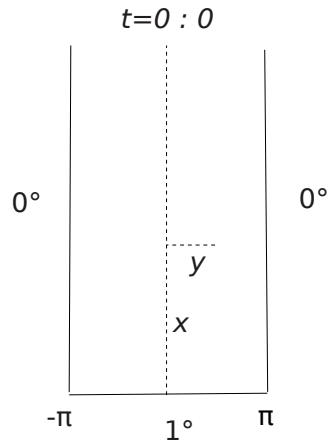

L'état permanent est déterminé par l'équation différentielle générale à *l'intérieur* de la lame, $\frac{\partial^2 v}{\partial x^2}+\frac{\partial^2 v}{\partial y^2}=0$ , les conditions aux bords et l'état initial[20]. En séparant les variables, Fourier obtient la *solution générale* de l'équation pour l'intérieur de la lame :

$$v = a e^{-x}\cos.y + b e^{-3x}\cos.3y + c e^{-5x}\cos.5y + d e^{-7x}\cos.7y + \ldots$$

Pour déterminer les coefficients *a*, *b*, *c*, etc de la solution générale, il utilise la température à la base fixée à 1 qui lui donne la relation suivante (*y* compris entre - π et π) :

$$1 = a\cos.y + b\cos.3y + c\cos.5y + d\cos.7y + \ldots$$

Déterminer les coefficients de la solution générale revient ainsi à développer la fonction constante égale à 1 en une série trigonométrique. Avant de consacrer la section suivante à la résolution de ce problème mathématique, Fourier termine cette section par l'interprétation physique de la solution générale (art. 170). La section II est entièrement consacrée à trouver le développement en série de cosinus de 1 en résolvant un système infini d'équations linéaires infinies à coefficients constants (arts 171-178). Le même résultat est obtenu dans section III en tronquant la série en une somme partielle et en passant à la limite (arts. 179-189). Fourier identifie ensuite les fonctions ( $\frac{x}{2}$ , $\log(2\cos(\frac{1}{2}x))$ , etc.) représentées par divers développements en séries de cosinus ou de sinus (arts. 179-184) et détermine les intervalles sur lesquels ces représentations sont valables (arts 185-188). Dans la section IV, il revient au problème physique et montre que l'expression de la solution trouvée permet de « *déterminer toutes les circonstances du mouvement permanent de la chaleur dans une lame rectangulaire échauffée à son origine* » (art. 192). Il donne ensuite l'expression de la distribution de la chaleur pour une base de longueur *2l* (art. 196)

$$v = \frac{4A}{\pi}(e^{-\frac{\pi x}{2l}}\cos\frac{\pi y}{2l} - \frac{1}{3}e^{-3\frac{\pi x}{2l}}\cos 3\frac{\pi y}{2l} + \frac{1}{5}e^{-5\frac{\pi x}{2l}}\cos 5\frac{\pi y}{2l} - \frac{1}{7}e^{-7\frac{\pi x}{2l}}\cos 7\frac{\pi y}{2l} + ...)$$

dont il donne aussi une interprétation physique[21] (arts 197-200). Il établit l'unicité de la solution (arts 200-204) en décomposant la lame en deux autres lames. Il revient au bilan des échanges moléculaires et fait appel entre autres *hypothèses et évidences physiques* au fait que la température finale doit être nulle (art. 201). Dans la section V, qui ne comprend que deux articles, Fourier remplace l'expression, comprenant une série, de la solution trouvée par une expression finie (art. 205).

---

20 Herivel [1975, 150] remarque que les conditions aux bords choisies évitent le problème de la diffusion de la chaleur par la surface (en l'occurrence, le bord) du solide. Sur le choix d'une température initiale nulle, voir Herivel [1975, 171-172].

La section VI, avant-dernière de ce chapitre et de loin la plus longue (51 pages), est celle consacrée au « *Développement d'une fonction arbitraire en séries trigonométriques* » (arts 207-235) :

> « on n'a traité qu'un seul cas d'un problème plus général, qui consiste à développer une fonction quelconque en une suite infinie de sinus ou de cosinus d'arcs multiples. Cette question est liée à la théorie des équations aux différences partielles et a été agitée dès l'origine de cette analyse. Il était nécessaire de la résoudre pour intégrer convenablement les équations de la propagation de la chaleur ; nous allons en exposer la solution. » art. 207

L'étude de la lame infinie dont la température à la base est maintenue égale à 1 a conduit au développement en série trigonométrique de la fonction constante égale à 1 (et donc de n'importe quelle constante). La lame et sa chaleur disparaissent à nouveau pour laisser la place à des considérations purement mathématiques et à la démonstration qu'une fonction $\varphi x$ impaire donnée par son développement de Taylor en 0 peut être développée en une série de sinus (arts 207-215).

Fourier résout pour cela un nouveau système infini d'équations linéaires infinies à coefficients constants et obtient les coefficients de la série trigonométrique sous la forme de séries de terme général les dérivées successives de $\varphi x$ en 0 (formule (A), art. 215), qu'il réduit ensuite pour obtenir l'expression suivante (formule (B), art. 217) :

$$\frac{1}{2}\pi\varphi x = \sin.x\left\{\varphi\pi - \frac{1}{1^2}\varphi^{II}\pi + \frac{1}{1^4}\varphi^{IV}\pi - \frac{1}{1^6}\varphi^{VI}\pi + etc.\right\}$$
$$-\frac{1}{2}\sin.2x\left\{\varphi\pi - \frac{1}{2^2}\varphi^{II}\pi + \frac{1}{2^4}\varphi^{IV}\pi - \frac{1}{2^6}\varphi^{VI}\pi + etc.\right\}$$
$$+\frac{1}{3}\sin.3x\left\{\varphi\pi - \frac{1}{3^2}\varphi^{II}\pi + \frac{1}{3^4}\varphi^{IV}\pi - \frac{1}{3^6}\varphi^{VI}\pi + etc.\right\}$$
$$-\frac{1}{4}\sin.4x\left\{\varphi\pi - \frac{1}{4^2}\varphi^{II}\pi + \frac{1}{4^4}\varphi^{IV}\pi - \frac{1}{4^6}\varphi^{VI}\pi + etc.\right\}$$
$$+ etc.$$

---

21 A titre d'exemple : « *Il est facile de voir, soit au moyen de cette équation, soit d'après l'art. 171, que la chaleur se propage dans ce solide, en s'éloignant de plus en plus de l'origine, en même temps qu'elle se dirige vers les faces infinies B et C.* » art. 197.

A ce stade, Fourier a *démontré* que les fonctions considérées pouvaient être développées en séries de sinus d'arcs multiples. Il doit maintenant aller au-delà :

> « On peut étendre les mêmes conséquences à des fonctions quelconques, même à celles qui seraient discontinues et entièrement arbitraires. Pour établir clairement la vérité de cette proposition, il est nécessaire de poursuivre l'analyse qui fournit l'équation précédente (B) et d'examiner quelle est la nature des coëffiçents qui multiplient sin.x, sin.2x, sin.3x, sin.4x. » (art. 219)

En remarquant que les coefficients de la formule précédente sont solution d'une équation différentielle linéaire du deuxième ordre, il les exprime sous forme intégrale et obtient la formule intégrale du développement en sinus (formule (D), art. 219) :

$$\frac{1}{2}\pi \varphi x = \sin.x \int \sin.x . \varphi x . dx + \sin.2x \int \sin.2x \varphi x \, dx$$
$$+ \sin.3x \int \sin.3x \varphi x \, dx \ldots + \sin.ix \int \sin.ix \varphi x \, dx + etc.$$

A partir de cette formule, la possibilité de développer en série de sinus une fonction impaire donnée par son développement de Taylor en 0 est étendue à une fonction arbitraire :

> « Cette remarque [c'est-à-dire le fait que les coefficients de la série trigonométrique soient donnés par l'expression intégrale] est importante, en ce qu'elle fait connaître comment les fonctions entièrement arbitraires peuvent aussi être développées en séries de sinus d'arcs multiples. En effet, si la fonction $\varphi x$ est représentée par l'ordonnée variable d'une courbe quelconque, dont l'abscisse s'étend depuis x=0 jusqu'à x=π, et si l'on construit sur cette même partie de l'axe la courbe trigonométrique connue dont l'ordonnée est y=sin.x, il sera facile de se représenter la valeur d'un terme intégral. Il faut concevoir que, pour chaque abscisse x à laquelle répond une valeur de $\varphi x$ et une valeur de sin.x, on multiplie cette dernière valeur par la première, et qu'au même point de l'axe on élève une ordonnée proportionnelle au produit $\varphi x . \sin.x$. On formera, par cette opération continuelle, une troisième courbe dont les ordonnées sont celles de la courbe trigonométrique, réduites proportionnellement aux ordonnées de la courbe arbitraire qui représente $\varphi x$ Cela posé, l'aire de la courbe réduite, étant prise depuis x=0 jusqu'à x=π, donnera la valeur exacte du coefficient de sin x ; et, quelle que puisse être la courbe donnée qui répond à $\varphi x$ soit qu'on puisse lui assigner une équation analytique, soit qu'elle ne dépende d'aucune loi régulière, il est évident qu'elle servira toujours à réduire d'une manière quelconque la courbe trigonométrique ; en sorte que l'aire de la courbe réduite a, dans tous les cas possibles, une valeur déterminée qui donne celle du coefficient de sin x dans le développement de la fonction. Il en est de même du coefficient suivant *b* ou $\int (\varphi x . \sin.2x . dx)$. » art. 220.

Fourier explique ici en détails comment obtenir à partir de la *courbe* de $\varphi x$ les coefficients de son développement en série de sinus : $\int_0^\pi \varphi x . \sin.x . dx$, $\int_0^\pi \varphi x . \sin.2x . dx$, etc. Ce n'est qu'ensuite qu'il utilise l' «orthogonalité» ( $\int_0^\pi \sin.ix . \sin.jx . dx = 0$ pour i≠j, $\int_0^\pi \sin.ix . \sin.ix . dx = \frac{1}{2}\pi$ ) pour établir que les coefficients d'une fonction arbitraire $\varphi x = a_1 \sin.x + a_2 \sin.2x + a_3 \sin.3x + \ldots a_i \sin.ix + \ldots etc$ développable en série de sinus sont donnés par la formule (art. 221) :

$$a_i = \frac{\int \varphi\, x.\sin.ix.dx}{\frac{1}{2}\pi}\ .$$

Cette formule est ensuite utilisée pour obtenir le développement en série de sinus de la fonction constante[22] (art. 222) et de la fonction *cos x* afin d'« *apporter un exemple qui ne laisse aucun doute sur la possibilité de ce développement* » (art. 223).

L' « orthogonalité » lui permet d'obtenir ensuite le développement en série de *cosinus* d'une fonction *arbitraire* (formule (n), art. 224) :

$$\frac{1}{2}\pi\, \varphi\, x = \frac{1}{2}\int_0^\pi \varphi\, x\, dx + \cos.x \int_0^\pi \varphi\, x \cos.x\, dx$$
$$+ \cos.2x \int_0^\pi \varphi\, x \cos.2x\, dx + \cos.3x \int_0^\pi \varphi\, x \cos.3x\, dx + etc.$$

Et il peut affirmer :

> « Ce théorème et le précédent conviennent à toutes les fonctions possibles, soit que l'on en puisse exprimer la nature par les moyens connus de l'analyse, soit qu'elles correspondent à des courbes tracées arbitrairement. » art. 224

Cette formule est à nouveau utilisée pour obtenir le développement en série de cosinus des fonctions x, sin.x (art. 225), puis celui de fonctions discontinues, comme la fonction égale à $\frac{\pi}{2}$ pour $0 \leq x \leq \alpha$ et nulle pour $\alpha < x \leq \pi$, et la fonction égale à sin.x pour $0 \leq x \leq \alpha$ et nulle pour $\alpha < x \leq \pi$ (art. 226), et ensuite à une fonction définie par un arc de courbe et une ligne droite (art. 227), un contour de trapèze (art. 228) et la surface d'une pyramide (art. 229).

Après quelques remarques sur ces séries, Fourier déduit le développement en série de sinus et de cosinus d'une fonction définie entre -π et π (art. 233), puis les développements en séries trigonométriques sur des intervalles de longueur quelconque :

(P)
$$r\, f\, x = \frac{1}{2}\int_{-r}^r f\, x\, dx + \cos.(\pi\frac{x}{r})\int_{-r}^r f\, x \cos.(\pi\frac{x}{r})dx + \cos.(2\pi\frac{x}{r})\int_{-r}^r f\, x \cos.(2\pi\frac{x}{r})dx + etc.$$
$$+ \sin.(\pi\frac{x}{r})\int_{-r}^r f\, x \sin.(\pi\frac{x}{r})dx + \sin.(2\pi\frac{x}{r})\int_{-r}^r f\, x \sin.(2\pi\frac{x}{r})dx + etc.$$

(N)
$$r\, f\, x = \frac{1}{2}\int_0^r f\, x\, dx + \cos.(\pi\frac{x}{r})\int_0^r f\, x \cos.(\pi\frac{x}{r})dx + \cos.(2\pi\frac{x}{r})\int_0^r f\, x \cos.(2\pi\frac{x}{r})dx + etc.$$

(M)
$$r\, f\, x = \sin.(\pi\frac{x}{r})\int_0^r f\, x \sin.(\pi\frac{x}{r})dx + \sin.(2\pi\frac{x}{r})\int_0^r f\, x \sin.(2\pi\frac{x}{r})dx + etc.$$

Cette section se conclut par l'énoncé inaugural, déjà cité, pour des fonctions définies sur un intervalle borné de longueur quelconque :

> « Il résulte de tout ce qui a été démontré dans cette Section concernant le développement des fonctions en séries trigonométriques que, si l'on propose une fonction *fx* dont la valeur est représentée, dans un intervalle déterminé, depuis x=0 jusqu'à x=X, par l'ordonnée d'une ligne courbe tracée arbitrairement, on pourra toujours développer cette fonction en une série qui ne contiendra que les sinus ou les cosinus, ou les sinus et cosinus des arcs multiples, ou les seuls cosinus des multiples impairs. » art. 235.

---

[22] Fourier dit ici avoir déjà donné ce développement, mais il s'agissait, comme on l'a vu, d'un développement en série de cosinus.

Dans la dernière section du chapitre, section VII, Fourier reprend le problème de la distribution de la chaleur dans la lame avec cette fois une distribution à la base définie par une *fonction arbitraire fx*.

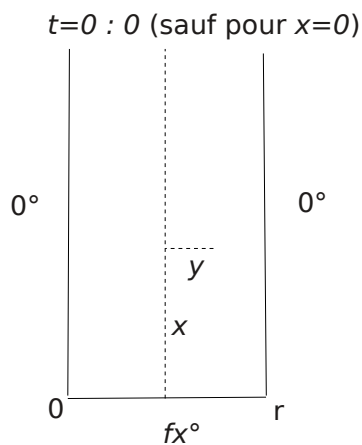

La relation (M) permet de déterminer dans ce cas les coefficients du développement en série trigonométrique de *fx*, et donc l'expression du régime permanent dans la lame pour une distribution à la base quelconque (art. 236)[23].

**b) Les séries trigonométriques et la théorie de la chaleur**

Dans la *Théorie analytique de la chaleur* Fourier inaugure à la fois une théorie de la chaleur et une représentation générale des fonctions. Seule la deuxième inauguration en fait un texte inaugural suivant notre acception[24], mais Fourier doit aussi présenter une théorie de la chaleur ayant une certaine *généralité*. Nous allons commencer par préciser la manière dont la généralité de la représentation par des séries trigonométriques intervient dans la généralité de la théorie. Nous le ferons en précisant le statut de l'équation générale de la chaleur établie pour un solide quelconque et celui des équations particulières établies pour les divers solides considérés.

---

23 Lors de l'étude de l'armille, Fourier déterminera l'expression de la fonction initiale à partir de l'expression de la distribution de chaleur en fonction du temps, retrouvant ainsi la possibilité de développer une fonction arbitraire en série trigonométrique (arts. 278 et 279).

24 La question n'est pas de savoir à qui l'on accorde le titre de « texte inaugural » mais d'introduire des distinctions pertinentes permettant de repérer des différences et des similitudes objectives. L'acception retenue privilégie l'inauguration de la représentation d'une totalité comme condition de théorèmes auparavant impossibles. Tel n'est pas l'enjeu de l'inauguration d'une théorie de la chaleur. Ce qui ne lui retire rien de son importance!

### *i) L'équation générale et les équations particulières*

L'équation *générale* désigne l'équation qui régit la propagation de la température à *l'intérieur* d'un solide *quelconque*. Or, l'exposé de la théorie n'est pas fondé sur l'étude d'un solide *quelconque* et sur la résolution de l'équation *générale* mais sur l'étude de solides *particuliers* et la résolution de leurs équations *particulières*[25]. Il s'agit d'expliciter ce qui détermine Fourier à développer sa théorie à partir de *cas de solide*, et de leurs équations *particulières* (et qui doivent être considérées comme telles même quand elles coïncident avec l'équation *générale*), plutôt que de la développer à partir de l'équation *générale*, qu'il a établie, pour un solide *quelconque* et de présenter ses différents cas de solides comme autant d'applications du cas général. Fourier, on l'a vu, remarque à ce propos que cela lui évite « *des calculs assez compliqués* » et, sans plus de précisions, qu' « *il est préférable de traiter de manière spéciale* » certaines questions (art. 101, déjà cité). Voyons plus précisément ce qu'il en est.
On commencera pour cela par rappeler le schéma de la description mathématique de la propagation de la chaleur dans un solide. Ce schéma comprend :

1) un système de deux équations différentielles :
   a) une équation différentielle décrivant la propagation à l'*intérieur* du solide[26] ;
   b) une équation différentielle décrivant la propagation à la *surface* du solide[27] ;
2) une fonction décrivant la distribution initiale de la chaleur dans le solide.

Pour obtenir la fonction décrivant la propagation de la chaleur dans le solide, en fonction du temps ou en régime permanent, Fourier commence par résoudre le système des deux équations différentielles. Il obtient ainsi la *solution générale* de ce système d'équations différentielles *particulières* sous la forme d'une série infinie (par séparation des variables et linéarité des équations). La forme du terme général de la série dépend alors des équations considérées. Chaque terme de la série comprend des constantes *a, b, c*, etc. qu'il reste à déterminer pour obtenir l'expression de la propagation de la chaleur dans le solide. Ces constantes sont déterminées grâce à la distribution initiale de la chaleur dans le solide.
Le calcul intégral et différentiel permet d'avoir une équation *générale* qui dépend des constantes physiques (la capacité spécifique, *C*, la conductibilité extérieure, *h*, et intérieure, *K*), mais qui est indépendante de la forme du solide. Toute la spécificité *géométrique* de l'objet est ainsi renvoyée à sa surface. L'équation générale à *la surface* dépend elle étroitement de la forme du solide et elle est inévitablement un obstacle à une description *générale* de la propagation de la chaleur [Dhombres & Robert 1998, 521]. Ceci suffit à rendre nécessaire un traitement au *cas par cas*. Mais cela ne suffit pas encore à rendre compte de la démarche de Fourier. En effet, il pourrait se dispenser de recourir aux équations *particulières* décrivant la chaleur à l'*intérieur* des solides et considérer plutôt l'équation *générale*. Or, il détermine au contraire toujours la solution *générale* à partir de l'équation *particulière*. Quand il cherche par exemple la distribution de la chaleur dans l'armille, il ne considère pas l'équation *générale* $\frac{dv}{dt} = \frac{K}{C.D}(\frac{d^2v}{dx^2} + \frac{d^2v}{dy^2} + \frac{d^2v}{dz^2})$ mais l'équation *particulière* $\frac{dv}{dt} = \frac{K}{C.D}\frac{d^2v}{dx^2} - \frac{h.l}{C.D.S}v$. Pourquoi? Parce que seule l'équation *particulière* peut donner une description *conforme* de la propagation de la chaleur dans un solide (voir ci-dessous sur la conformité). Partir de l'équation générale conduirait en effet à la *même* expression de la solution générale, c'est-à-dire à la même expression de la distribution de chaleur à l'intérieur de *tous* les solides. Or, l'expression de la distribution de chaleur dans un corps est donnée par l'expression de la solution générale (ces deux expressions ne diffèrent que par les coefficients, indéterminés dans un cas, déterminés dans l'autre). Fourier, on l'a vu, ne manque pas d'interpréter thermiquement cette solution générale (c'est aussi, inversement, un contrôle de la validité de la solution)[28]. La conformité des *valeurs* de la fonction intervient bien sûr dans la vérification expérimentale mais ces valeurs, exprimées par une série infinie, sont d'abord généralement impossibles à calculer exactement, et surtout, la conformité ne se réduit pas à cette coïncidence numérique : avant tout calcul, c'est (la forme de) *l'expression* de la solution qui, pour Fourier, représente le phénomène physique[29]. C'est l'*expression* de la solution qui a un sens physique et qui est interprétée par Fourier. La validité de la théorie est fondée sur la conformité de l'*expression* des solutions aux phénomènes. La conformité des valeurs en est une conséquence. Si l'on considère l'expression de la solution générale obtenue à partir de l'équation générale, la distribution de chaleur à l'intérieur des corps serait indépendante de la forme du corps : la conformité de l'expression de la distribution de la chaleur, même interne, n'est garantie et ne peut être obtenue qu'à partir de l'équation *particulière* du solide La conformité de la théorie, telle qu'elle est conçue par Fourier, exclut de recourir à l'équation générale, et donc : la généralité de la théorie ne peut en rien venir de la généralité de l'équation générale.

---

25 C'est nous qui introduisons le syntagme «équation particulière » que Fourier n'emploie pas. Il utilise « l'équation qui exprime le mouvement de la chaleur dans [le solide] ». Il emploie bien en revanche le syntagme « équation générale » (avec le pluriel pour désigner les deux équations différentielles pour l'intérieur et pour la surface d'un solide).
26 Cette équation peut être l'équation *générale*, c'est-à-dire la même pour tous les cas de solide, ou *particulière,* c'est-à-dire déterminée à partir des caractéristiques géométriques du solide considéré. Fourier adopte systématiquement l'équation *particulière*.
27 Le cas d'un solide de dimensions infinies traité au chapitre IX est à cet égard plus simple puisqu'il fait disparaître les conditions relatives à la surface (art. 461). Il existe aussi, on l'a vu, une équation différentielle générale pour la surface, mais celle-ci fait bien sûr intervenir les fonctions décrivant les coordonnées de celle-ci.

Voyons-en à présent les conséquences sur le statut de la représentation d'une fonction arbitraire par une série trigonométrique. Les séries trigonométriques sont introduites lors de l'étude de la diffusion de la chaleur dans la lame. On vient de voir que l'équation interne déterminait, en tant qu'équation *particulière*, l'expression *particulière* de la solution générale interne, et donc le terme général ou la *base* (le terme, bien sûr, n'est pas de Fourier) de la série, en l'occurrence *cos.nx*. La distribution de la chaleur dans *chaque* cas de solide n'est donc *pas* et ne peut pas être exprimée par une série trigonométrique mais par une série dont la forme du terme général dépend de l'*équation particulière* du solide (et qui n'est trigonométrique que si l'équation particulière coïncide avec l'équation générale). L'équation générale ne jouant aucune rôle dans la généralité de la théorie, *le fait que les séries trigonométriques interviennent dans la résolution de cette équation ne leur confère aucun rôle privilégié dans la théorie*. L'équation générale ne rend donc pas non plus compte du rôle pourtant privilégié de ces séries dans la théorie. Confondre l'équation particulière avec l'équation générale, même quand leurs expressions sont les mêmes, peut conduire à donner à la représentation trigonométrique une généralité qu'elle n'a pas dans la *Théorie analytique de la chaleur*.

---

28 A comparer avec Dhombres & Robert [1998, 528] qui écrivent à ce propos «*Mais hésitant à procéder directement à l'explication qui est attendue, c'est-à-dire à fournir un procédé de calcul, il passe sans transition à la formulation générale. Comme si le problème était déjà résolu avec en place de la constante 1 une fonction quelconque, ou plutôt comme si la résolution effective importait peu. Parler de maladresse logique, c'est manquer le sens même de la démarche : une forme générale a été obtenue avec l'écriture* [ $1 = a\cos y + b\cos 3y + c\cos 5y + d\cos 7y + ...$ ] *et c'est son sens qu'il faut d'abord acquérir. Il s'agit pour Fourier de montrer non seulement que la progression de la pensée est bonne, c'est-à-dire orientant la suite, mais aussi et surtout que l'analyse a effectivement trouvé la voie du réel* ».

29 L'aspect expérimental a été rendu marginal dans l'exposé de la théorie dans la *Theorie analytique de la chaleur*.

Introduire ou reconnaître ici un caractère arbitraire à l'expression de la solution générale revient à introduire ou à reconnaître un caractère arbitraire à toute interprétation de cette expression : la conformité de la solution *exclut un choix libre* de la base. La conformité exclut l'arbitraire[30]. Le rapport des mathématiques aux phénomènes se retrouve (nécessairement) dans les caractéristiques sémiotiques des expressions utilisées et l'épistémologie dans le sémiotique.

La lame est le premier des six solides étudiés et sa distribution de chaleur est déterminée par une équation différentielle qui n'est autre que l'équation générale (réduite à deux dimensions). La résolution de cette équation, avec les conditions aux bords choisies, conduit à une distribution de chaleur exprimée par une série trigonométrique. On pourrait dès lors penser avoir là la raison du rôle des séries trigonométriques dans la théorie générale : ces séries expriment la distribution de chaleur dans une lame, associée à l'équation générale, elles devraient donc pour cette raison servir à exprimer la distribution de chaleur dans les autres cas de solide. On vient de voir que ça n'était pas le cas. L'intérêt de commencer par l'étude de la lame ne vient donc pas de son rapport à l'équation générale. Il nous faut donc encore comprendre à la fois pourquoi Fourier commence par ce cas de solide et le rôle dans la théorie générale de la représentation des fonctions par des séries trigonométriques. Les deux questions sont liées.

---

30 Fourier l'exprime on ne peut plus clairement : « *Les fonctions que l'on obtient par ces solutions sont donc composées d'une multitude de termes, soit finis, soit infiniment petits : mais la forme de ces expressions n'a rien d'arbitraire ; elle est déterminée par le caractère physique du phénomène. C'est pourquoi, lorsque la valeur de la fonction est exprimée par une série où il entre des exponentielles relatives au temps, il est nécessaire que cela soit ainsi, parce que l'effet naturel dont on recherche les lois, se décompose réellement en parties distinctes, correspondantes aux différents termes de la série. Ces parties expriment autant de mouvements simples compatibles avec les conditions spéciales ; pour chacun de ces mouvements, toutes les températures décroissent en conservant leurs rapports primitifs. On ne doit pas voir dans cette composition un résultat de l'analyse dû à la seule forme linéaire des équations différentielles, mais un effet subsistant qui devient sensible dans les expériences. Il se présente aussi dans les questions dynamiques où l'on considère les causes qui anéantissent le mouvement ; mais il appartient nécessairement à toutes les questions de la théorie de la chaleur, et il détermine la nature de la méthode que nous avons suivie pour les résoudre.* » art. 428-7°. Fourier remarque qu'il pourrait aussi prendre, au lieu des fonctions trigonométriques, « une infinité d'autres fonctions » pour exprimer une fonction arbitraire : « *On parvient ainsi à exprimer une fonction arbitraire f(x) sous diverses formes très remarquables ; mais nous ne faisons point usage de ces transformations dans la recherche spéciale qui nous occupe.* » art. 423. Voir aussi art. 424 et art. 428-4° et 10°.

**ii) L'expression de l'état initial**

La représentation des fonctions par des séries trigonométriques a bien un rôle privilégié dans l'ensemble de la *Théorie analytique de la chaleur,* Fourier indique lui-même explicitement lequel :

> « On ne peut résoudre entièrement les questions fondamentales de la théorie de la chaleur, sans réduire à cette forme [série trigonométrique] les fonctions qui représentent l'état initial des températures. » art. 235[31]

Les séries trigonométriques servent en effet à exprimer l'*état initial* des températures. On peut vérifier que même dans les cas de solide pour lesquels la solution générale du système d'équations différentielles *n*'est *pas* exprimée par une série trigonométrique, c'est-à-dire pour la sphère (art. 290) et le cylindre (arts 310, 311), la température *initiale* est exprimée par une série trigonométrique, ce qui permet à Fourier de déterminer les coefficients de la solution générale (ce que rien *a priori* ne garantit, comme cela ressort de la diversité des méthodes employées à ces endroits par Fourier). Exprimer sous forme de série trigonométrique la distribution initiale, qui n'est autre que la solution générale considérée pour t=0, permet de fait à chaque fois, mais de manière *ad hoc*, de déterminer les coefficients de la solution générale (alors même, il faut y insister, que celle-ci n'est pas elle-même toujours exprimée par une série trigonométrique).

D'un point de vue physique, la distribution de chaleur dépend de sa distribution initiale dans le solide. Pour déterminer la distribution de la chaleur dans des solides *particuliers* mais avec une distribution initiale *quelconque* il est nécessaire de disposer d'une expression pour une fonction quelconque. Mais cette expression doit aussi permettre de trouver la suite infinie des coefficients qui figurent dans l'expression de la solution générale. Les séries *entières* ont une infinité de coefficients indéterminés, mais elles ne permettent pas d'exprimer une fonction quelconque (ni de résoudre les équations que l'on obtiendrait). Les séries trigonométriques font elles les deux : elles offrent une représentation générale des fonctions qui permet aussi la résolution de ces systèmes d'équations. Faire l'un sans l'autre n'aurait guère d'intérêt. Les séries trigonométriques satisfont des conditions de nature différente : sémiotiques (représentation), mathématique (résolution d'un système d'équations) et physique (interprétation thermique de leur expression).

---

31 On peut encore citer : « *Ces transformations des fonctions en suites trigonométriques sont des éléments de la Théorie analytique de la chaleur ; il est indispensable d'en faire usage pour résoudre les questions qui dépendent de cette théorie.* » art. 362. A propos de l'armille Fourier écrit : « *Pour former [la forme] qui exprime le mouvement de la chaleur dans une armille, il était nécessaire de résoudre une fonction arbitraire en une série de sinus et cosinus d'arcs multiples* » art. 282. La représentation d'une fonction arbitraire par une série trigonométrique intervient aussi, mais cette fois subrepticement, dans le cas du prisme rectangulaire et du cube. Voir dans l'édition des œuvres la note de Darboux aux articles 322 et 335. Pour le chapitre IX, voir en particulier arts 352 et 362. C'est *a contrario* parce que Daniel Benroulli n'a pas correctement pris en compte l'état initial que sa solution n'était, selon Fourier, pas satisfaisante : « *L'imperfection de cette solution consiste, ce me semble, en ce que l'auteur ne déterminait point les coëfficients par la comparaison avec l'état initial, et se contentait d'alléguer qu'ils pouvaient l'être.* » Fourier 1807/1972, 252.

Fourier justifie de commencer par l'étude de la lame parce qu'elle est « *le premier exemple de l'analyse qui conduit à ces solutions* [complètes et d'une application facile] » et qu'elle lui paraît « *plus propre qu'aucune autre à faire connaître les éléments de la méthode que nous avons suivie*» (art. 163). L'histoire des sciences est elle-même invoquée : « [L]*a considération des questions simples et primordiales est un des moyens les plus certains de découvrir les lois des phénomènes naturels, et nous voyons, par l'histoire des Sciences, que toutes les théories se sont formées suivant cette méthode.* » (art. 164). Et plus loin : « *l'ordre et l'espèce des questions ont été tellement choisis que chacune d'elles présentât une difficulté nouvelle et d'un degré plus élevé* » (art. 341). Tout cela est évidemment juste, mais ne rend néanmoins pas complètement compte de la position de cette étude : à ces raisons pédagogiques et épistémologiques il faut ajouter une nécessité sémiotique[32].

La lame n'est pas seulement, suivant les justifications données par Fourier, le cas de solide dont l'étude présenterait le moins de difficultés. La nécessité de commencer ne vient pas non plus, on l'a vu, de son rapport à l'*équation générale* mais de la nécessité de disposer d'une représentation d'une fonction arbitraire. La généralité de la représentation par des séries trigonométriques est ainsi un élément constitutif de la généralité de sa théorie[33], elle en est une condition de possibilité : sans elle, la distribution de la chaleur ne pourrait être déterminée que pour des conditions initiales *particulières* (et qui devrait encore être données sous une forme permettant de déterminer les coefficients cherchés). Si l'exposé de la théorie est fondé sur l'étude de solides *particuliers*, grâce à cette représentation d'une fonction arbitraire Fourier peut néanmoins s'affranchir de la contrainte de ne considérer que des distributions initiales elles-mêmes *particulières*. La contribution de cette représentation à la généralité de la théorie a ainsi été

---

32 Remarquons que dans l'étude de la lame ce n'est pas la distribution initiale, en l'occurrence nulle à l'intérieur, qui est représentée par une série trigonométrique mais la distribution à la base. C'est là une autre singularité de l'étude de la lame par rapport aux autres solides.

33 Fourier [1829/oeuvres II, 175] souligne lui-même que « *des températures initiales arbitraires (...) augmente beaucoup la généralité de la question*» et plus loin qu' « *il est nécessaire de ne point particulariser l'état initial. En effet, l'état qui se forme après que la continuité s'est établie dépend lui-même et très prochainement de la disposition initiale, qui est entièrement arbitraire. La continuité est compatible avec une infinité de formes qui différeraient extrêmement du paraboloïde; et l'on ne peut pas restreindre à cette dernière figure celle du petit corps immergé sans latérer, dans ce qu'elle a d'essentiel, la généralité de la question.* » Fourier [1829/oeuvres II, 178-179].

précisée[34]. Fourier, à nouveau, le dit aussi clairement[35] :

> « Pour que ces solutions [des équations générales de la propagation de la chaleur] fussent générales et qu'elles eussent une étendue équivalente à celle de la question, il était nécessaire qu'elles pussent convenir avec l'état initial des températures qui est arbitraire. L'examen de cette condition fait connaître que l'on peut développer en séries convergentes, ou exprimer par des intégrales définies, les fonctions qui ne sont point assujéties à une loi constante, et qui représentent les ordonnées des lignes irrégulières ou discontinues. Cette propriété jette un nouveau jour sur la Théorie des équations aux différences partielles, et étend l'usage des fonctions arbitraires en les soumettant aux procédés ordinaires de l'analyse. » art. 14

---

34 Cela permet à Fourier d'affirmer aussi : « *On connaîtra, par la lecture de cet ouvrage, que la chaleur affecte dans les corps une disposition régulière, indépendante de la distribution primitive, que l'on peut regarder comme arbitraire.* » (art. 19). Sans une représentation d'une distribution initiale arbitraire, ce résultat ne pourrait être obtenu par une analyse mathématique. Cette remarque est déjà dans le mémoire de 1807 : « *Dans la théorie dont nous nous occupons, la forme des intégrales est déterminée par la nature même des conditions physiques, ainsi qu'on le reconnaîtra dans la suite de ce mémoire; toute recherche d'autres intégrales serait ici entièrement infructueuse; mais il était nécessaire de faire coïncider les résultats avec un état initial quelconque.* » Fourier [1807/1972, 253]. Nous retrouverons l'importance du rôle de la condition initiale dans la deuxième partie consacrée à l'analyse de l'introduction des séries trigonométriques lors de la controverse sur les cordes vibrantes. Le rôle de cette condition initiale semble échapper à G. Bachelard [1927, 40] quand il écrit : « *Il est assez remarquable que le phénomène arrive finalement à suivre une loi conditionnée uniquement par la forme géométrique des corps où il se manifeste et qu'il se libère complètement d'une distribution qui est pourtant, au début, la cause déterminante du mouvement calorifique. On assiste là à une géométrisation automatique qui mérite l'attention du philosophe. La forme triomphe en quelque sorte des données physiques disparates ; l'irrationalité est étouffée ; des causes différentes produisent de mêmes effets.* ».

35 Poisson [1835, 167-168] donnera aussi une présentation remarquablement claire de la manière dont on obtient ainsi *diverses* représentations des fonctions : "*La fonction F(x,y,z), tout-à-fait arbitraire, et qui peut être continue ou discontinue, se trouve donc ainsi développée en une série de quantités qui sont toutes de la même forme ; et comme cette forme dépendra de la figure du corps et des conditions relatives, il en résultera une infinité de développemens différens d'une même fonction, propres à en exprimer les valeurs dans une étendue déterminée par celle de chaque corps. Parmi ces développemens, il y en aura qui procéderont suivant les sinus ou cosinus de ces variables multipliées par les racines d'une équation transcendante, d'autres suivant des fonctions transcendantes d'un ordre plus élevé, qui ne pourront même s'exprimer, sous forme finie, que par des intégrales définies. Chaque développement devra être considéré comme une formule d'interpolation, applicable à une étendue limitée des valeurs de la fonction F(x,y,z) ; et quoiqu'on ne puisse pas toujours démontrer directement ces diverses formules, il ne pourra cependant rester aucun doute sur leur exactitude. En effet, la formule (20) représente certainement, d'après les considérations précédentes, la valeur la plus générale de u qui puisse satisfaire à l'équation L=0; par hypothèse, on a déterminé les quantités arbitraires que cette série renferme, au moyen des autres conditions du problème qui a conduit à l'équation L=0, et d'après l'état initial du système: si ces conditions ne sont point incompatibles, et que le problème soit susceptible d'une solution, il faut donc qu'après cette détermination la série (20) exprime la valeur de u à un instant quelconque et en un point quelconque du système; par conséquent, en faisant t=0 dans cette série, elle devra représenter les valeurs initiales de u dans l'étendue où elles auront concouru à la former; ou, autrement dit, l'équation (22) devra être regardée comme une conséquence rigoureuse de la solution complète du problème.*"

# 4 - La diffusion de la chaleur dans la lame et la représentation d'une fonction arbitraire par une série trigonométrique

### *a) Introduction de la représentation*

Au moment d'analyser son introduction, il vaut peut-être la peine de remarquer que la représentation d'une fonction arbitraire par une série trigonométrique est introduite sans avoir vraiment été annoncée. Il n'en n'est pas question dans le « Discours préliminaire » et il faut attendre l'article 14 du chapitre I, dans le passage qui vient d'être cité, pour apprendre qu'il serait possible de « *développer en séries convergentes, ou exprimer par des intégrales définies, les fonctions qui ne sont point assujéties à une loi constante, et qui représentent les ordonnées des lignes irrégulières ou discontinues* ». Encore n'est-il pas même fait mention dans ce passage des séries trigonométriques! Cette représentation n'est pas non plus annoncée au début du chapitre III qui lui est pourtant largement consacré.
Fourier, on l'a aussi vu, justifie de commencer par la lame parce que son étude est la « *plus propre qu'aucune autre à faire connaître les éléments de la méthode que nous avons suivie*» (art. 163). La méthode suivie pour la lame est pourtant différente des suivantes puisque c'est au cours de celle-ci que la représentation d'une fonction par une série trigonométrique est inaugurée. Dans tous les autres cas, Fourier, disposant de cette représentation, pourra considérer d'emblée une distribution initiale arbitraire. Dans le cas de la lame, il commence au contraire par considérer une distribution particulière et ne traite le cas général qu'à la fin. L'introduction de la représentation d'une fonction par une série trigonométrique n'est donc pas simplement insérée au milieu de l'étude de ce cas : elle modifie la formulation du problème et la manière d'arriver au cas général traité. Il s'agit donc maintenant de préciser quel est le rôle de la lame dans l'introduction de cette représentation.

### *b) Le développement de 1 en série trigonométrique et la conformité*

Les séries trigonométriques s'introduisent dans la *Théorie analytique de la chaleur* avec le développement sous cette forme de la fonction constante 1[36]. Chercher un tel développement n'est alors mathématiquement ni évident ni naturel. Il faut une raison pour se lancer dans la recherche d'un développement aussi mathématiquement contre-intuitif et qui n'est obtenu qu'au terme d'un calcul aussi long qu'original (Bernkopf 1968, 313-315).
L'équation que résout Fourier
$$1 = a\cos(y) + b\cos(3y) + c\cos(5y) + d\cos(7y) + \dots$$

n'est en rien un problème arbitraire, elle est l'expression mathématique d'un *fait* physique : par continuité, l'expression de la distribution de chaleur à l'intérieur de la lame doit, pour y=0, coïncider avec la distribution de chaleur à la base fixée à 1. Avant de se lancer dans la détermination de cette série, Fourier est revenu, on l'a vu, à l'expression de la solution générale pour en montrer la signification physique (art. 170)[37]. Avant d'être un problème d'analyse mathématique, cette équation est l'expression d'un fait physique : 1 est développable en série trigonométrique *parce que* d'après la solution générale trouvée pour l'intérieur de la lame, il *doit* en être ainsi. Cette équation n'est pas une question, « 1 est-il développable en série trigonométrique ? », mais une affirmation : « 1 est développable en série trigonométrique ». La recherche de ce développement est commandée par la résolution de cette équation qui est elle-même l'expression d'un fait thermodynamique : si 1 n'était pas développable en série de cosinus, la conformité de l'analyse mathématique avec les phénomènes thermiques serait mise en défaut. La *conformité* de l'analyse est ainsi en jeu dans la *possibilité* de développer 1 en une série de cosinus et trouver ce développement a pu en effet être un moment fondateur pour Fourier [Grattan-Guinness & Ravetz 1972, 147]. Mais la question de conformité n'aurait pas non plus lieu d'être sans la lame...

---

[36] Dès 1805, Fourier étudie la distribution de la chaleur dans la lame, avec ces conditions au bord [Charbonneau 1976, 84, 97].

**c) Fonction constante et fonction arbitraire**

La lame est ainsi une sorte de dispositif qui établit la possibilité de développer 1 en série trigonométrique (sans en déterminer les coefficients). Mais cette possibilité est de cette manière établie aussi bien pour une distribution *quelconque*. De ce point de vue, le choix d'une distribution à la base particulière (constante et égale à 1) apparaît *arbitraire*[38] : le système d'expression (le schéma géométrique) qui représente la lame *ne* permet tout simplement *pas* d'exprimer le choix d'une distribution *particulière*. C'est uniquement en y juxtaposant des expressions d'un autre type (comme « 1 ») que l'on y parvient. Bien entendu, ce système d'expression est tout autant incapable d'exprimer un développement en série trigonométrique! Mais il ne s'agit pas ici de rendre compte de l'introduction de l'expression des séries trigonométriques mais, celle-ci ayant de fait été introduite (par la solution générale), de rendre compte de sa mise en relation avec une fonction *quelconque*. Il s'agit d'identifier les conditions (sémiotiques) de l'introduction de la *généralité* de la représentation par des séries trigonométriques. La généralité a besoin d'une expression, sa nature et son introduction en dépendent, il s'agit dès lors de l'identifier. Et c'est ici un *segment,* le segment de la base de la lame, qui est l'expression de la généralité.

---

37 Une équation similaire apparaissant pour les autres cas de solide, il est ainsi possible d'obtenir d'autres représentations d'une fonction arbitraire. Liouville [1830] exploitera cette possibilité.

38 Dans le manuscrit de 1807, la possibilité de considérer une distribution à la base *quelconque* est considérée *avant* même que ne soit présenté le calcul donnant le développement de 1 ([Fourier 1807/1972, 139 ; Dhombres & Robert 1998, 528]). Le développement de 1 apparaît plus clairement dans le manuscrit comme un cas *particulier* de celui d'une fonction *arbitraire* : « *Il nous reste à déterminer les constantes &c. qui entrent dans l'équation générale*

$$z=a_1 e^{-\frac{1}{2}\pi x}\cos(\frac{1}{2}\pi y)+a_2 e^{-\frac{3}{2}\pi x}\cos(\frac{3}{2}\pi y)+a_3 e^{-\frac{5}{2}\pi x}\cos(\frac{5}{2}\pi y)+a_4 e^{-\frac{7}{2}\pi x}\cos(\frac{7}{2}\pi y)+...etc$$

*On traitera en premier lieu le cas qui se rapporte à la question présente où tous les points qui se rapportent à la première arête ont une température commune.* » (Fourier 1807/1972, 147]

Ainsi, l'expression (géométrique[39]) de la lame, considérée dans son ensemble maintenant, est la représentation d'un solide chauffé qui, en tant que telle, met en relation l'expression générale d'une série trigonométrique et celle d'une fonction particulière/quelconque (par son segment de base). Elle sert de fait de lieu de rencontre de l'expression (analytique) générale d'une série trigonométrique et d'une fonction particulière/quelconque[40].

## 5 - L'inauguration des séries trigonométriques

Les analyses précédentes ont permis de préciser le rapport entre la représentation d'une fonction arbitraire par des séries trigonométriques et l'étude de la distribution de la chaleur dans les différents cas de solide considérés pas Fourier. Nous avons ainsi montré la singularité de la lame et son rôle dans l'*introduction* de la représentation des fonctions par des séries trigonométriques dont Fourier dispose ensuite, ce qui lui permet de considérer après des distributions de chaleur initiales arbitraires. Nous allons à présent considérer l'inauguration proprement dite, c'est-à-dire les arguments avancés par Fourier pour soutenir la possibilité de représenter une fonction arbitraire par une série trigonométrique en les rapportant aux problèmes sémiotiques que pose l'inauguration de cette représentation.

### *a) Les problèmes sémiotiques*

L'inauguration proprement dite a principalement lieu dans la section VI du chapitre III intitulée « *Développement d'une fonction arbitraire en séries trigonométriques* ».
On peut distinguer essentiellement trois arguments :
- une *démonstration* de l'énoncé pour les développements de Taylor impairs ;

- la formule intégrale des coefficients ;

---

39 On ne peut dissocier la manière de poser le problème de la diffusion de la chaleur de la représentation de la lame (comme des autres solides).
40 La question suivante se pose dès lors : un autre objet que la lame peut-il permettre d'introduire la possibilité de représenter n'importe quelle fonction par une série trigonométrique? Nous devons nous restreindre aux solides dont la solution générale prise pour t=0 donne une série trigonométrique (en remplaçant la distribution à la base par la distribution initiale), c'est-à-dire l'armille, le parallélipipède rectangle et le cube. On remarque alors que dans aucun de ces trois cas l'équation du développement de la distribution initiale en série trigonométrique n'est l'expression d'un phénomène physique évident comme dans le cas de la lame : développer cette distribution en série trigonométrique ne résulte pas de la solution générale, c'est à chaque fois, comme on l'a vue, une initiative de Fourier. Autrement dit, l'expression de la généralité est dans ces cas introduite par Fourier, et non par les séries trigonométriques. Les conditions de l'introduction de la représentation d'une fonction arbitraire par les séries trigonométriques ne se retrouvent pas dans ces autres solides. Mais ces autres cas de solide pourraient et ont pu introduire l'idée d'*autres* représentations que celle par les séries trigonométriques.

- des exemples de développements en séries trigonométriques.

C'est avant tout la part de *nécessité* de ces arguments que nous voulons mettre en évidence.

Il est ainsi en premier lieu *inévitable* que Fourier ne puisse donner au mieux, indépendamment des critiques qui peuvent lui être adressées, qu'une démonstration pour une *partie* de la totalité des fonctions. Les séries de Taylor ne représentent qu'une partie de la totalité des courbes auxquelles se rapporte l'énoncé inaugural et elles sont, contrairement aux courbes (géométriques), homogènes aux séries trigonométriques. Cette restriction l'affranchit des problèmes sémiotiques (incommensurabilité) de l'énoncé inaugural. Il lui faudrait sinon disposer d'une représentation alternative de *toutes* les fonctions[41].

Mais il lui faut ensuite soutenir l'énoncé inaugural *au-delà* de ce cas restreint et se confronter à ses problèmes sémiotiques. Cette extension est un moment nécessaire et repérable dans une inauguration. Elle est en l'occurrence réalisée par la formule intégrale des coefficients. Introduite à l'issue de la démonstration dans le cas restreint des séries de Taylor impaires, cette formule est explicitement introduite pour « *étendre les mêmes conséquences* [i.e. la possibilité de développer en série de sinus une fonction impaire développable en série de Taylor] *à des fonctions arbitraires, même à celles qui seraient discontinues et entièrement arbitraires* » (art. 219). Elle est ainsi le principal argument donné par Fourier pour réaliser l'extension requise[42].

Les exemples particuliers de développements en séries trigonométriques répondent aussi à la nécessité de cette extension. La plupart de ces exemples sortent du cadre couvert par la démonstration. C'est leur intérêt : ils établissent que la représentation considérée va au-delà. Et cela *doit* être établi. Il faut convaincre de la validité de l'énoncé inaugural *au-delà* de ce que les représentations disponibles (séries entières) permettent de démontrer. Plus ils le font, mieux ils jouent leur rôle. Particuliers, ils ont toujours aussi un sens général. Le développement de la fonction cosinus en série de sinus n'est pas seulement le développement de la fonction cosinus en série de sinus, c'est la preuve qu'une série de fonctions toutes impaires peut représenter une fonction paire. De même, le développement de la fonction égale à $\frac{\pi}{2}$ entre 0 et α et nulle entre α et π n'est pas tant le développement de cette fonction particulière que la preuve qu'une série de fonctions toutes continues peut représenter une fonction discontinue, etc. Ce sont autant d'objections réfutées. Ces objections étaient générales, si elles sont réfutées, il n'y a plus d'obstacle à ce que toutes les fonctions paires et discontinues soient développables ainsi.

Il est enfin nécessaire de sortir du cadre couvert par une démonstration *par des exemples*. C'est une conséquence de l'absence de représentations alternatives couvrant intégralement la totalité considérée. Autrement, il ne s'agirait pas d'un énoncé inaugural mais d'un théorème de représentation[43].

---

41 L'expression fx qui en est une intervient dans la formule intégrale des coefficients considérée plus loin.
42 Il est à noter que la lame *n*'est *pas* un argument énoncé par Fourier.
43 A comparer avec Dhombres & Robert [1998, 558] : « *L'idée n'est pas encore claire chez Fourier*

La disparité des arguments utilisés pour l'inauguration est le pendant de l'uniformisation introduite par la représentation inaugurée : *les défauts de la « démonstration » de l'énoncé inaugural sont en partie le reflet des vertus de la représentation inaugurée.*

### *b) La formule intégrale du développement en série trigonométrique*

#### *i) Une démonstration défectueuse*

L'extension *de* la formule intégrale démontrée pour une fonction impaire développable en série de Taylor à une fonction arbitraire n'est mathématiquement pas correcte. L'extension *par* la formule intégrale de la possibilité de développer une fonction impaire développable en série de Taylor à une fonction arbitraire ne l'est pas non plus. Ce manque de rigueur a souvent été dénoncé[44]. Fourier s'y expose d'autant plus qu'il qualifie son énoncé de « *théorème* » (art. 362), qu'il fait valoir qu'il « *a donné plusieurs démonstrations de ce théorème dans le cours de [son] Ouvrage* » (art. 503) et que ses « *démonstrations ne laissent aucun doute sur l'exacte vérité de ces équations* » (art. 428-11°)[45]. Les défauts de ses arguments sont incontestables et il est sans doute intéressant de les relever. Il ne faudrait pas pour autant laisser entendre qu'il suffirait de corriger ces défauts pour arriver à une démonstration correcte qui serait aujourd'hui disponible quelque part[46]. Tel n'est pas le cas : cet énoncé *ne peut pas* être démontré. Il importe avant tout de le reconnaître[47] : *compte tenu des moyens d'expression disponibles, cette impossibilité est inhérente à l'énoncé*. Ce sont les conditions sémiotiques de

---

*et il va y revenir à sa façon très particulière qui est l'enveloppement par des exemples. C'est qu'il n'est pas encore sorti de l'induction* ». Ou encore : « *Il le pallie* [le manque de rigueur] *par la multiplication des exemples* » Dhombres & Robert [1998, 561] ; « *Malheureusement, la méthode originale qui guide Fourier dans sa « démonstration » n'a pas été caractérisée par la pratique mathématique, puisque cette méthode est celle de l'exhaustivité des exemples traités, donc au fond du réel inventorié.* » Dhombres & Robert [1998, 565].

44 Cette critique a d'emblée été adressée par les examinateurs du Mémoire de 1807 (qui comprenaient Lagrange, Laplace, Monge et Lacroix) : « *la manière dont l'Auteur parvient à ses équations n'est pas exempte de difficultés, et que son analyse, pour les intégrer, laisse encore quelque chose à désirer, soit relativement à la généralité, soit même du côté de la rigueur* » Oeuvres, I, vii-viii. Voir aussi par exemple : Truesdell [1961, 77] ; Kline [1972/1990, 677] ; Dhombres & Robert [1998, 557-561, 564] ; Grattan-Guinness [2005, 360].

45 On peut encore citer : « Nous avons démontré plusieurs fois le théorème exprimé par l'équation

$$\psi(x) = \frac{2}{\pi} \sum \sin(ix) \int_0^\pi \psi(r) \sin(ir) dr$$

on y peut parvenir de différentes manières, et la formule se déduit très facilement de l'intégration définie; mais ce qu'il importe surtout de reconnaître distinctement, c'est que la série est toujours convergente, et que la valeur attribuée à la variable x doit ici être comprise dans l'intervalle de 0 à π. » Fourier 1829/oeuvres II, 154.

46 Dhombres & Robert [1998, 556] écrivent par exemple : « *Cauchy a pris soin d'une définition précise d'une fonction continue, d'une fonction dérivable, et d'une fonction analytique ou localement développable en série entière. Ces distinctions ne trouvent aucune place dans la* Théorie analytique de la chaleur *dernière mouture, alors qu'elles en permettraient une lecture « rigoureuse ».* » Grattan-Guinness (2005, 360) considère que « *[t]he best solution to the problem of generality would be a proof of their convergence* ».

47 Poisson [1835, 168] pointe, nous semble-t-il, cette impossibilité sans clairement en indiquer les raisons.

possibilité de cette démonstration qui font ici avant tout défaut[48]. Pour démontrer que les séries trigonométriques sont une représentation de la totalité des fonctions, il faudrait une autre représentation de toutes les fonctions susceptible de permettre cette démonstration. L'intérêt de l'énoncé inaugural et l'impossibilité de le démontrer sont indissociables. C'est ce qui justifie de distinguer ces énoncés et ces textes .

### *ii) La formule intégrale des coefficients : l'incommensurabilité éliminée ?*

On a montré que la « démonstration » proposée par Fourier était en partie déterminée par les problèmes sémiotiques inhérents à l'énoncé inaugural. Mais inversement, certains de ses arguments, en l'occurrence la formule intégrale du développement en série trigonométrique, peuvent donner l'impression qu'il a aussi pu s'affranchir de ces problèmes. Cette formule tout à fait remarquable semble en effet donner le moyen de développer en série trigonométrique (à peu près...) n'importe quelle fonction sans rencontrer de problèmes sémiotiques. On ne peut résoudre cette difficulté sémiotique en invoquant simplement les défauts mathématiques des arguments de Fourier ; il faut la résoudre sémiotiquement et montrer que les problèmes sémiotiques qui ne peuvent être résolus ne l'ont pas été.

D'un point de vue sémiotique, la manière dont la formule intégrale permet de transformer l'expression d'une fonction quelconque en série trigonométrique, c'est-à-dire les expressions d'une des représentations de l'énoncé inaugural en l'autre, est remarquable. Une comparaison avec un autre texte inaugural peut aider à mettre en évidence ce qu'il y a là de remarquable.
Descartes soutient dans *La Géométrie* qu'il est possible de mettre en équation n'importe quel problème de géométrie. Fourier soutient quant à lui qu'il est possible de mettre sous forme de série trigonométrique n'importe quelle fonction. Pour en convaincre son lecteur, Descartes décrit comment transformer l'énoncé d'un problème de géométrie en une équation. Sa description de cette transformation est elle-même exprimée d'une certaine manière : c'est en l'occurrence une expression *mixte* qui mêle (inévitablement) l'énoncé d'un problème de géométrie (indéterminé) et une expression polynomiale (indéterminée) qui va lui correspondre[49]. L'expression de cette procédure n'est pas elle-même l'expression d'un problème de géométrie. Ce n'est pas non plus une expression polynomiale. Elle ne relève exclusivement d'aucun des deux types de représentation en présence. Elle en est au contraire un mixte. Il peut même sembler étrange d'envisager qu'il en soit autrement compte tenu de l'incommensurabilité de la représentation des problèmes de géométrie et des polynômes : les caractéristiques d'un énoncé inaugural (incommensurabilité) se répercutent nécessairement sur toute expression les mettant en relation. La formule intégrale du développement trigonométrique est elle-même l'expression de la procédure permettant de passer d'une fonction quelconque à son développement en série trigonométrique. Elle est pourtant très homogène et semble s'affranchir de l'incommensurablité. C'est ce qu'il faut à présent examiner.

---

48 Grattan-Guinness [1970, 11-13 ; 1990, 598-600 ; 2005, 359-360] pointe un « problème de représentabilité » (*representability problem*) mais la périodicité des séries trigonométriques qu'il met en avant ne rend pas compte des problèmes sémiotiques posés par l'énoncé inaugural de Fourier.

49 Il en est de même de l'inauguration des machines logiques de Turing. Il faut pour chaque nombre calculable donner la machine logique qui le calcule. Or, il n'y a pas de machine (logique?!) qui déterminerait à partir d'un nombre donné la machine logique qui le calcule. Il peut éventuellement en exister une pour toute une famille de nombres calculables, mais il n'y en a pas pour tous les nombres calculables, autrement la thèse de Turing ne serait pas un énoncé inaugural mais un théorème de représentation. On retrouve la nécessité de procéder au cas par cas. Frege [1899], Whitehead & Russell [1910] doivent aussi construire pas à pas leurs formules faute d'une expression en mesure de le faire pour toutes les propositions logiques.

Il est utile ici de distinguer d'une part la formule intégrale exprimant l'égalité entre deux expressions analytiques et d'autre part en tant que procédé qui permet de représenter une courbe par une série trigonométrique.
Suivant la première acception, les occurrences de $\phi x$ dans la formule sont associées à des expressions analytiques. Ainsi la démonstration que Fourier donne de cette formule se rapporte à une fonction impaire donnée par son développement en série de Taylor. La formule et sa démonstration n'impliquent aucune courbe. Elle n'établit donc pas qu'une courbe peut être représentée par une série trigonométrique, elle établit seulement qu'une expression analytique (qui ne couvre pas la totalité des fonctions) peut être transformée en une autre expression analytique (une série trigonométrique). Il n'y a là aucune incommensurabilité : l'incommensurabilité n'a pas été résolue, elle a d'emblée été éliminée en considérant une série entière. La formule intégrale est un *théorème de représentation* pour les fonctions impaires développables en série de Taylor. Elle contribue à l'inauguration, comme on s'y attend, en *démontrant* l'énoncé inaugural pour une *partie* de la totalité des fonctions.

Considérons maintenant la formule intégrale comme *procédé* qui associe à une courbe désignée par $\phi x$ une série trigonométrique dont les coefficients $\frac{2}{\pi}\int_0^\pi \phi x \sin.ix\,dx$ dépendent de l'indice $i$ et d'une intégrale faisant intervenir la fonction $\phi x$ et une fonction trigonométrique (sans nous arrêter là non plus aux problèmes mathématiques posés par l'affirmation de l'égalité entre les deux). Ce procédé s'applique à une fonction arbitraire dans la mesure où les intégrales qui y figurent sont définies pour une fonction *arbitraire*. Or ces intégrales désignent pour Fourier des mesures d'aires[50] : il n'y pas d'intégrale s'il n'y a pas d'aire. L'intégrale n'a en l'occurrence pas d'autre sens. Elles se rapportent à des *courbes*. C'est en considérant une *courbe quelconque* que Fourier peut inférer l'existence, pour une fonction quelconque, des intégrales qui figurent dans la formule. Ou plutôt, Fourier ne peut considérer une *fonction* quelconque qu'à partir d'une *courbe* quelconque : au moment d'inaugurer les séries trigonométriques il n'a pas d'autre expression d'une fonction *quelconque*. La possibilité de considérer que ces intégrales sont définies pour une *fonction arbitraire* fait intervenir la possibilité de se référer à la *totalité des courbes* ou à n'importe quelle courbe. Les fonctions et de courbes sont l'objet d'autres associations dans ce texte, mais celle-ci serait sans doute l'une des plus difficiles à éliminer. L'histoire de la théorie de l'intégration après Fourier, et dès Cauchy, en donne une bonne indication.

---

50 «*Les coefficients des séries trigonométriques sont des aires définies, et ceux des séries de puissances sont des fonctions données par la différentiation, et dans lesquelles on attribue aussi à la variable une valeur définie.*» art. 235. Voir aussi ci-dessus, la citation de l'art. 220.

La formule intégrale implique donc à la fois les expressions analytiques des fonctions *et* leurs expressions géométriques (courbes). Il est impossible de rendre complètement compte de l'usage de l'expression $\phi x$ au sein de la formule intégrale sans lui reconnaître ces *deux* contenus[51]. Cette formule est souvent utilisée par Fourier comme une relation fonctionnelle pour laquelle il *n*'est pas *nécessaire* d'associer un contenu géométrique aux expressions qui y figurent. Mais si l'on veut rendre compte *complètement* de la formule, et notamment l'appliquer à une fonction *quelconque*, il est ici nécessaire de se référer aux courbes, et ainsi de reconnaître aux fonctions un double contenu, analytique et géométrique[52].

La dualité et avec elle l'incommensurabilité n'ont donc pas disparu de la formule intégrale, elles ont été enfouies dans le signe de fonction[53].

### *iii) La généralité de la formule intégrale*

La formule intégrale[54] est le principal argument donné par Fourier pour étendre à des fonctions arbitraires la possibilité d'une représentation par des séries trigonométriques démontrée pour une partie d'entre elles. Nous venons aussi de voir qu'elle n'éliminait pas les problèmes sémiotiques que pose une telle extension à des fonctions arbitraires. Elle peut jouer ce rôle en raison de caractéristiques remarquables[55]. La généralité de cette formule, ce qui en fait une « formule », vient, d'un point de vue sémiotique, de l'expression $\varphi x$. Or, cette expression a des caractéristiques sémiotiques tout à fait remarquables : elle exprime la *totalité des fonctions* par une expression qui peut elle-même être considérée comme *une fonction* : on peut l'additionner, la multiplier, la dériver, l'intégrer, etc. A *contrario*, l'expression « fonction » est aussi un moyen d'exprimer la totalité des fonctions mais elle n'est pas elle-même considérée comme une fonction : on ne peut l'additionner, la multiplier, la dériver, l'intégrer, etc. L'expression « fonction » fait toujours *référence* à une fonction, on ne peut la prendre elle-même pour la fonction. Avec l'expression $\varphi x$ cela est possible. Disposer d'une telle expression c'est disposer d'un moyen d'expression remarquable d'une totalité : elle peut être utilisée comme l'expression d'une fonction *particulière,* comme par exemple *cos(x)*, mais sans que les expressions dans lesquelles elle entre soient ensuite entachées de ce caractère particulier[56]. Elle permet de faire des calculs/démonstrations comme s'il s'agissait d'une fonction particulière mais en donnant la preuve *de fait* que l'on n'a pas fait intervenir les caractéristiques particulières de cette fonction. Elle permet ainsi d'exprimer la *généralité* avec une expression qui a les caractéristiques d'une fonction *particulière*. L'induction est ici légitime.

---

51 Pour la description d'autres signes en mathématiques faisant intervenir plusieurs plans de contenu voir Herreman [2000].
52 « Ce théorème et le précédent conviennent à toutes les fonctions possibles, soit que l'on en puisse exprimer la nature par les moyens connus de l'Analyse, soit qu'elles correspondent à des courbes tracées arbitrairement. » art. 224.
53 Le contenu géométrique des intégrales de cette formule a déjà été souligné par Hawkins [1970, 7-8], Grattan-Guinness [1970, chp. 2], Grattan-Guinness & Ravetz [1972, 191-192], Charbonneau [1976, 117-118], Dhombres & Robert [1998, 561]. Nous avons voulu montrer ici que la présence de ce contenu géométrique répondait à une *nécessité* découlant des caractéristiques sémiotiques de l'énoncé inaugural.
54 Sur cette formule, voir Grattan-Guinness & Ravetz [1972, 191], Grattan-Guinness [1990, 597-600], Dhombres & Robert [1998, 552-561].
55 Dhombres & Robert [1998, 555] : « *l'arbitraire est le fruit de la formule, non celui de la démonstration* ».

Ces caractéristiques ne son évidemment pas *dans* l'expression $\varphi x$ : elles découlent de la possibilité de substitutions (substitution des fonctions particulières les unes aux autres, substitution de fonctions particulières à $\varphi x$). Elles supposent donc d'autres expressions *dans* lesquelles se font ces substitutions. Il s'agit donc de caractéristiques d'une expression dans un *système d'expressions* (d'où la nécessité d'abandonner une sémiotique du signe)[57].

### iv) La formule intégrale et la lame

D'un point de vue sémiotique, la formule intégrale et la lame rendent finalement compte l'une et l'autre de la possibilité de représenter une fonction arbitraire par une série trigonométrique[58]. Seule la formule intégrale participe de l'inauguration. Inversement la lame rend sans doute mieux compte de la manière dont Fourier a été amené à envisager cette possibilité et à s'en convaincre. Leur rôle est bien distinct dans le texte, mais elles ont des caractéristiques sémiotiques communes : l'une et l'autre comprennent une expression de la totalité des fonctions qu'elles mettent en rapport avec l'expression d'un développement en série trigonométrique. Pour l'une c'est l'expression d'un segment, pour l'autre l'expression $\varphi x$.

Le *double* rôle de la lame a aussi été mis en évidence. La lame est d'abord la représentation d'un solide dont Fourier étudie la distribution de chaleur. Elle est/représente alors un de ces objets ou une de ces représentations tout à fait typiques de la physique-mathématique qui circulent entre les protagonistes et dont l'examen structure l'exposé de la théorie. Elle est, à cet égard, tout à fait semblable aux autres solides étudiés. Mais elle est aussi un dispositif qui, comme la formule intégrale, met en relation les séries trigonométriques et une distribution arbitraire de chaleur. En cela elle diffère des autres solides et a un rôle singulier[59]. Ce double rôle en fait une expression fondamentalement ambivalente, à la fois objet et dispositif, sémiotique et méta-sémiotique. Cette ambivalence est la condition du rapprochement des séries trigonométriques et d'une fonction arbitraire.

---

56 Ainsi, si la formule intégrale à des limites, elles viennent des intégrations qui y figurent (ou encore des produits et des additions).

57 Les caractéristiques indiquées définissent ce que l'on peut appeler une *expression d'invariante génératrice* [Herreman 2005, 169-172]. Disposer d'expression d'invariante génératrice suppose, on le voit, des moyens d'expression élaborés dont Fourier, de toute évidence, dispose. Une histoire de l'élaboration de ces conditions sort évidemment du cadre de cet article.

58 Dans Dhombres & Robert [1998, 554] l'introduction de la possibilité de la représentation d'une fonction arbitraire par une série trigonométrique semble n'être rapportée qu'à la formule intégrale : « *la relation porte un sens que Fourier donne aussitôt, et ce sens est celui de la généralité. La généralité porte sur la fonction. Nous y sommes !* ». Voir aussi Dhombres & Robert [1998, 555]. Cette formule fait bien en effet « *connaître comment les fonctions entièrement arbitraires peuvent aussi être développées en séris de sinus d'arcs multiples* » Fourier [1807/1972, 261]. Connaître ce « *comment* » est essentiel, mais la lame a déjà établi la possibilité et la généralité de ces développements.

59 On pourrait faire de la lame un instrument thermique servant à développer une fonction en série trigonométrique comme il en existe pour calculer leur intégrale.

Une ambivalence semblable apparaissait dans l'inauguration de la représentation des problèmes et des courbes par des polynômes par Descartes. Comme la lame est la représentation typique d'un solide, le problème de Pappus est un problème typique de géométrie. Et comme la lame, le problème de Pappus est aussi dans *La Géométrie* un dispositif particulier qui fait se rencontrer problèmes de géométrie et polynômes. Le problème de Pappus est un *problème mathématique* qui permet de parcourir *l'ensemble des problèmes de géométrie* (du point de vue de leur représentation par des polynômes, c'est faux, mais Descartes le croit). Et il en est encore de même avec l'inauguration de la représentation des courbes par des polynômes : le mésolabe est à la fois un instrument parmi d'autres (répondant à la caractérisation retenue par Descartes pour définir une courbe géométrique) et un instrument qui trace des courbes géométriques de *tous* les genres.

Il y a bien sûr aussi des différences. La lame n'intervient pas ici comme *un solide* qui serait aussi le moyen d'expression de *tous les solides*[60] : elle intervient en tant qu'expression d'*un solide* qui contient un segment de base qui sert d'expression de *toutes les distributions de chaleur* (courbes). Mais son ambivalence est essentielle à l'exposé de la théorie : elle s'inscrit dans la série d'objets étudiés (dont elle est même le premier terme) et elle contient le segment qui sert d'expression de toutes les fonctions.

Ces similitudes ne sauraient être fortuites. Les attribuer à des influences ne convient évidemment pas (sauf à considérer des « influences sémiotiques »...). On peut à nouveau se tourner vers les conditions de possibilité, ni nécessaires ni contingentes, de l'émergence de ces représentations générales.

La lame, la formule intégrale, le problème de Pappus et le mésolabe ont en commun d'être à la fois *une* expression et le lieu d'une variabilité. La lame est l'expression d'un solide et le lieu de la variabilité de la distribution de la chaleur à sa base. Dans la formule intégrale fx, on l'a vu, est à la fois l'expression d'une fonction et le lieu d'une variabilité des fonctions. Le problème de Pappus est l'expression d'un problème et le lieu de la variabilité par le nombre et la configuration de droites qui entrent dans son énoncé. Le mésolabe est un instrument et une chaîne indéfiniment prolongeable d'instruments. On voit ainsi à chaque fois le double rôle de ces expressions et la fonction de celui-ci dans l'émergence de nouvelles représentations : chacune de ces expressions intervient en tant qu'expression (d'un solide, d'un problème, d'un instrument), mais aussi en tant que lieu d'une variabilité associée à une totalité. La double fonction sémiotique de ces expressions s'explique ainsi par le fait de comprendre une expression qui sert de support de l'expression d'*une totalité* ainsi mise en relation avec une autre expression faisant *partie* de cette totalité[61].

---

60 Le calcul intégral et différentiel permet de passer de l'étude de la propagation de la chaleur dans un solide (infinitésimal) de forme géométrique *particulière* (en l'occurrence un parallélipipède, dont la lame est un cas particulier) à celle dans un corps de forme *quelconque*. C'est là que réside la possibilité d'établir une équation *générale*. Cette expression de la généralité propre au calcul intégral et différentiel, bien sûr très importante en mathématiques, joue en tant que tel comme on l'a vu peu de rôle dans la théorie de la chaleur de Fourier.

61 On peut trouver là un élément qui rend compte sémiotiquement de l'impact sur les mathématiques

Il convient aussi de souligner la diversité des types d'expression qui ont joué un rôle dans l'émergence de ces représentations puisque l'on trouve aussi bien des schémas (de la lame, d'un instrument), une formule mathématique (la formule intégrale) et une expression en langue naturelle faisant référence à une configuration géométrique (l'énoncé du problème de Pappus).

### *v) Intervalle de définition et de représentation*

Fourier a été amené à considérer l'intervalle de validité de la représentation d'une fonction par une série trigonométrique et à introduire à la même occasion la notation $\int_a^b$ , dans laquelle les bornes de la variable sont indiquées en haut et en bas du symbole d'intégration (arts 222, 231 et 431, Fourier 1829/oeuvres II, 150, Fourier 1816, 361)[62], notation qui a été très vite adoptée et diffusée par Cauchy [1814/1827/1882, 340 ; 1823a/1899, 126 ; 1823b/1905, 334].

Fourier consacre en effet d'assez longs développements à déterminer la « *limitation à laquelle les séries trigonométriques sont assujéties* » (art. 184, arts 185-188), c'est-à-dire à déterminer l'intervalle sur lequel une série trigonométrique représente une fonction donnée par une autre expression. De toute évidence, la série trigonométrique ne peut représenter partout où elle est elle-même définie une fonction donnée (ne serait-ce qu'en raison de sa périodicité)[63]. La représentation trigonométrique doit nécessairement avoir des limites et il est donc nécessaire de les déterminer pour soutenir la généralité de cette représentation. D'où les développements que Fourier consacre à cette question.

Considérons à présent plus précisément ce qui le conduit à s'intéresser à cet intervalle. D'un point de vue sémiotique, une particularité remarquable de la représentation par les séries trigonométriques est sa dépendance à l'intervalle : si l'on change l'intervalle, ce n'est pas seulement l'intervalle mais le terme général (la base) de la série qui change. C'est là une différence remarquable avec un développement en série entière : même quand il n'est que local, les coefficients de la série peuvent changer, mais le terme général du développement, c'est-à-dire la *base de la représentation*, ne change pas. Dans le cas des séries trigonométriques, la base même de la représentation change. La possibilité de représenter une fonction par une « série trigonométrique » est indépendante de tout intervalle, mais le terme général de la représentation en dépend. Il est tout à fait original et remarquable d'avoir une représentation dont la base soit ainsi « variable »[64]. Ainsi, l'intervalle n'est pas tant associé aux fonctions représentées qu'au *système de représentation. Ce n'est pas la définition de la fonction qui dépend de l'intervalle, c'est le système d'expression.* L'intervalle est un intervalle de représentation et non

---

de leur application à la physique.
62 Cela a bien sûr été souvent souligné, voir notamment : Grattan-Guinness & Ravetz [1972, 241] ; Dhombres & Robert [1998, 552, 554 , 559-561] ; Grattan-Guinness [1990, 629] ; Grattan-Guinness [2005, 359].
63 Fourier joue ici d'une *non* conformité : l'intervalle sur lequel la série trigonométrique définit une fonction *ne* coïncide *pas* avec l'intervalle sur lequel cette série sert de représentation.
64 Les exposés actuels privilégient l'intervalle [0,1] et tendent ainsi à masquer cette variabilité. La variabilité de la base dans la représentation des nombres (base décimale, binaire, hexadécimale, etc.) est aussi remarquable mais la totalité (des nombres) représentée reste dans ce cas la même.

un intervalle de définition[65].

La représentation trigonométrique présente donc cette caractéristique sémiotique tout à fait originale de n'être pas seulement un système d'expression mais une *variété* de systèmes d'expression[66]. Fourier a à résoudre divers problèmes de mathématiques, de physique et de physique-mathématiques, mais il a aussi à résoudre les problèmes qui découlent des caractéristiques de la représentation qu'il inaugure, problèmes qui comprennent ceux communs à toutes les inaugurations mais ne s'y réduisent pas (et les séries trigonométriques en détermineront encore bien d'autres après lui, notamment celui, typiquement sémiotique, de l'unicité de la représentation) : les caractéristiques sémiotiques inhabituelles de la représentation conduisent à des considérations mathématiques inhabituelles.

On peut aussi à présent considérer d'un point de vue sémiotique la concomitance de la prise en compte de l'intervalle de représentation par des séries trigonométriques et l'introduction de la notation $\int_a^b$. L'intervalle sur lequel une série trigonométrique représente une fonction donnée ne coïncide pas, on l'a vu, avec celui sur lequel la série définit elle-même une fonction. Cette limitation de la représentation d'une fonction par une série trigonométrique *doit apparaître* dans la formule intégrale qui serait sinon manifestement fausse. L'introduction de la notation $\int_a^b$, qui en tant que telle n'a plus de rapport aux séries trigonométriques, répond à la nécessité d'exprimer l'intervalle de validité de la formule intégrale et celle de la représentation de la fonction par la série

---

65 Dhombres & Robert [1998, 556] y voient un intervalle de définition : « *[l]a notation des bornes d'une fonction dans une intégrale – un coefficient de Fourier – correspond exactement à la précision du domaine de définition qui doit accompagner toute fonction, et la réduire à une commodité c'est passer à côté de ce qu'apporte Fourier* » . La prise en compte de cet intervalle est décrite en termes de « *prise de conscience* » et de compréhension (« *il a compris* ») : « *Ce qui nous apparaît aujourd'hui banal nécessita une prise de conscience pour laquelle Fourier créa sinon un concept nouveau, du moins ou mieux, un paradigme autre. Par cette distinction, il opère une coupure importante dans l'analyse mathématique, quelque chose qui lui fait tourner le dos à la tradition incarnée par Euler (...). S'il a compris que, contrairement à l'égalité de deux formes, celle de deux fonctions ne signifie pas une coïncidence partout, mais peut ne valoir que par un domaine dont il importe de donner l'étendue, c'est bien grâce aux séries de Fourier* » Dhombres & Robert [1998, 559]. Le point de vue épistémologique adopté par ces auteurs les conduit à rapporter au *concept* de fonction ce qui, selon nous, doit l'être en l'occurrence à leur représentation, et plus précisément à une caractéristique bien particulière de celle-ci. De même, les problèmes sémiotiques inhérents à l'inauguration des séries trigonométriques, et notamment celui de l'expression d'une totalité, est ramené par eux à un défaut conceptuel : « *Il « démontre » quant à lui la totalité en faisant des recensions complètes d'exemples. Faute de la notion conceptuelle* ad hoc *qui assurerait cette totalité!* » Dhombres & Robert 1998, 565. Analyser la *Théorie analytique de la chaleur* en considérant qu'un concept (« *la notion conceptuelle* ad hoc ») « *assurerait cette totalité* » c'est à notre avis manquer ce qu'un concept, en l'occurrence celui de fonction, doit aux caractéristiques des moyens qui servent à son expression. Leur analyse est aussi pour une grande part fondée sur le concept de « rigueur » (Dhombres & Robert 1998, 564, 556-557, 657-662] sans que ne soient pris en compte les conditions de possibilité sémiotique de cette rigueur. L'épistémologie et la philosophie des mathématiques sont ainsi souvent des auxiliaires efficaces pour entretenir l'ignorance des conditions de possibilité sémiotiques des mathématiques et avec elle celle de la nécessité de les inaugurer.

66 Reconnaître par exemple que deux fonctions coïncident (i.e. ont les mêmes valeurs) sur un intervalle à partir de leur développement en séries trigonométriques sur des intervalles différents, sans être disjoints, est un problème difficile parce qu'il s'agit en fait d'expressions de deux systèmes d'expression différents.

trignométrique.

Une notation doit être générale et Fourier introduit ici une notation de l'intégrale qui spécifie l'intervalle d'intégration. En dépit de sa généralité contestable et contestée, en raison d'une extension indéterminée, c'est *en tant que formule générale* que la formule intégrale intervient dans l'introduction de cette notation toujours en vigueur. Cette nouvelle notation, compte tenu de la généralité qui lui est nécessairement associée, ne renvoie pas à une fonction particulière mais à une *totalité* de fonctions, la totalité des fonctions définies (et intégrables) sur l'intervalle. Ainsi, la *totalité* des fonctions définies (et intégrables) sur un intervalle a bien été inaugurée et intervient à nouveau ici en tant que telle. La notation $\int_a^b$ en est à la fois une conséquence et une manifestation. En s'arrêtant à ce que la généralité de la représentation trigonométrique et la formule intégrale ont de contestable on risque de manquer certaines conséquences de leur généralité comme ici son rôle dans l'introduction de la notation $\int_a^b$ : les moyens d'expression utilisés en mathématiques, et notamment les moyens d'expression de la généralité, ont des conséquences effectives sur les mathématiques dont les conséquences mathématiques reconnues, et *a fortiori* les conséquences logiques, ne rendent pas compte.

*c) La conformité*

L'inauguration ne consiste pas seulement à soutenir une correspondance réciproque entre les courbes et les séries trigonométriques : toutes les propriétés des courbes doivent pouvoir être exprimées par celles des séries trigonométriques et toutes les propriétés des séries trigonométriques doivent inversement avoir une interprétation géométrique. Soutenir cette conformité au-delà de là correspondance pose aussi des problèmes spécifiques.

*i) L'inauguration des séries et la conformité de l'Analyse mathématique*

La conformité est sans aucun doute la caractéristique la plus remarquable des représentations introduites par les textes inauguraux. C'est elle qui justifie et permet de distinguer ces représentations et les textes qui les inaugurent. Prétendre à une telle conformité pose un problème d'expression évident : comment exprimer la *totalité* des propriétés que ce soit des courbes ou des séries trigonométriques? Ce problème se pose et c'est bien un problème d'expression semblable, mais distinct, de celui de l'expression de la *totalité* des fonctions. Il convient donc aussi d'examiner comment Fourier le résout.

Il faut pour commencer remarquer que Fourier *dispose* de la conformité des courbes et des polynômes déjà soutenue par Descartes, étendue aux séries entières et au-delà, disons aux fonctions algébriques, et surtout amplement soutenue par une longue pratique mathématique. Il n'a ainsi plus à convaincre son lecteur, ni lui-même, comme Descartes avait dû le faire, que la tangente et la normale à une courbe peuvent être obtenues à partir de l'expression analytique de la courbe (Herreman 2012), qu'il en est de même de l'aire, etc. Plus de cent cinquante ans de mathématiques ont permis d'établir cette conformité et sans doute même déjà de faire oublier la nécessité de l'établir. C'est elle notamment qui permet à Fourier de considérer, ce qui n'est autrement rien moins qu'évident, que la solution déterminée par l'Algèbre ou l'Analyse d'un système d'équations (différentielles ou non) établies à partir d'un problème de géométrie, de physique, etc. donne la réponse géométrique, physique, etc. de ce problème[67]. Descartes avait dû consacrer un livre entier à l'établir pour la géométrie...

Mais la conformité ne se réduit pas ici aux courbes et aux séries trigonométriques. Tout au long de son exposé Fourier s'assure systématiquement de la complète interprétabilité physique des expressions mathématiques qu'il considère. Il s'agit pour lui de soutenir que l'Analyse mathématique « *réduit* toutes [c'est nous qui soulignons] *les recherches physiques sur la propagation de la chaleur à des questions de Calcul intégral dont les éléments sont donnés par l'expérience.* » (art. 1). Il fait ainsi valoir que ses expressions mathématiques donnent des valeurs de la température en un point de la lame conformes aux expériences :

> « Il restait encore à comparer les faits avec la Théorie. On a entrepris, dans cette vue, des expériences variées et précises dont les résultats sont conformes à ceux du calculs et lui donnent une autorité qu'on eût été porté à lui refuser dans une matière nouvelle et qui paraît sujette à tant d'incertitude. » art. 15

On a vu qu'avant de déterminer les coefficients de la solution générale obtenue pour la lame, il s'attachait à montrer « *l'effet que représente chacun des termes de la série* » (art. 170). Et une fois la résolution terminée, il conclut en soutenant que « *toutes les conditions physiques de la question sont exactement remplies* » (art. 190). L'expression des séries offre selon lui une représentation de la chaleur qui se diffuse, de son mode fondamental, de ses modes propres etc :

> « Les intégrales que nous avons obtenues ne sont point seulement des expressions générales qui satisfont aux équations différentielles ; elles représentent de la manière la plus distincte l'effet naturel qui est l'objet de la question. C'est cette condition principale que nous avons eu toujours en vue, et sans laquelle les résultats du calcul ne nous paraîtraient que des transformations inutiles. Lorsque cette condition est remplie, l'intégrale est, à proprement parler, l'équation du phénomène ; elle en exprime clairement le caractère et le progrès, de même que l'équation finie d'une ligne ou d'une surface courbe fait connaître toutes les propriétés de ces figures. Pour découvrir ces solutions, nous ne considérons point une seule forme de l'intégrale ; nous cherchons à obtenir immédiatement celle qui est propre à la question. C'est ainsi que l'intégrale, qui exprime le mouvement de la chaleur dans une sphère d'un rayon donné, est très-différente de celle qui exprime ce mouvement dans un corps cylindrique, ou même dans une sphère d'un rayon supposé infini. Or, chacune de ces intégrales a une forme déterminée qui ne peut pas être suppléée par une autre. Il est nécessaire d'en faire usage, si l'on veut connaître la distribution de la chaleur dans le corps dont il s'agit. En général, on ne pourrait apporter aucun changement dans la forme de nos solutions, sans leur faire perdre leur caractère essentiel, qui est de représenter les phénomènes.» art. 428

---

67 Cette conformité sous-tend le premier principe énoncé par Fourier pour décrire sa méthode : « *On considère à-la-fois la condition générale donnée par l'équation aux différences partielles, et toutes les conditions singulières qui déterminent entièrement la question, et l'on se propose de former l'expression analytique qui satisfait à toutes ces conditions.* » art. 428-1°.

Ainsi, l'expression $\quad v = a e^{-x} \cos. y + b e^{-3x} \cos. 3y + c e^{-5x} \cos. 5y + d e^{-7x} \cos. 7y + ...$
est *intégralement* signifiante et son analyse est *conforme* à l'analyse des phénomènes thermiques[68]. Elle rend en quelque sorte la propagation de la chaleur visible à nos yeux comme les résonances harmoniques d'un corps sonore sont audibles par nos oreilles (xxiv). Cette expression ne « *laisse rien de vague et d'indéterminé* » (xxii). On peut avec elle « *déterminer toutes* [c'est moi qui souligne] *les circonstances du mouvement permanent de la chaleur dans une lame rectangulaire échauffée à son origine* » (art. 192), elle « *représente exactement* toutes [c'est moi qui souligne] *les circonstances du mouvement de la*

---

[68] Citons : « *L'irradiation de la chaleur a une relation manifeste avec les Tables de sinus ; car les rayons, qui sortent d'un même point d'une surface chauffée, diffèrent beaucoup entre eux, et leur intensité est rigoureusement proportionnelle au sinus de l'angle que fait leur direction avec l'élément de la surface. Si l'on pouvait observer pour chaque instant, et en chaque point d'une masse solide homogène, les changements de température, on retrouverait dans la série de ces observations les propriétés des séries récurrentes, celles des sinus et des logarithmes ; on les remarquerait, par exemple, dans les variations diurnes ou annuelles des températures des différents points du globe terrestre qui sont voisins de la surface.* » art. 20. Ou encore : « *Cette forme du résultat est nécessaire, parce que le mouvement variable, qui est l'objet de la question, se compose de tous ceux qui auraient lieu séparément, si chaque point du solide était seul échauffé, et que la température initiale de tous les autres fût nulle.* » art. 428-5°. Et enfin : « *chacune de ces solutions donne l'équation propre du phénomène, parce qu'elle le représente distinctement dans toute l'étendue de son cours, et qu'elle sert à déterminer facilement en nombre tous les résultats* » art. 428-7°

Ou encore : « L'équation

$$v = \frac{4}{\pi} e^{-x} \cos(y)$$

représente ainsi un état du solide qui se conserverait sans aucun changement, s'il était d'abord formé ; il en serait de même de l'état exprimé par l'équation

$$v = \frac{4}{3\pi} e^{-3x} \cos(3y)$$

et en général, chaque terme de la série correspond à un état particulier qui jouit de la même propriété. Tous ces systèmes partiels existent à la fois dans celui que représente l'équation (α) ; ils se superposent, et le mouvement de la chaleur a lieu pour chacun d'eux de la même manière que s'il était seul. Dans l'état qui répond à l'un quelconque de ces termes, les températures fixes des points de la base A diffèrent d'un point à un autre, et c'est la seule condition de la question qui ne soit pas remplie ; mais l'état général qui résulte de la somme de tous les termes satisfait à cette même condition.

(...)

On voit donc que les valeurs particulières

$$a e^{-x} \cos(y), b e^{-3x} \cos(3y), c e^{-5x} \cos(5y), ...$$

prennent leur origine dans la question physique elle-même et ont une relation nécessaire avec les phénomènes de la chaleur. Chacune d'elle exprime un mode simple suivant lequel la chaleur s'établit et se propage dans une lame rectangulaire, dont les côtés infinis conservent une *chaleur* » (art. 206)[69]. La conformité ne saurait être affirmée de manière plus explicite. Le chapitre III se conclut sur une phrase dans laquelle la conformité de la solution trouvée est à nouveau affirmée :

> « En général il n'y a aucune propriété du mouvement uniforme de la chaleur dans une lame rectangulaire, qui ne soit exactement représentée par cette solution. » art. 238.

Il ne fait pas de doute que Fourier soutient la conformité des séries trigonométriques aussi bien aux courbes qu'aux phénomènes thermiques auxquelles elles servent d'expression. Mais il ne la soutient pas que pour les séries trigonométriques. Au contraire, nous avons vu que c'était aussi par souci de conformité qu'il partait des équations *particulières* plutôt que de *l'équation générale.* C'est toute la théorie de la chaleur qui doit offrir une expression conforme des phénomènes thermiques, et c'est, au-delà, toute l'Analyse mathématique qui est conforme aux phénomènes sensibles[70] :

> « On reconnaitrait encore les mêmes résultats et tous les éléments principaux de l'Analyse générale dans les vibrations des milieux élastiques, dans les propriétés des lignes ou des surfaces courbes, dans les mouvements des astres et dans ceux de la lumière ou des fluides. C'est ainsi que les fonctions obtenues par des différentiations successives, et qui servent au développement des séries infinies et à la résolution numérique des équations, correspondent aussi à des propriétés physiques. La première de ces fonctions, ou la fluxion proprement dite, exprime, dans la Géométrie, l'inclinaison de la tangente des lignes courbes, et, dans la Dynamique, la vitesse du mobile pendant le mouvement varié : elle mesure, dans la Théorie de la chaleur, la quantité qui s'écoule en chaque point d'un corps à travers une surface

---

> température constante. Le système général des températures se compose toujours d'une multitude de systèmes simples, et l'expression de leur somme n'a d'arbitraire que les coefficients a, b, c, d, ... » art. 191

Citons enfin : « La décomposition dont il s'agit n'est point un résultat purement rationnel et analytique ; elle a lieu effectivement et résulte des propriétés physiques de la chaleur. » Fourier, "Extrait du mémoire sur la chaleur", 1807, Hérivel 1980, 55

69 A propos des termes qui composent les séries : « *Ainsi les valeurs particulières que nous avons considérées précédemment, et dont nous composons la valeur générale, tirent leur origine de la question elle-même. Chacune d'elles représente un état élémentaire qui peut subsister de lui-même dès qu'on le suppose formé ; ces valeurs ont une relation naturelle et nécessaire avec les propriétés physiques de la chaleur.* » art. 241.

70 A ce propos, voir Dhombres & Robert [1998]. On peut citer : « *[Fourier] est persuadé que son analyse se confond avec la nature* » Dhombres & Robert [1998, 731]. Voir aussi Dhombres & Robert [1998, 475 ; 478 ; 502 ; 519 ; 532 ; 566 ; 594 ; 665]. Friedman [1977, 83] : « *Fourier's conception of the science contains more than the mathematical equations for heat propagation; for him the physical aspects remain, even when the equations can stand alone.* ». Grattan-Guinness & Ravetz [1972, 307] : « *We see again in this section [Fourier 1807, art. 114] Fourier's strong desire to develop his results as statements about heat diffusion as well as mathematical theorems. One of the most impressive features of the manuscripts is the constants interplay between new mathematical results and their interpretation in the physical problem at hand; it is this quality which gives it its position as a landmark in pure and applied mathematics alike.* »

> donnée. L'Analyse mathématique a donc des rapports nécessaires avec les phénomènes sensibles ; son objet n'est point créé par l'intelligence de l'homme ; il est un élément préexistant de l'ordre universel et n'a rien de contingent et de fortuit ; il est emprunt dans toute la nature. » art. 20

Le système d'expression de l'Analyse mathématique est selon lui tel que d'une part toutes les questions relatives à la propagation uniforme de la chaleur peuvent y être exprimées et d'autre part que tous les problèmes auxquels ces questions sont réduites peuvent y être complètement résolues, leurs solutions donnant les réponses aux questions posées. Les phénomènes thermiques peuvent ainsi être entièrement saisis par des équations différentielles et décrits par leurs solutions obtenues par une analyse mathématique qui ne doit plus rien aux phénomènes considérés, menée suivant ses propres règles, mais qui néanmoins restitue tous les aspects des phénomènes. La conformité de l'analyse mathématique et les phénomènes thermiques va de pair avec leur séparation (dualité) et l'existence propre et autonome de chacune :

> « Les questions relatives à la propagation uniforme ou au mouvement varié de la chaleur dans l'intérieur des solides sont réduites, par ce qui précède, à des problèmes d'Analyse pure, et les progrès de cette partie de la Physique dépendront désormais de ceux que fera la science du calcul. Les équations différentielles que nous avons démontrées contiennent les résultats principaux de la théorie ; elles expriment, de la manière la plus générale et la plus concise, les rapports nécessaires de l'analyse numérique avec une classe très étendue de phénomènes, et réunissent pour toujours aux sciences mathématiques une des branches les plus importantes de la Philosophie naturelle. Il nous reste maintenant à découvrir l'usage que l'on doit faire de ces équations pour en déduire des solutions complètes et d'une application facile. » art. 163

La conformité suppose la dualité. Pour que l'Analyse mathématique soit conforme aux phénomènes physiques elle *ne* doit *pas* en faire partie ; ses expressions ne peuvent se substituer à celles qui composent les phénomènes physiques.

On peut mentionner ici la démonstration de l'unicité de la solution trouvée de l'équation différentielle pour la lame infinie (l'argument est repris au chapitre VII pour le prisme rectangulaire, arts 326-327) (arts 200-204)[71]. Fourier consacre un développement assez long (cinq articles) à cette unicité qu'il aurait pu inférer de l'unicité phénoménale (la lame n'hésite pas entre plusieurs états permanents...). En montrant que l'un comme l'autre ont une seule solution, il contribue à établir, au lieu de la présupposer, la conformité du problème mathématique (équation différentielle avec conditions aux bords) à la question physique ainsi que la capacité de l'Analyse mathématique à résoudre elle-même les problèmes que l'on peut y formuler. Pour démontrer cette unicité il introduit un théorème sur l'additivité des lois d'évolution de la distribution de la chaleur dans deux lames mais doit, pour le démontrer, faire une analyse moléculaire des échanges thermiques et ainsi revenir aux phénomènes physiques au cours de l'analyse mathématique. C'est là une des rares entorses dans ce chapitre à la séparation des deux analyses. Fourier ne réussit pas ici à reproduire exactement *dans le cours de son exposé* la séparation entre les phénomènes et leur représentation mathématique. Mais cette entorse intervient dans la démonstration d'un problème qui n'est posé que pour mieux établir cette séparation qui, si elle n'est pas complète, a pu ainsi être poussée un peu plus loin. Elle n'en met pas moins en évidence un défaut d'indépendance du système d'expressions mathématiques qui apparaît en l'occurrence incapable d'établir par ses propres moyens cette unicité[72].

---

71 La démonstration de l'unicité est inévitablement confrontée au problème de l'expression d'une totalité, en l'occurrence la totalité des répartitions de chaleur dans la lame. C'est la lame elle-même, comme précédemment le segment de sa base, et l'expression $\phi(x,y)$ qui en font ici office. Pour une autre analyse voir Dhombres & Robert [1998, 541-544].

**ii) Le Discours préliminaire : une philosophie du langage mathématique**

La *Théorie analytique de la chaleur* est précédée d'un « Discours préliminaire » dont la conformité de l'Analyse mathématique aux phénomènes naturels est un thème central. Son examen va permettre de retrouver dans la philosophie de la nature défendue par Fourier cette caractéristique, et les autres, des énoncés inauguraux.

Le Discours commence par un panorama des progrès de la mécanique rationnelle, en particulier depuis Newton, et des phénomènes naturels dont elle a pu rendre compte : le mouvement et la forme des astres, les marées, les phénomènes vibratoires, la propagation de la lumière, etc. Mais cette mécanique n'a pas jusqu'à présent, selon Fourier, été appliquée aux phénomènes thermiques qui « *composent un ordre spécial de phénomènes* » (xvi) dont il fait valoir l'importance à la fois domestique (chauffage) et climatique : répartition de la température à la surface de la Terre et dans le sol, dans l'atmosphère et dans les mers, et leur évolution dans le temps (en supposant toujours que les températures sont devenues permanentes...). Comme Newton, il entend faire une analyse exempte d'hypothèse[73] sur la nature des phénomènes même si, à nouveau comme celui-ci, il adopte le modèle moléculaire adapté au calcul différentiel. Ces phénomènes sont aussi supposés être régis par des lois dont il convient de trouver « *l'expression mathématique* » (xix), de sorte que « *tous les phénomènes* » soient interprétés « *par le même langage* » (xxiv).

L'Analyse mathématique est ainsi conçue comme un *langage*[74]. Ce langage donne une expression conforme des phénomènes dont il est tenu séparé. Il est aussi nécessaire à l'expression de leurs lois[75]. Il doit être constitué et ne saurait l'être

---

72 Ce n'est pas la seule entorse dans ce chapitre à la séparation entre l'analyse mathématique et physique : Fourier justifie de limiter dans les séries la sommation à des entiers *positifs* par le fait que les termes $e^{-mx}$ de la solution générale ne doivent pas, *pour des raisons physiques,* devenir infinis (art. 168).

73 « Les principes de cette théorie sont déduits, comme ceux de la Mécanique rationnelle, d'un très petit nombre de faits primordiaux, dont les géomètres ne considèrent point la cause, mais qu'ils admettent comme résultant des observations communes et confirmées par toutes les expériences. » (xxi).

74 Sur l'analyse mathématique comme langage à la fin du 18ème siècle, voir Brian [1994].

75 « Les effets de la chaleur sont assujettis à des lois constantes que l'on ne peut découvrir sans le secours de l'Analyse mathématique. La Théorie que nous allons exposer a pour objet de démontrer ces lois ; elle réduit toutes les recherches physiques sur la propagation de la chaleur à des questions de Calcul intégral dont les éléments sont donnés par l'expérience. » (art. 1). Remarquons incidemment que Fourier utilise lui-même les verbes « exprimer » et « représenter » (indices d'une sans une étude des phénomènes[76] :

> « Il ne peut y avoir de langage plus universel et plus simple, plus exempt d'erreurs et d'obscurités, c'est-à-dire plus digne d'exprimer les rapports invariables des êtres naturels.
>
> Considérée sous ce point de vue, l'Analyse mathématique est aussi étendue que la nature elle-même ; elle définit tous les rapports sensibles, mesure les temps, les espaces, les forces, les températures » xxiii

On retrouve dans ce Discours un dualisme, avec d'un côté la totalité des phénomènes naturels et de l'autre ce langage, l'Analyse mathématique. On retrouve la conformité. Leur différence de nature les rend aussi incommensurables.

Devant rendre compte du fait que ce langage est commun à des phénomènes aussi différents que le mouvement de la lumière dans l'atmosphère, la diffusion de la chaleur et les probabilités, Fourier stipule qu'il nous en fait connaître « *les éléments fondamentaux* » (xxiii) : s'il y a une expression commune c'est, selon lui, qu'il y a un contenu commun. Et inversement, cette expression commune est là « *comme pour attester l'unité et la simplicité du plan de l'univers, et rendre encore plus manifeste cet ordre immuable qui préside à toutes les causes naturelles* » (xxiv). Ainsi, en dépit de l'incommensurabilité de ces divers phénomènes entre eux, une totalité uniforme est constituée qui permet de préserver la conformité et la dualité avec l'Analyse mathématique.

La philosophie du langage mathématique développée ici est en adéquation complète avec les caractéristiques de la *Théorie analytique de la chaleur* qui ont été dégagées indépendamment. Il y a ainsi une homologie remarquable entre la philosophie ou l'épistémologie déployées et les caractéristiques sémiotiques des textes qu'elles introduisent[77].

### iii) L'inauguration d'une totalité

L'inauguration des séries trigonométriques par Fourier a ceci de particulier comparée à d'autres que les expressions de la représentation inaugurée *ne sont pas*, en un sens, nouvelles. Fourier n'inaugure pas les séries trigonométriques dans la mesure où il n'a pas à expliquer la formation de ces expressions ni à introduire

---

dualité) : « l'équation (...) exprime le mouvement uniforme de la chaleur etc. » (art. 321), « l'équation qui exprime le mouvement de la chaleur etc. » (art.238), « le mouvement de la chaleur dans un cylindre solide (...) est représenté par les équations etc. » (art. 306), « l'équation (...) qui représente le mouvement de la chaleur dans un solide etc. » (art. 333) etc. La phrase déjà citée « *En général, on ne pourrait apporter aucun changement dans la forme de nos solutions, sans leur faire perdre leur caractère essentiel, qui est de représenter les phénomènes.* » art. 428, reprend un passage d'une lettre à Laplace : « *on ne pourrait point faire de changement à cette solution sans qu'elle cesse d'exprimer le phénomène* » où « représenter » [Herivel 1980, 26] est remplacer par « exprimer ».

76 « L'étude approfondie de la nature est la source la plus féconde des découvertes mathématiques. (...) [Elle] est encore un moyen assuré de former l'Analyse elle-même, et d'en découvrir les éléments qu'il nous importe le plus de connaître, et que cette science doit toujours conserver : ces éléments fondamentaux sont ceux qui se reproduisent dans tous les effets naturels » (xxii)

77 Sur ce point Brian [1994, 49-71].

les conventions nécessaires à leur compréhension. Toutes leurs caractéristiques sont déjà reçues et familières. Descartes devait au contraire indiquer à son lecteur comment mettre sous forme d'équation polynomiale un problème de géométrie et comment tracer la courbe solution à partir de cette équation et préciser leur rapport. Cela, Fourier n'a pas à le faire. Cette particularité de l'inauguration des séries trigonométriques fait cependant bien ressortir un trait en revanche commun à toutes ces inaugurations qui est l'inauguration d'une *totalité en tant que totalité*. Et en cela, il y a bien inauguration des séries trigonométriques : c'est bien à la *totalité* des séries trigonométriques que Fourier a à faire. C'est leur totalité qui acquiert un sens nouveau, c'est elle qui va pouvoir dès lors exister en tant que telle et qui va devenir une nouvelle totalité pré-établie. L'inauguration de cette totalité diffère de l'introduction de ses termes. Il est sans doute intéressant de recenser les usages de séries trigonométriques *avant* Fourier ([Grattan-Guinness 1970, 19]), mais il convient aussi de distinguer ces occurrences de celles de leur totalité : les conditions de leur introduction et leurs implications ne sont pas les mêmes. Il y a là une différence objective et significative : à partir de Fourier les séries trigonométriques ont pu être considérées en tant que *totalité* et des développements leur ont été consacrés *en tant que telle*. L'étude de la convergence de ces séries en est une exemple ([Sachse 1880]) : il ne s'agit pas d'étudier la convergence d'une série trigonométrique mais celle *des* séries trigonométriques prises comme totalité. C'est une forme d'anachronisme de voir dans une série trigonométrique la totalité de ces séries avant que cette totalité n'ait été considérée, voire inaugurée, comme telle. Cet anachronisme ignore autant qu'il présuppose cette inauguration. Il peut dès lors en devenir un indice comme le sont tous les développements, mathématiques ou historiques, qui impliquent cette totalité.

### *iv) Conformité et inauguration*

La postérité de la représentation des fonctions arbitraires par les séries trigonométriques, pas plus que celles des problèmes de géométrie et des courbes géométriques par des polynômes, n'a pas tenu à leur *conformité totale*[78]. Les

---

78 Le jeune Liouville [1830, 135] crut semble-t-il à "*la possibilité de représenter une fonction quelconque, entre des limites données, par une série dont les termes sont les intégrales d'une équation linéaire du second ordre*". Lejeune-Dirichlet [1829, 157] évoque aussi « *les séries de sinus et de cosinus, au moyen desquelles on peut représenter une fonction arbitraire dans un intervalle donné* ». Pour G. Libri [1831, 224] la « *question de la discontinuité des fonctions arbitraires, qui complètent les intégrales des équations aux différentielles partielles, agitée d'abord entre Daniel Bernoulli, Euler et d'Alembert, et discutée depuis par les plus grands géomètres, parait avoir été résolue par M. Fourier qui a montré le premier comment l'on pouvait déterminer, dans chaque cas particulier, les fonctions arbitraires de manière à satisfaire aux conditions initiales du problème, même lorsque celles-ci n'obéissent pas aux lois de continuité. Les formules que cet illustre géomètre a trouvées sont propres à exprimer une fonction discontinue quelconque, dont les diverses parties, comprises entre des limites données de la variable, suivent une marche dissemblable, et sont représentées par des expressions différentes* ». Voir aussi Libri [1833, 303] où cette fois le nom de Fourier est associé à celui de Poisson qui est lui-même convaincu de la possibilité de développer une fonction arbitraire de plusieurs manières, et notamment en série trigonométrique [Poisson 1835, 167-168, 185-233]. Philip Kelland [1837, 56] expose le « *general process of expressing a discontinuous function by means of a* polynômes ou les séries trigonométriques ont néanmoins acquis le statut de totalités pré-établies et de nombreux développements leur ont en tant que tels été ensuite consacrés. Ces représentations ont depuis largement *fait leur preuve*. Leurs multiples usages suffisent depuis longtemps à justifier leur introduction. Mais leur conformité totale n'y a plus guère sa part. Il est même devenu facile de la mettre en défaut. Un premier risque est alors d'ignorer la part de conformité acquise et le rôle qu'elle continue de jouer : Fourier dispose pour les séries trigonométriques d'une certaine conformité aux courbes dont il tire toujours parti. Un second risque est de ne voir au mieux dans la présence de cette conformité au moment de l'inauguration qu'une « erreur de genèse», et d'en appeler à un contexte historique et des partis pris philosophiques plus ou moins propres à ceux qui inaugurent ces représentations. Ces contextes et ces partis pris philosophiques existent et jouent un rôle, mais la récurrence de la conformité invite à lui accorder plus d'attention et à la considérer comme une condition sinon nécessaire tout du moins remarquablement récurrente de l'introduction de nouvelles représentations en mathématiques. Le point est ici de considérer *en tant que tels* ces moments d'introduction de nouveaux moyens d'expression, c'est-à-dire d'en reconnaître la spécificité, spécificité que les usages que nous faisons de ces représentations, souvent à l'origine de notre intérêt pour leur histoire, nécessairement ultérieurs et qui supposent pourtant leur introduction, ne permettent pas de reconnaître voire nous conduisent à ignorer.

---

*trigonometric series* ». William Thomson [1839-1841/1882, 1] considère aussi que le développement d'une fonction « *completely arbitrary* » en série trigonométrique a été « *rigorously demonstrated* » d'abord par Fourier puis par Poisson. Riemann [1854, trad. fr, 226] est plus prudent : « *Les séries trigonométriques, ainsi appelées par Fourier, c'est-à-dire de la forme (…) jouent un rôle considérable dans la partie des Mathématiques où l'on rencontre des fonctions entièrement arbitraires ; on est même fondé à dire que les progrès les plus essentiels de cette partie des Mathématiques, si importante pour la Physique, ont été subordonnés à la connaissance plus exacte de la nature de ces séries. Dès les premières recherches mathématiques qui ont conduit à la considération des fonctions arbitraires, s'est posée la question de savoir si une fonction entièrement arbitraire pouvait se représenter par une série de la forme ci-dessus.* »

# III - Courbes et fonctions arbitraires : les séries trigonométriques dans la controverse des cordes vibrantes

Nous avons montré comment la possibilité de la représentation d'une fonction quelconque par une série trigonométrique s'introduisait et était aussi inaugurée dans *La théorie analytique de la chaleur*. La possibilité de représenter toute fonction par une série trigonométrique avait déjà été envisagée. Fourier évidemment ne l'ignorait pas[79] :

> « Si l'on applique ces principes à la question du mouvement des cordes vibrantes, on résoudra les difficultés qu'avait d'abord présentées l'analyse de Daniel Bernoulli. La solution donnée par ce géomètre suppose qu'une fonction quelconque peut toujours être développée en série de sinus ou de cosinus d'arcs multiples. Or, de toutes les preuves de cette proposition, la plus complète est celle qui consiste à résoudre en effet une fonction donnée en une telle série dont on détermine les coefficients. » art. 230

Il savait aussi que c'était une « c*onclusion que le célèbre Euler a toujours repoussée* » (Fourier 1805/1972, 183 ; 1822, 588)[80]. Les exemples de développements en séries trigonométriques qu'il donne dans ses mémoires ou dans sa correspondance ([Herivel 1980]) sont autant de contre-exemples aux objections qu'il savait avoir été avancées contre la généralité de ces séries[81]. L'énoncé inaugural soutenu par Fourier a en effet déjà été énoncé dans les années

---

[79] Voir notamment Fourier [1807/1972, art. 75] « *Remarques diverses sur la nature des développements précédents et sur les difficultés que présente l'équation du mouvement des cordes sonores* ».

[80] « *Les solutions que l'on obtient par cette méthode sont complètes et consistent dans des intégrales générales. Aucune autre intégrale ne peut avoir plus d'étendue. Les objections qui avaient été proposées à ce sujet sont dénuées de tout fondement ; il serait aujourd'hui superflu de les discuter.* » art. 428-6°. Ou encore : « *Cette notion n'est point opposée aux principes généraux du Calcul ; on pourrait même en trouver les premiers fondements dans les écrits de Daniel Bernoulli, de Clairaut, de Lagrange et d'Euler. Toutefois on avait regardé comme manifestement impossible d'exprimer en séries de sinus d'arcs multiples, ou du moins en séries trigonométriques convergentes, une fonction qui n'a de valeurs subsistantes que si celles de la variable sont comprises entre certaines limites, et dont toutes les autres valeurs seraient nulles. Mais ce point d'analyse est pleinement éclairci ; et il demeure incontestable que les fonctions séparées, ou parties de fonctions, sont exactement exprimées par des séries trigonométriques convergentes, ou par des intégrales définies. Nous avons insisté sur cette conséquence dès l'origine de nos recherches jusqu'à ce jour, parce qu'il ne s'agit point ici d'une question abstraite et isolée, mais d'une considération principale, intimement liée aux applications les plus utiles et les plus étendues. Rien ne nous paru plus propre que les constructions géométriques à démontrer la vérité de ces nouveaux résultats, et à rendre sensibles les formes que l'Analyse emploie pour les exprimer.* » 428-13°. Voir aussi art. 230.

[81] « *A l'égard des recherches de D'Alembert et d'Euler ne pourrois-je point ajouter que s'ils ont connu ces développements il n'en ont fait qu'un usage bien imparfait, car ils étoient persuadés l'un et l'autre qu'une fonction arbitraire et discontinue ne pourroit jamais être résolue en séries de ce genre.* » Fourier à Lagrange (?)[Herivel 1980, 21].

1750 au cours de la controverse sur les cordes vibrantes.

Nous voudrions maintenant préciser le statut de cet énoncé dans cette controverse en identifiant comme nous l'avons fait pour *La théorie analytique de la chaleur* les moyens d'expression de la généralité des courbes et des fonctions en présence et leurs relations avérées au cours de celle-ci. Nous verrons notamment le rôle joué par l'équation différentielle d'une corde vibrante établie par D'Alembert. Nous verrons l'incidence de l'expression générale de la solution qu'il en a donnée. Nous verrons aussi que l'assimilation du problème à son équation différentielle change le statut de la forme initiale de la courbe et conduit, comme la distribution initiale de chaleur dans la barre chez Fourier, à l'introduction de la possibilité de représenter une fonction quelconque par une série trigonométrique. Nous verrons enfin qu'il s'agit d'un énoncé inaugural qui n'a été introduit que pour être contesté et qui n'a vraiment été soutenu, pour des raisons parfois très différentes, par aucun des mathématiciens ayant pris part à cette controverse.

## 1 - L'énoncé inaugural dans la controverse

La controverse sur les cordes vibrantes opposa principalement D'Alembert, Euler, Daniel Bernoulli et Lagrange[82]. La description mathématique d'une corde vibrante conduit en particulier à la question de savoir s'il est possible de donner une expression mathématique unique décrivant le mouvement continu d'une corde qui peut avoir une forme initiale quelconque, avec des tangentes discontinues (corde pincée) et avoir à la fin de son excursion la forme d'une droite. Cette controverse mit aux prises quelques-uns des plus grands mathématiciens de l'époque qui à cette occasion s'opposèrent notamment sur la représentation générale des fonctions et l'extension de leur représentation par des séries entières ou trigonométriques[83].

Au début du siècle, Taylor avait déjà appliqué sa méthode des incréments à l'étude du mouvement d'une corde vibrante[84]. Il avait montré, en se limitant à des oscillations infiniment petites, que la compagne de la cycloïde, c'est-à-dire une sinusoïde, était une solution de ce problèmes et avait défendu que la corde devait rapidement prendre cette forme. D'Alembert [1749] revint ensuite sur ce problème. Il donna l'équation différentielle partielle qui le régit[85]. Il montra qu'il

---

82 Le terme de « controverse » ne renvoie pas ici à un genre historiographique qui a pu être à la mode ; il est en l'occurrence employé par les protagonistes, ainsi D'Alembert écrivait à Lagrange : « *J'ai parlé de vous dans ma Préface [Opuscules mathématiques, t. IV] comme je le dois à tous égards, à l'occasion de notre controverse sur les vibrations des cordes.* » D'Alembert à Lagrange, 29 avril 1768, Oeuvres de Lagrange, vol. 13, 111.

83 Ils s'opposèrent aussi sur leurs conceptions du rapport entre mathématiques et physique, et en particulier sur leurs conceptions de la théorie musicale et de l'acoustique [Darrigol 2007].

84 Voir Taylor [1713, 1715] ; Truesdell [1960, 129-132] ; Cannon & Dostrovsky [1981, 15-22].

85 L'expression $(\frac{ddy}{dt^2}) = cc(\frac{ddy}{dx^2})$ est donnée par Euler [1755, 211]. D'Alembert écrit quant à lui simplement $\alpha = \beta$ où $\alpha$ et $\beta$ sont les fonctions de $t$ et de $s$ définies par les différentielles complètes $dp = \alpha dt + \nu ds$ et $dq = \nu dt + \beta ds$, les fonctions $p$ et $q$ étant elles-mêmes les fonctions définies par la différentielle complète $d[\phi(t,s)] = p dt + q ds$ où $\phi(t,s)$ est l'expression de l'ordonnée de la corde. Il n'y a pas dans ces mémoires de D'Alembert de notation propre aux différentielles partielles d'une fonction. Je parlerai néanmoins admettait une infinité de solutions et contesta qu'elles prennent rapidement la forme d'une sinusoïde. Le même problème fut peu après repris par Euler [1749] avec des conditions initiales plus générales que celles que D'Alembert avait admises (et qu'il était prêt à admettre). Euler démontra que si les séries trigonométriques étaient bien des solutions, elles n'en étaient cependant qu'une partie. Daniel Bernoulli est ensuite intervenu et a défendu dans deux mémoires paru dans le même volume de l'*Histoire de l'Académie Royale des Sciences et des Belles-Lettres de Berlin* que les solutions, pour des oscillations infiniment petites, s'obtenaient toutes comme somme, éventuellement infinie, des solutions de Taylor, c'est-à-dire, et en se servant des « *dénominations de M. Euler* » [Bernoulli 1755, 156], sous la forme de séries trigonométriques [Bernoulli 1755, 157] :

$$y = \alpha \sin\left(\frac{\pi x}{a}\right) + \beta \sin\left(\frac{2\pi x}{a}\right) + \gamma \sin\left(\frac{3\pi x}{a}\right) + \delta \sin\left(\frac{4\pi x}{a}\right) + etc.$$

Dans un mémoire inséré dans le même volume que ceux de Bernoulli, Euler contesta que celui-ci ait pu ainsi donner *toutes* les solutions :

> « M. Bernoulli tire toutes ces excellentes réflexions uniquement des recherches, que feu M. Taylor a faites sur le mouvement des cordes, & soutient contre M. D'Alembert & moi, que la solution de Taylor est suffisante à expliquer tous les mouvements, dont une corde est susceptible ; de sorte que les courbes, qu'une corde prend pendant son mouvement, soit toujours, ou une trochoïde allongée simple, ou un mêlange de deux ou plusieurs courbes de la même espece. » Euler [1755, 196-7]

Il précise un peu plus loin sa critique :

> "Mais il y a plus : je n'avois donné cette équation,
>
> $$y = \alpha \sin\left(\frac{\pi x}{a}\right) + \beta \sin\left(\frac{2\pi x}{a}\right) + \gamma \sin\left(\frac{3\pi x}{a}\right) + \delta \sin\left(\frac{4\pi x}{a}\right) + etc.$$
>
> que comme une solution particulière de la formule, qui contient en général toutes les courbes, qui peuvent convenir à une corde mise en mouvement ; & il y a une infinité d'autres courbes, qui ne sauroient être comprises dans cette équation. Si M. Bernoulli tomboit d'accord là dessus, il n'auroit pas avancé, que toutes les courbes d'une corde frappée résultoient uniquement de la combinaison de deux ou plusieurs courbes Tayloriennes ; & il auroit reconnu, que le raisonnement fondé sur cette combinaison n'est pas suffisant à fournir une solution complette de la question dont il s'agit." Euler [1755, 198]

Dénoncer l'incomplétude des solutions de Bernoulli est l'objet de son mémoire :

> « La question principale, que j'ai à développer, est donc : si toutes les courbes d'une corde mise en mouvement sont comprises dans l'équation rapportée, ou non? » Euler [1755, 198].

Euler rappelle ici qu'il a aussi trouvé ces solutions mais qu'il ne les tient que pour des solutions particulières et conteste que ce soient là toutes les solutions. Il émet alors l'hypothèse suivante :

> « Mais peut-être repliquera-t-on, que l'équation y=αsin πx/2a+ &c. à cause de l'infinité de coëfficiens indéterminés, est si générale, qu'elle renferme toutes les

---

dans les deux cas d'*équations aux dérivées partielles*, les différences éventuelles entre les moyens d'expression utilisés n'ayant pas d'incidence sur mon propos ici.

courbes possibles » Euler [1755, 200]

Euler énonce ici l'*énoncé inaugural* que l'on a vu avoir été soutenu par Fourier. Dans un manuscrit de 1755, « Observations sur deux Mémoires de M. Euler et Daniel Bernoulli, insérés dans les mémoires de 1753 »[86], D'Alembert répondit aux solutions et objections de Bernoulli et d'Euler. La série de critiques qu'il adresse à Bernoulli commence par le même énoncé inaugural, qu'il énonce et conteste à son tour :

> « Il me paroit impossible[87] de prouver que toutes les courbes de la corde vibrante puissent être renfermées dans l'équation
>
> $$y = \alpha \sin(\frac{\pi x}{a}) + \gamma \sin(\frac{3 \pi x}{a}) + \delta \sin(\frac{4 \pi x}{a}) + etc$$
>
> en supposant même que cette courbe de la corde vibrante soit méchanique et composée de parties égales et semblables situées alternativement au dessus et au dessous de l'axe et liées par la loy de continuité. » D'Alembert , 24 ; Jouve 2007, II, 97-8.

Plus tard, Fourier considèrera lui aussi que la solution de Bernoulli suppose l'énoncé inaugural qu'il soutiendra :

> « La solution donnée par ce géomètre [Daniel Bernoulli] suppose qu'une fonction quelconque peut toujours être développée en séries de sinus ou de cosinus d'arcs multiples. » Fourier [1822, art. 230]

Ces citations d'Euler, de D'Alembert et de Fourier suggèrent que Daniel Bernoulli aurait soutenu que toute courbe pourrait être représentée par une série trigonométrique. Cependant, cet énoncé n'a pas été introduit par Bernoulli mais par Euler et D'Alembert (la priorité de l'un ou de l'autre étant sans importance ici), et seulement pour le contester, chacun pour des raisons différentes[88].
Nous allons à présent montrer que cet énoncé a été introduit par Euler (les mêmes arguments valant aussi, à des nuances sans importance pour notre propos, pour D'Alembert) et en mettre en évidence les raisons. Avant cela, nous rappellerons néanmoins le principe des systèmes oscillants introduit et défendu par Bernoulli afin d'indiquer que ses mémoires se rapportent bien plutôt à ce principe qu'à l'énoncé inaugural sur les séries trigonométriques avec lequel il ne peut être confondu.

## 2 - Le principe des systèmes oscillants de Bernoulli

L'analyse des cordes vibrantes proposée par Bernoulli est fondée sur le principe suivant :

> « dans tout système les mouvements réciproques [périodiques] des corps sont

---

86 Ce manuscrit a été édité dans Jouve [2007, II, 77-113]. Soumis en 1755 à l'Académie de Berlin, le Mémoire, jugé trop polémique, sera refusé et servira de base au premier mémoire des *Opuscules* [D'Alembert 1761 ; 2008]. Pour l'histoire de ce mémoire, voir Jouve [2008].
87 La suite du texte rend clair que « impossible » veut dire ici que cela est faux.
88 Cette position est aussi celle défendue par Burkhardt [1908, 20 ; Jouve 2007, II, 29] .

toujours un mélange de vibrations simples, régulières et permanentes de différentes espèces. » Bernoulli [1755, 195]

C'est sur ce principe « *d'une très grande étendue* » [Bernoulli 1758, 160] qu'il fonde sa « *théorie sur la pluralité des vibrations coexistentes dans un seul & même système* » [Bernoulli 1758, 157][89] et à partir duquel il dérive sa « *méthode [qui] s'étend à déterminer les vibrations & mouvemens réciproques dans tous les systèmes de corps pour lesquels on peut déterminer les vibrations simples ; c'est-à-dire, telles que toutes les parties fassent des vibrations parfaitement synchrones, & chacune la sienne, suivant la nature des petites oscillations d'un pendule simple* » [Bernoulli 1758, 158]. Ce principe de composition des vibrations est pour lui « *le seul & vrai principe* », il est « *général, sans contredit, pour les cordes uniformément épaisses (…) jusqu'aux cordes inégalement épaisses* » [Bernoulli 1767, 306]. Sa validité va d'ailleurs au-delà du problème des cordes vibrantes. Suivant ce principe, tous les sons, harmonieux ou non, se composent d'harmoniques qui peuvent être isolées comme Newton a séparé les rayons primitifs qui composent la lumière [Bernoulli 1758, 158]. Mais il s'agit là d'une *principe physique,* un « *principe sur la coexistence des vibrations* » [Bernoulli 1758, 160], et non d'un *énoncé inaugural*. Il y a bien réalisme, il y a bien inauguration, mais il n'y aucun dualisme ni aucune incommensurabilité. La question de la conformité, si elle se pose, se pose dès lors aussi différemment.

Le mémoire de Bernoulli publié en 1755, ainsi que les autres sur le même sujet, sont en grande partie consacrés à défendre et illustrer ce *principe*. Bernoulli donne ainsi divers arguments tendant à montrer qu'un son est bien la superposition de sons fondamentaux qui s'entendent distinctement et se combinent sans se mélanger[90]. Il défend qu'il obtient ainsi toutes les solutions

---

89 « *J'ai évité jusque ici les calculs, & j'ai fondé tout mon raisonnement sur le principe confirmé par l'expérience (§6) qu'il peut se faire un mêlange de vibrations dans un seul & même corps sonore, qui soient absolument indépendantes les unes des autres.* » Bernouilli [1755, 160] ; « *il semble que la Nature n'agit très souvent que par les principes des vibrations isochrones imperceptibles, & infiniment diversifiées, pour produire un grand nombre de Phénomènes.* » Bernoulli [1755b, 173] ; « *J'ai démontré de plus dans les Mémoires de Berlin, que les vibrations de différens ordres, quels qu'on les prenne, peuvent coéxister dans une seule & même corde, sans se troubler en aucune façon, ces différentes especes de vibration coéxistentes étant absolument indépendantes les unes des autres. De là cette pluralité de sons harmoniques qu'on entend à la fois d'une seule & même corde. Si toutes especes de vibration commencent au même instant, il arrivera que la premiere vibration du premier ordre, la seconde vibration du second ordre, la troisième vibration du troisième ordre &c finiront au même instant. C'est là un synchronisme apparent dans un certain sens, & qui n'est rien moins que général, puisqu'il y a une infinité d'autres vibrations qui ne finissent pas au même instant.* » Bernoulli [1767, 283].

90 « *Effectivement tous les Musiciens conviennent, qu'une longue corde pincée donne en même tems, outre son ton fondamental, d'autres tons beaucoup plus aigus ; ils remarquent surtout le mêlange de la douzième & de la dix-septième majeure : s'ils ne remarquent pas aussi distinctement l'octave & la double octave, c'est n'est qu'à cause de la trop grande ressemblance de ces deux tons avec le ton fondamental. Voilà une preuve évidente, qu'il peut se faire dans une seule & même corde un mêlange de plusieurs sortes de vibrations Tayloriennes à la fois. On entend pareillement dans le son des grosses cloches un mêlange de tons differens. Si l'on tient par le milieu une verge d'acier, & qu'on la frappe, on entend à la fois un mêlange confus de plusieurs tons, lesquels étant appréciés par un habile Musicien se trouvent extrêmement desharmonieux, de sorte qu'il se forme un concours de vibrations, qui ne commencent & ne finissent jamais dans un même instant, sinon par un grand hazard : d'où l'on voit que l'harmonie des sons, qu'on entend dans une même corps sonore à la fois, n'est pas essentielle à cette matière, & ne doit pas servir de principe pour les systèmes de Musique. L'air n'est pas exemt de cette multiplicité de sons* possibles et affirme que les solutions proposées par D'Alembert et Euler sont de cette forme[91], qu'elles sont donc déjà contenues dans les siennes, et qu'il peut ainsi leur donner un sens physique[92]. L'exemple d'une corde avec 1001 ventres lui permet de montrer aussi la diversité de ses courbes[93]. Il défend inversement de manière récurrente que tous les sons correspondant à ses séries peuvent être produits par des corps sonores[94]. Et c'est fort de la conviction que son principe rend compte de la composition des sons qu'il discute les solutions proposées par les deux Géomètres[95].

Ainsi, Bernoulli fonde-t-il sa résolution du problème des cordes vibrantes sur son principe de composition de vibrations fondamentales, ce problème contribuant inversement à en établir la validité et l'intérêt. C'est ce principe physique que Bernoulli défend. Ses arguments sont eux-mêmes essentiellement physiques, étayés par des expériences sonores. Quand ils ne sont pas exclusivement physiques, comme quand il esquisse une démonstration de l'interpolation d'une courbe par des sinusoïdes pour répondre à ses détracteurs, ils sont alors exclusivement géométriques. Il n'y a pas plus de dualité.

Un énoncé inaugural est remarquable par le fait de soutenir l'existence d'une représentation uniforme et conforme d'une totalité. Il y a là une prétention exorbitante. Un énoncé inaugural, même pour celui qui le soutient, reste un énoncé incroyable. Celui qui l'énonce doit en répondre et d'une certaine manière y faire face. Les textes inauguraux y répondent explicitement. Fourier, on l'a vu, ne manque pas de marquer son étonnement et de développer une succession d'arguments spécifiques pour soutenir son énoncé inaugural. Son texte témoigne d'une démarche inaugurale dont la finalité ne fait pas de doute. De même, tout son « Discours préliminaire » présente une philosophie du langage mathématique qui répond à cette prétention exorbitante, étendue à toute l'Analyse mathématique. Il en est de même des articles de Church, de Turing, de la *Begriffsschrift* de Frege ou encore de *La Géométrie* de Descartes. On ne trouve en revanche dans les mémoires de Bernoulli aucune reconnaissance de ce caractère exorbitant, aucun argument n'est donné pour y répondre. Ses arguments sont donnés pour défendre un *principe physique* sur la nature des phénomènes vibratoires et non pour inaugurer la représentation conforme des fonctions par des séries trigonométriques.

---

*coëxistans : il arrive souvent qu'on tire deux sons differens à la fois d'un tuyau ; mais, ce qui prouve le mieux, combien peu les différentes ondulations de l'air s'entre-empêchent, est qu'on entend distinctement toutes les parties d'un concert, & que toutes les ondulations causées par ces différentes parties se forment dans la même masse d'air sans se troubler mutuellement, tout comme les rayons de la lumière, qui entrent dans une chambre obscure à travers une petite ouverture, ne se troublent point.* » Daniel Bernoulli [1755, 152-3]

91 « *si j'ai bien compris leurs énoncés, toutes les nouvelles courbes qu'ils donnent, sont comprises dans notre construction, & sont un simple mêlange de plusieurs especes de vibrations, dont chacune à part se fait suivant les loix décrites par M. Taylor. Mais il me semble que ce n'est là qu'une espece de composition de mouvement, qui ne peut donner aucune amplification à la théorie de M. Taylor.* » Bernoulli [1755, 155] ; « *Voyons encore si les nouvelles courbes trouvées par M. Euler, sont comprises dans notre remarque.* » Bernoulli [1755, 156].

92 « *mon intention n'a été principalement que d'exposer ce que les nouvelles vibrations de Mrs D'Alembert & Euler ont de physique* » Bernoulli[1755, 158].

93 « *Pour mieux sentir l'incongruïté d'une telle amplification, nous combinerons la courbe fondamentale de M. Taylor, qui est représentée par la première figure, avec la figure anguiforme Taylorienne qui auroit 1001 ventres* » Bernoulli [1755, 155]

94 « *Ma conclusion est, que tous les corps sonores renferment en puissance une infinité de sons, & une infinité de manières correspondantes de faire leurs vibrations régulières ; enfin, que dans chaque différentes espece de vibrations les inflexions des parties du corps sonore se font d'une manière différente.* » Bernoulli [1755, 151]

95 « *aussi n'est-ce à mon avis que sous cette forme que les vibrations peuvent devenir régulières, simples, et isochrones, malgré l'inégalité des excursions. Avec cette idée, que j'ai toujours euë, je ne pouvois qu'être surpris de voir dans les Mémoires des années 1747 & 1748, une infinité d'autres courbures comme douées de la même propriété ; il ne me falloit pas moins que les grands noms de Mrs D'Alembert & Euler ; que je ne pouvois soupçonner d'aucune inattention, pour examiner s'il n'y auroit pas quelque équivoque dans l'aggrégation de toutes ces courbes avec celle de M. Taylor, & dans quel sens on pourroit les admettre. J'ai vû aussi-tôt, qu'on ne pouvait admettre cette multitude de courbes que dans un sens tout à fait impropre ; je n'en estime pas moins les calculs de Mrs D'Alembert & Euler, qui renferment certainement tout ce que l'Analyse peut avoir de plus profond & de plus sublime ; mais qui montrent en même tems, qu'une analyse abstraite, qu'on écoute sans aucun examen synthétique de la question proposée, est sujette à nous surprendre plutôt qu'à nous éclairer. Il me semble à moi, qu'il n'y avoit qu'à faire attention à la nature des vibrations simples des cordes, pour prévoir sans aucun calcul tout ce que ces deux grands Géomètres ont trouvé par les calculs les plus épineux & les plus abstraits, dont l'esprit analytique se soit encore avisé.* » Bernoulli [1755, 148].

## 3 - L'énoncé inaugural d'Euler et l'introduction de l'expression d'une fonction arbitraire

L'énoncé inaugural sur les séries trigonométriques a pourtant bien été énoncé au cours de cette controverse. Mais il l'a été par Euler après que Bernoulli ait répondu aux mémoires de D'Alembert et à celui d'Euler publié à leur suite. Nous voudrions maintenant rendre compte du fait que cet énoncé ait été introduit par Euler.
Dans son mémoire de 1749, D'Alembert a montré que les solutions de l'équation aux dérivées partielles décrivant le mouvement d'une corde vibrante sont toutes les courbes d'équation $y=\Psi(t+x)-\Gamma(t-x)$ où $\Psi$ et $\Gamma$ sont des *fonctions quelconques*. Comme les extrémités de la corde sont aussi fixes, les solutions doivent de plus être de la forme $y=\Psi(t+x)-\Psi(t-x)$, avec $\Psi$ $2l$-périodique ($l$ longueur de la corde). D'Alembert considère aussi le cas particulier où la figure initiale de la corde est une droite (pour des oscillations infinitésimales), ce qui impose en plus que la fonction $\Psi$ soit paire. Il appelle « *courbe génératrice* » la courbe notée $\Psi$ [D'Alembert 1749, 220]. Il peut ainsi montrer comme il le voulait qu'il est possible d'engendrer une *infinité* de fonctions satisfaisant ces conditions en reportant au-dessus et en-dessous de l'axe une portion de courbe définie par une fonction génératrice qui peut ne pas être sinusoïdale. Euler [1749, 1755] reprend ces résultats et expose en détail la construction d'une solution en reportant alternativement sous l'axe et au-dessus *n'importe quelle* courbe s'annulant à ses deux extrémités.
L'objectif de D'Alembert était de montrer qu'il existait une infinité de solutions en plus de celles trouvées par Taylor. Ce qu'il réussit à faire. Euler a repris sa démarche mais en faisant valoir qu'il est non seulement possible d'obtenir une infinité de solutions mais aussi de partir de *n'importe quelle* courbe génératrice pour les obtenir. Il souligne lui-même à plusieurs reprises cet aspect :

> « Or, sans faire encore attention à ces propriétés, & m'arrêtant uniquement à l'équation (ddy/dt2)=cc(ddy/dx2), il est important de remarquer, que les deux courbes ES & FT sont absolument arbitraires, & qu'on les peut prendre à volonté ; car, quelles que soient ces deux courbes, si nous y prenons les abscisses EQ=x+ct & FR=x+ct [il faut lire FR=x-ct], & que nous posions y=QS+RT, ou bien y=n(QS+RT), il est certain que cette valeur satisfait à notre équation ; ce qu'il serait aussi aisé de prouver indépendamment de l'analyse, que je viens de développer. Or, ce qui est le principal, ces deux courbes appliquées de la manière enseignée, satisfont également, soit qu'elles soient exprimables par quelque équation, ou qu'elles soient tracées d'une manière quelconque, de sorte qu'elles ne puissent être assujetties à aucune équation. Le Lecteur est prié de réfléchir bien sur cette circonstance, qui contient le fondement de l'universalité de ma solution, contestée par M. D'Alembert. » Euler [1755, 213-4]

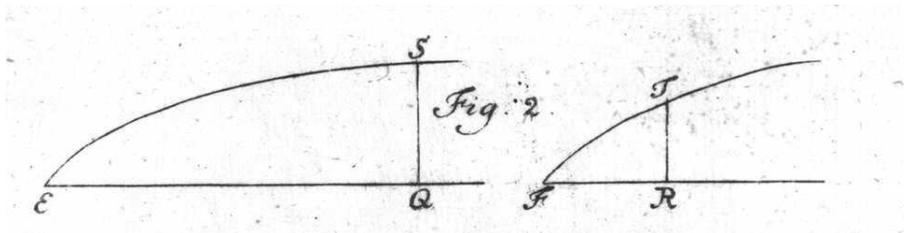

Il souligne le caractère *quelconque* de la courbe :

> « Mais, de la manière que je viens de conduire la solution, il n'est pas nécessaire, que la courbe directrice soit exprimée par quelque équation, & la seule considération du trait de la courbe suffit à nous faire connoitre le mouvement de la corde, sans l'assujettir au calcul. » Euler [1755, 217]

L'extension opérée par Euler est aussi reconnue par Lagrange :

> « M. Euler est, je crois, le premier qui ait introduit dans l'Analyse ce nouveau genre de fonctions, dans sa solution du problème de *chordis vibrantibus* » Lagrange [1760, œuvres, 1, 158].

La solution du problème donnée par D'Alembert, $y=\Psi(t+x)-\Psi(t-x)$, comprend l'expression d'une fonction *arbitraire*[96]. Cette expression de la solution établit une correspondance entre la totalité des solutions de ce *problème particulier* et la *totalité des fonctions* quoi que cela recouvre pour chacun des mathématiciens. Le problème de la corde vibrante est donc un problème particulier dont la solution implique *toutes* les fonctions.
La solution $y=\Psi(t+x)-\Psi(t-x)$ implique cette totalité sans que celle-ci n'ait ici besoin et n'ait d'autre expression qu'une variable de fonction, $\Psi$, c'est-à-dire un simple symbole de fonction tenant lieu de n'importe quelle fonction. Il est ainsi possible d'exprimer n'importe quelle fonction sans avoir de système d'expressions (analytiques) établi pour leur totalité. Dès lors quand Bernoulli affirme que son principe lui permet de déterminer que toutes les solutions du problème sont des séries de sinus, cela implique immédiatement *pour D'Alembert et Euler* que n'importe quelle fonction devrait être une somme de sinus. Ils en tirent l'un et l'autre cette conséquence et introduisent l'énoncé afférent, mais sans pour autant le soutenir. L'introduction de cet énoncé apparaît ainsi comme une *conséquence* du fait que la solution donnée par D'Alembert implique l'expression d'une fonction quelconque. C'est l'introduction d'une variable parcourant *toutes* les fonctions qui transforme l'affirmation de Bernoulli en une affirmation sur la possibilité de représenter toutes les fonctions par des séries trigonométriques. La *totalité des fonctions* s'introduit à cette occasion.

---

[96] L'apparition d'une fonction arbitraire dans la résolution d'une équation aux dérivées partielles n'est bien sûr pas propre à ce problème ou à cette équation. Une telle solution apparaît déjà, par exemple, dans D'Alembert [1747a], voir Demidov [1982, 15].

Dans le mémoire de Bernoulli [1755] la *totalité des fonctions* n'intervient d'aucune manière. Là où quelqu'un qui inaugurerait une telle représentation y aurait sans doute fait référence, Bernoulli écrit lui « *Après ces remarques il sera facile de construire une infinité de courbes initiales* [et non toutes!] *à la corde AB avec cette condition, que chaque point de la corde arrive quelquefois en même tems à un point de repos instantané, & de donner la loi générale pour toutes ces courbes sans aucun calcul préalable.* » [Bernoulli 1755, 153].

Euler introduit quant à lui son énoncé au cours de la *discussion* des différentes solutions proposées. Cet énoncé n'est utile ni à la mise en équation, ni à la résolution de celle-ci. Il est une conséquence du rapprochement, c'est-à-dire de la considération conjointe de la solution proposée par Euler, à la suite de D'Alembert, et de celle proposée par Bernoulli[97]. Il n'est utile ni à l'une, ni à l'autre.

L'expression d'une courbe ou d'une fonction *arbitraire* intervient aussi avec la *figure initiale* de la corde. Avant de considérer le rôle que cette figure initiale *arbitraire* peut aussi jouer, il convient de préciser le rapport du problème des cordes vibrantes à son équation aux dérivées partielles.

## 4 - Le problème des cordes vibrantes et son équation aux dérivées partielles

La différence de statut de l'équation aux dérivées partielles est l'indice de certaines différences entre les quatre mathématiciens. Elle permet en particulier de préciser la différence entre ce que Lagrange appelle l'« *examen synthétique* » de Bernoulli [Lagrange 1759, œuvres, 1, 95] et l'« *Analyse abstraite* » proposée par D'Alembert, Euler et lui-même [Lagrange 1759, Oeuvres, I, 70][98]. Le problème considéré par Daniel Bernoulli a toujours été et restera celui de la description mathématique du mouvement d'une corde vibrante (en liaison avec des considérations acoustiques). Le souci d'en donner une solution mathématique conditionne l'interprétation et la formulation du problème et l'oblige par exemple à considérer des oscillations infinitésimales. Mais la représentation du problème au moyen d'une équation aux dérivées partielles obtenue par D'Alembert modifie le traitement du problème et l'acception qui peut être donnée de sa résolution.

L'introduction de l'équation aux dérivées partielles contribue à séparer deux temps en même temps qu'elle marque cette séparation : celui de la mise en équation (différentielle) et celui de la résolution de celle-ci. Au cours de la mise en équation sont prises en compte les données et les hypothèses physiques. Il faut se rapporter à la physique du problème. Au terme de cette étape, le problème est réduit à une ou plusieurs équations aux dérivées partielles avec des conditions initiales. La résolution de ces équations devient un nouveau problème qui se substitue au précédent (sans être nécessairement tenu pour conforme), et qui doit être résolu suivant les principes du calcul différentiel et intégral[99]. Les données physiques utilisées sont supposées avoir toutes été converties en grandeurs et relations mathématiques, et ce faisant éliminées : l'énoncé du problème (déterminer les fonctions satisfaisant telles équations différentielles avec telles conditions initiales) et sa résolution sont dès lors coupés du problème physique d'origine, et aucune considération physique n'intervient plus, sinon marginalement, dans sa résolution. Une telle séparation est par exemple manifeste chez Euler quand il déclare :

> « La question mécanique proposée se réduit à ce problème analytique, de chercher des fonctions $r$ & $s$ de $x$ & $t$, telles que ces formules differentielles $rdx + sdt$, & $sdx + \frac{2M}{FA} rdt$, deviennent intégrables » Euler [1750, 76-7].

---

97  Cette observation est déjà dans Sachse [1880, 45].
98  Sur D'Alembert et les équations aux dérivées partielles, voir Guilbaud & Jouve [2009].

Les équations aux dérivées partielles marquent aussi bien dans le processus de résolution que dans le texte qui l'expose une séparation entre des développements physico-mathématiques et des développements purement mathématiques.

Les équations aux dérivées partielles peuvent donner lieu à des équations autonomes : une même équation ou un même système d'équations étant associé à *différents* problèmes physiques. Lagrange nous en donne un exemple en réduisant aussi bien le problème de la corde vibrante que celui de la propagation du son au même système d'équations du premier ordre[100]. Après avoir dégagé ces équations par une analyse physico-mathématique nécessairement propre à chacun des deux problèmes et exposée dans deux chapitres distincts, il peut ensuite en proposer dans *un même* chapitre une résolution commune qui ne repose plus que sur des développements strictement mathématiques [Lagrange 1759, Oeuvres, I, 72-90]. Ce traitement commun à la fois suppose et assure que la suite de la résolution sera strictement mathématique et pure de considérations physiques.

L'autonomie de l'équation aux dérivées partielles se retrouve dans le mémoire d'Euler [1755] quand, après avoir ré-exposé l'analyse qui permet d'obtenir l'équation $(\frac{ddy}{dt^2}) = \frac{Fa}{2M}(\frac{ddy}{dx^2})$, il écrit :

> « Voilà donc à quoi le problème sur le mouvement de la corde est réduit : il s'agit de trouver pour $y$ une telle fonction des deux variables $x$ & $t$, qui satisfasse à cette équation : $(\frac{ddy}{dt^2}) = \frac{Fa}{2M}(\frac{ddy}{dx^2})$, outre qu'elle renferme les propriétés marquées cy-dessus » Euler [1755, 208].

---

99 Le rapport des relations différentielles au problème dont elles sont issues est sans aucun doute différent chez D'Alembert et Euler, mais cette différence n'a pas ici d'incidence.

100 « Il est visible que toutes ces équations [du problème de la corde vibrante, $\frac{d^2 y_{m-1}}{dt^2} = \frac{2Ph}{MT^2} \frac{-2y_{m-1} + y_{m-2}}{r}$ ] sont entièrement semblables à celles que nous avons trouvées pour les mouvements des corps élastiques, et qu'il n'y a qu'à faire P=E, pour qu'elles deviennent tout à fait les mêmes ; d'où il s'ensuit que les deux problèmes qui y répondent sont de même nature, et qu'en en résolvant un on résout l'autre en même temps. » Lagrange [1759, Oeuvres, vol. 1, 60].

Un peu plus loin, Euler parle de « *notre équation* [ $(\frac{ddy}{dt^2}) = cc(\frac{ddy}{dx^2})$ ], *qui* détermine [c'est moi qui souligne] *le mouvement de la corde* » Euler [1755, 211].

Les équations aux dérivées partielles marquent une coupure nette entre deux moments bien distincts de la résolution du problème physique initial. Il n'est pas difficile de donner d'autres manifestations de cette séparation chez D'Alembert, Euler ou Lagrange.

Cette séparation est aussi nettement marquée chez Fourier [1822]. L'autonomie de l'équation ressort par exemple clairement de cette citation :

> « On voit, par exemple, qu'une même expression, dont les géomètres avaient considéré les propriétés abstraites et qui, sous ce rapport, appartient à l'Analyse générale, représente aussi le mouvement de la lumière dans l'atmosphère, qu'elle détermine les lois de la diffusion de la chaleur dans la matière solide, et qu'elle entre dans toutes les questions principales de la Théorie des probabilités. » Fourier "Discours préliminaire", xxiii.

Cette séparation se retrouve aussi dans l'organisation du livre. La *Théorie analytique de la chaleur* présente une hiérarchie linéaire d'unités supra-phrastiques bien distinctes : chapitre, section, article. Ce sont bien les seules unités distinguées à ces niveaux puisque les théorèmes, les définitions et les démonstrations ne font l'objet d'aucune séparation systématique sur le plan de l'expression (ils ne sont pas signalés par un « théorème » ou un « définition » ou par des mises en page particulières). Les chapitres ont un titre et un numéro les sections de même, mais leur numérotation reprend à partir de 1 au début de chaque chapitre et les articles sont dépourvus de titre mais numérotés de manière continue du début jusqu'à la fin du livre. Les chapitres, c'est-à-dire la plus grande unité infra-textuelle, répondent à deux séries d'oppositions. Ils suivent d'une part les trois étapes du processus de mise en équation : définition des grandeurs physiques (chp. I), mise en équation (chp. II), résolution (chp. III-IX) dont la séparation est ainsi nettement marquée. A partir du troisième, la succession des chapitres suit celle des corps considérés. Ainsi, la plus grande unité infra-textuelle reproduit les étapes de la mathématisation avec ses trois moments bien marqués, appliquée à des corps particuliers. Elle domine dans les trois premiers chapitres l'opposition entre les différents corps considérés qui, dans le chapitre II, est maintenue au niveau des sections, avant de remonter ensuite à la surface, c'est-à-dire au niveau des chapitres (III-IX). Les unités textuelles de rang inférieur au chapitre accueillent d'autres oppositions. Elles accueillent en particulier l'entrelacement des développements purement mathématiques et la vérification de leur conformité physique.

On ne retrouve aucune séparation semblable dans les mémoires de Bernoulli. Une équation est donnée, mais celle-ci ne constitue pas un problème mathématique complet et autonome susceptible de se substituer au problème de la corde vibrante. Elle ne marque pas une séparation entre deux moments distincts et autonomes de la résolution. Euler devra, comme on le verra, reformuler le problème en y intégrant les conditions initiales pour qu'il en soit ainsi, faisant de l'équation un énoncé mathématique complet, et même conforme, au problème physique (reformulé...). Elle n'est pour Bernoulli qu'une étape parmi bien d'autres d'une analyse qui reste physico-mathématique d'un bout à l'autre de la résolution.

La séparation opérée ou non par l'équation dans la résolution du problème de la corde vibrante est un indice du statut de la démarcation entre analyse physique et mathématique. La lecture de ces mémoires et de la correspondance entre ces mathématiciens ne laisse pas de faire apparaître une opposition marquée à cet égard entre Daniel Bernoulli et les trois autres mathématiciens (sans que le statut de cette démarcation soit pour autant le même chez ces derniers). La controverse implique et révèle aussi bien sûr bien d'autres désaccords, y compris entre D'Alembert, Euler et Lagrange sur le principe de continuité et plus généralement sur les principes mêmes du calcul différentiel et intégral. La diversité des oppositions en jeu participe du caractère exceptionnel de cette controverse. Mais l'opposition entre Daniel Bernoulli et les trois autres mathématiciens tient, en partie, à la nature de leurs démarches respectives, et en particulier au rôle différent des mathématiques dans celles-ci[101] et dont la démarcation par les équations différentielles est un des indices objectifs ayant ses conséquences, notamment sur la prise en compte ou non de la courbe initiale. D'Alembert, Euler et Lagrange ne s'entendent pas entre eux sur la place exacte ou le statut de cette démarcation. D'Alembert en particulier accepte volontiers des limites plus étroites et fixes aux mathématiques que les deux autres géomètres. Mais ils partagent tous les trois l'idée d'une telle démarcation qui intervient dans leur résolution du problème de la corde vibrante. Si la présence de cette séparation dans un texte ne suffit évidemment pas à établir son caractère inaugural, il semble difficile d'inaugurer une représentation comme celle des fonctions par des séries entières, voire seulement de l'envisager, sans que celle-ci ne s'y manifeste.

## 5 - Courbes et conditions initiales

Nous pouvons à présent revenir sur la figure initiale de la corde vibrante. La figure initiale introduit aussi l'expression d'une *courbe quelconque* dans le problème des cordes vibrantes. Il importe donc d'apprécier son rôle dans l'introduction de l'énoncé inaugural. Le statut des conditions initiales des relations différentielles dépend de celui de ces relations, et en particulier de leur autonomie vis-à-vis du problème considéré. Nous avons vu qu'il n'était pas le même chez Bernoulli et les trois autres mathématiciens.
Bernoulli s'intéresse aux cordes vibrantes en tant qu'elles produisent des sons. Or, ces sons ne dépendent pas de la figure initiale de la corde. Celle-ci n'est de ce fait pas considérée par Bernoulli. Comme Taylor, il ne s'intéresse à l'aspect de la corde qu'*après* qu'elle ait atteint un certain régime vibratoire *indépendant de la*

---

[101] Lagrange décrit la voie suivie par Bernoulli comme « une espèce d'induction » Lagrange [1759, œuvres, I, 69].

*figure initiale*[102] :

> « Au reste je crois que quelque courbure initiale qu'on donne à la corde, elle ne manquera pas de faire ses vibrations presque aussitôt suivant la simple uniformité des mouvements isochrones, & conformément à la nature de la trochoïde prolongée exposée par M. Taylor » Bernoulli [1755, 158].

Ainsi, le problème des cordes vibrantes est pour Bernoulli indépendant de la figure initiale de la courbe. Le caractère arbitraire de cette courbe ne saurait donc l'amener à inaugurer une nouvelle représentation des courbes ou des fonctions.
Il n'en est pas de même pour D'Alembert et Euler qui eux résolvent des équations aux dérivées partielles. Une solution mathématique consiste pour eux à décrire le mouvement de la corde *depuis sa figure initiale* jusqu'à son évanouissement complet. La courbe initiale doit pour eux vérifier l'équation associée au problème. *Elle fait partie des solutions*. Même si leurs analyses mathématiques ne sortent pas non plus du cadre infinitésimal[103], et supposent des oscillations infiniment petites, ils y font entrer la figure initiale. Cette différence ressort nettement du passage remarquable suivant dans lequel Euler peut *remonter* du moment où le système satisfait aux conditions de l'analyse de Taylor et de Bernoulli jusqu'à sa figure initiale qui « *dépend de notre bon plaisir* » :

> « Il est effectivement prouvé d'une manière suffisante, que si une seule vibration est conforme à cette règle [avoir une forme sinusoïdale], toutes les suivantes doivent l'observer aussi. On voit en même tems par là, comment l'état des vibrations suivantes dépend des précédentes, & peut être déterminé par elles; comme réciproquement, par l'état des suivantes, on peut conclure la disposition de celles qui ont précédé. C'est pourquoi, si les vibrations suivantes sont régulières, il ne sera en aucune manière possible que les précédentes se soient écartées de la règle : d'où résulte aussi évidemment, que si la première vibration a été irrégulière, les suivantes ne peuvent jamais parvenir à une parfaite régularité. Or la première vibration dépend de notre bon plaisir, puisqu'on peut, avant que de lacher la corde, lui donner une figure quelconque ; ce qui fait que le mouvement vibratoire de la même corde peut varier à l'infini, suivant qu'on donne à la corde telle ou telle figure au commencement du mouvement. » Euler [1750, 70].

A la suite de ce passage Euler reformule l'énoncé du problème en y incluant explicitement la figure initiale :

> « *Si une corde de longueur, & de masse donnée, est tendüe par une force, ou un poids donné ; qu'au lieu de la situation droite, on lui donne une figure quelconque, qui ne differe cependant de la droite qu'infiniment peu, & qu'ensuite on la lache*

---

[102] Cette manière de considérer le problème, déjà chez Taylor, est aussi soulignée par Euler : « *A l'egard de l'autre limitation, qui suppose toutes les vibrations régulieres, on tâche de la défendre en disant, que bien qu'elles s'ecartent de cette loi au commencement du mouvement, elles ne laissent pas de s'assujettir au bout d'un très court espace de tems à l'uniformité, de sorte qu'à chaque vibration la corde s'étend tout à la fois, & ensemble en ligne droite, affectant hors de cette situation la figure d'une trochoïde prolongée.* » Euler [1750, 70].

[103] « *A la vérité la première limitation, par laquelle les vibrations de la corde sont regardées comme infiniment petites, quoique réellement elles conservent toujours une raison finie à la longueur de la corde, cela ne dérange presque en rien les conclusions qu'on en tire, parce qu'en effet ces vibrations sont pour l'ordinaire si petites, qu'elles peuvent être prises pour infiniment petites, sans qu'il en résulte d'erreur sensible. D'ailleurs on n'a pas encore poussé assez loin, ni la Mechanique, ni l'Analyse, pour être en état de déterminer les mouvements dans les vibrations finies.* » Euler [1750, 69-70].

*tout à coup ; déterminer le mouvement vibratoire total, dont elle sera agitée.* »
Euler [1750, 70]

Que ce soit par ce mouvement à rebours qui remonte jusqu'à la courbe initiale ou en raison de la constitution de l'équation aux dérivées partielles en un problème mathématique autonome, la courbe initiale est intégrée à l'énoncé et aux solutions du problème. La courbe initiale est ainsi, pour eux, et à la différence de Bernoulli, une donnée du problème et surtout une solution de celui-ci.
Or, avec cette courbe initiale s'introduit l'expression d'une courbe arbitraire :

> « Réciproquement donc, dès qu'on connoit la figure initiale, qu'on aura donnée à la corde avant que de la relâcher, rien ne sera plus aisé que de décrire notre courbe infinie B'D'ADB'D & qui nous fera connoitre le mouvement, que la corde poursuivra. On tracera la courbe AMDB, égale & semblable à la figure initiale de la corde, & on en réïtérera la construction, tant vers la gauche au delà du point A, que vers la droite au delà du point B, alternativement au dessus & au dessous de l'axe, en sorte que partout les bouts qu'on a joints ensemble soient les mêmes. Cette construction a toujours lieu, de quelque nature que soit la figure initiale proposée de la courbe, & il ne s'agit que de la portion ADB; laquelle quand elle-même auroit d'autres continuations de part & d'autre en vertu de sa nature, elles n'entrent en aucune considération. Ainsi, si la figure ADB étoit un arc de cercle, sans se soucier de la continuation naturelle du cercle, on répétera la description de ce même arc de cercle ADB à l'infini alternativement au dessus & au dessous de l'axe ; & la même régle a toujours lieu, *de quelque nature que puisse être la figure initiale de la corde* [c'est moi qui souligne]. » Euler [1755, 216-7]

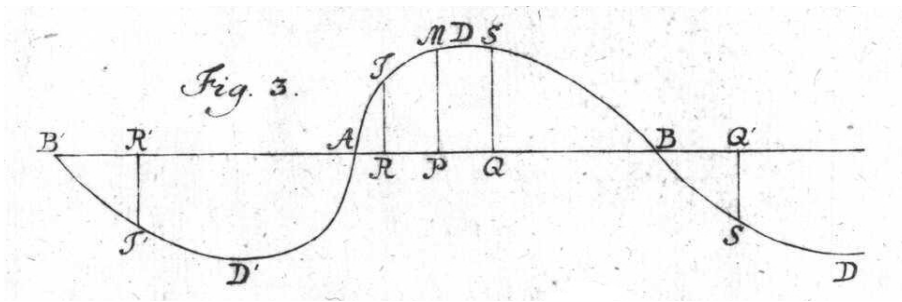

Pour D'Alembert et Euler la courbe initiale qui peut être une courbe quelconque fait partie des solutions. Dès lors, quand Bernoulli affirme que toutes les solutions sont données par des séries trigonométriques, cela implique, *pour D'Alembert et Euler*, qu'une courbe quelconque, quoi que chacun entende par cela, devrait pouvoir être exprimée par une série trigonométrique.
Après l'introduction d'une fonction Ψ quelconque dans l'expression des solutions de l'équation différentielle, la figure initiale de la corde introduit une seconde expression générale des courbes ou des fonctions (pour ceux qui, comme D'Alembert et Euler, réduisent le problème à une équation aux dérivées partielles). Un second lien est ainsi effectivement établi entre les solutions de ce problème et le fait de considérer la totalité des courbes ou des fonctions.

## 6 - Un énoncé inaugural sans inauguration

La controverse des cordes vibrantes offre l'exemple d'un énoncé inaugural qui n'a été énoncé que pour être dénoncé. Si l'on considère les critères retenus pour caractériser un texte inaugural, on constate qu'ils ne sont pas satisfaits voire même explicitement contestés. Seul un certain réalisme se retrouve chez les quatre mathématiciens. La dualité est avérée dans les mémoires de D'Alembert, Euler et Lagrange, mais ne l'est pas dans ceux de Bernoulli. Il n'y a d'inauguration qu'avec Bernoulli, mais il s'agit alors de défendre la validité d'un principe physique, ce qui ne correspond pas à l'acception donnée ici à l'inauguration. Les trois autres ne sauraient évidemment inaugurer un énoncé qu'ils contestent. De même la validité du principe de Bernoulli et les exemples d'expériences qu'il donne pour l'attester ne correspondent pas à l'acception de la conformité adoptée ici.

D'Alembert, Euler et Bernoulli ont chacun eu un rôle dans l'introduction de cet énoncé. D'Alembert a établi l'équation aux dérivées partielles du problème et en a exprimé la solution au moyen de l'expression d'une fonctions *quelconque*. Une expression a ainsi été donnée aux solutions qui ne fait plus intervenir les expressions analytiques reçues (algébriques ou transcendantes), et qui est indépendante de l'acception donnée par chacun à une fonction quelconque. Mais D'Alembert n'accorde qu'un pouvoir expressif limité à l'Analyse mathématique : tous les problèmes physiques n'ont pas selon lui nécessairement une expression mathématique adéquate[104] qu'il limite essentiellement aux développements en séries entières[105] [Jouve 2007, I, 144]. Il reste, à cette époque, attaché à ce cadre[106]. Euler a quant à lui repris l'expression de la solution donnée par

---

[104] Par exemple : « *On objectera peut-être, qu'il est impossible d'expliquer dans ma théorie, pourquoi la corde frappée d'une maniere quelconque, rend toujours à peu près le même son ; puisque ses vibrations, selon moi, peuvent être très-irrégulieres en plusieurs cas [163] . J'en conviens ; mais je suis persuadé que la solution de cette question n'appartient point à l'Analyse : elle a fait tout ce qu'on étoit en droit d'attendre d'elle ; c'est à la Physique à se charger du reste.* » D'Alembert [1761, 40 ; Jouve 2007, II, 171] ; « *l'expérience & la théorie ne sauroient être d'accord dans la détermination du mouvement des cordes vibrantes, ne fût-ce que par cette seule raison, que le mouvement d'une corde vibrante cesse bientôt, quoique par la théorie il doive durer continuellement.* » D'Alembert [1761, 67 ; Jouve 2007, II, 205]. Un autre exemple du refus de mélange des genre : « *Avant que de répondre à ce Géometre [Daniel Bernoulli], nous allons rendre sa proposition encore plus générale. Nous prouverons qu'en effet la somme moyenne entre toutes les sommes de la serie cos .x + cos .2x + cos .3x, &c. à l'infini = − ½ ; mais nous le prouverons rigoureusement, & non par le calcul probabilités, dont l'usage, nous l'osons dire, est dans le cas présent tout-à-fait chimérique. Il ne s'agit pas ici de conjecturer, mais de démontrer ; & il seroit dangereux, (quoiqu'à la vérité ce malheur soit peu à craindre) qu'un genre de démonstration si singulier s'introduisît en géométrie. Ce qui pourra seulement paroître surprenant, c'est que de pareils raisonnemens soient employés comme démonstratifs par un Mathématicien célèbre, & dans un écrit où il s'exprime avec très-peu de ménagement sur le dixiéme Mémoire de mes Opuscules, concernant la théorie des probabilités ; théorie qui n'a pourtant (ce me semble) rien de plus choquant que l'usage qu'il fait ici de l'analyse des jeux de hasard. Quoi qu'il en soit, après avoir prouvé, non par cette analyse, mais par un calcul rigoureux, que la somme moyenne dont il s'agit, est en effet = − ½ dans tous les cas, nous prouverons en suite qu'on n'en peut rien conclure pour la somme réelle de la serie, que cette fraction − ½ ne représente point.* » D'Alembert [1761, 157-8 ; Jouve 2007, II, 244-5]

[105] Par exemple, quand une fonction est impaire, D'Alembert infère qu'elle doit n'avoir que des puissances impaires de x [D'Alembert 1761, 14].

[106] Par exemple : « *le mouvement de la corde ne peut être soumis à aucun calcul analytique, ni*

D'Alembert en y intégrant des « *fonctions irrégulières* » dont il trouve dans la courbe tracée, à défaut d'expression analytique, une expression géométrique. Cette extension à des « *fonctions irrégulières* » est récurrente avec les équations aux dérivées partielles [Euler 1767 ; Dhombres 1988]. Elle est aussi un thème récurrent dans la correspondance d'Euler et un terrain d'entente entre Euler et Lagrange contre Bernoulli mais aussi contre D'Alembert[107]. Mais s'il accepte cette extension de l'Analyse au-delà des expressions analytiques reconnues ou existantes, Euler n'a pas pour autant de système d'expressions de substitution à proposer susceptible d'exprimer *toutes* les « *fonctions irrégulières* ». Il s'accommode de la nécessité d'appliquer l'Analyse mathématique à des courbes qui n'ont pas d'expression analytique. L'Analyse mathématique a dès lors une portée géométrique (courbes) qui excède ses moyens d'expression analytiques[108].

En exprimant le problème d'une corde vibrante par une équation aux dérivées partielles D'Alembert fait de la courbe initiale une solution du problème physique. Or avec cette courbe initiale s'introduit à nouveau l'expression d'une courbe quelconque. De ce fait, mais aussi à cette condition, prétendre donner toutes les solutions du problème de la corde vibrante revient à prétendre donner une représentation de *toutes* les fonctions. Daniel Bernoulli est celui qui dans ce contexte va déclencher l'introduction de l'énoncé en affirmant, à partir d'un principe physique sur la décomposition des sons en modes propres, que les séries trigonométriques tenues par D'Alembert et Euler pour des solutions particulières sont en fait *toutes* les solutions du problème de la corde vibrante. La présence de l'expression d'une fonction quelconque dans la solution des relations différentielles et surtout l'expression de la courbe initiale qui pour D'Alembert comme pour Euler doit être intégrée à l'énoncé du problème les conduisent à récuser l'affirmation de Bernoulli et à formuler l'énoncé inaugural qui en résulterait si elle était vérifiée. L'introduction de cet énoncé est la conséquence,

---

*représenté par aucune construction, quand la courbure fait un saut en quelque point* » [D'Alembert 1761, 22 ; Jouve 2007, II, 148].

107 Ainsi, Lagrange fait valoir auprès d'Euler l'intérêt de la nouvelle approche qu'il a trouvée au problème de la corde vibrante : « *Cette solution, étant d'un genre tout à fait nouveau, ne sera peut-être pas aussi indigne de votre attention, et elle servira encore plus à établir l'usage des fonctions irrégulières et discontinues dans une infinité d'autres problèmes.* » Lagrange à Euler, 24 novembre 1759, Oeuvres de Lagrange, vol. 14, 171-2. Dans une autre lettre, Euler réfute une objection de D'Alembert contre les fonctions arbitraires : « *c'est une propriété essentielle des équations différentielles à trois et à plusieurs variables, que leurs intégrales renferment des fonctions arbitraires qui peuvent aussi bien être discontinues que continues.* » Euler à Lagrange, 16 février 1765, Oeuvres de Lagrange, 14, 203. On peut encore citer Lagrange : « *D'ailleurs, les phénomènes de la propagation du son ne peuvent s'expliquer qu'en admettant les fonctions discontinues, comme je l'ai prouvé dans ma seconde dissertation.* » Lagrange à D'Alembert, 20 mars 1765, Oeuvres de Lagrange, vol. 13, 38. Voir encore Lagrange 1760, Oeuvres, 1, 158.

108 C'est là typiquement un exemple de *problème d'expression* en l'occurrence consécutif à la considération de nouveaux problèmes exprimés au moyen d'équations aux dérivées partielles. C'est un point déjà abordé par Jean Dhombres [1988]. Le problème se trouve renforcé quand la considération de grandeurs imaginaires met aussi en défaut l'expression par une courbe : « *Mais, si l'on demandait une semblable intégrale complète pour le cas où a serait une quantité négative -b, je ne vois pas comment on la pourrait représenter par des courbes arbitraires, puisqu'on n'y saurait assigner les appliquées qui répondent à des abscisses imaginaires.* » Euler à Lagrange, 9 novembre 1762, Oeuvres de Lagrange, vol. 14, 203. L'expression de la fonction quelconque pourrait aussi apparaître dans l'équation différentielle, et non seulement dans l'énoncé du problème comme condition initiale. C'est ce qui arrive avec le cas d'une corde d'épaisseur non uniforme.

dans ce contexte et avec ces conditions, de la forme des solutions proposées par Bernoulli et de la conviction de celui-ci qu'il a ainsi donné *toutes* les solutions. Mais si Euler formule cet énoncé, il ne le soutient pas. Son mémoire publié en 1755 est ainsi largement consacré à le dénoncer. Bernoulli serait le seul à le soutenir, mais il ne l'a pas énoncé. Il se doit néanmoins de répondre aux objections. Ses réponses vont surtout consister en une défense de sa solution, et non en un soutien de l'énoncé inaugural, qu'il va considérer comme un théorème d'interpolation dont il esquissera la démonstration [Bernoulli 1758, 165]. Ainsi, son « soutien » se fait soit de manière indirecte sous la forme d'objections aux objections d'Euler soit en l'assimilant à un théorème qu'il s'agirait de démontrer. Si Euler a bien écrit un mémoire de type anti-inaugural pour réfuter cet énoncé, il n'y a pas parmi les textes cités de Bernoulli de texte inaugural écrit pour le soutenir. Aucun de ces mathématiciens finalement ne soutient l'énoncé : celui qui l'énonce le réfute, celui qui pourrait le soutenir, ni ne le soutient ni ne l'énonce[109]. Ainsi, cette controverse donne l'exemple d'un énoncé inaugural qui n'a été introduit que pour être dénoncé. Elle aura permis de comparer un texte qui défend un principe physique et un autre, celui de Fourier, qui soutient un énoncé inaugural. Cela confirme leur différence et donne un exemple de principe physique là où l'on pouvait s'attendre, compte tenu des citations données, à trouver un texte inaugural[110].

## 7 - Lame chauffée et corde vibrante

L'identification des expressions de la totalité des courbes et des fonctions dans *La théorie de la chaleur* de Fourier et dans la controverse des cordes vibrantes met en évidence la similarité du rôle de la base de la lame chez Fourier et de l'état initial de la corde vibrante. Elles sont l'une et l'autre l'expression d'une courbe quelconque. Leurs caractéristiques sémiotiques sont évidemment différentes : la représentation d'une courbe par la forme initiale de la corde est nettement iconique, ce qui n'est pas le cas de la base de la lame, etc. Mais l'une et l'autre sont néanmoins considérées comme des expressions d'une courbe quelconque. De même, la lame et la corde vibrante jouent le rôle d'un dispositif mettant en relation l'expression d'une courbe quelconque et le développement en série trigonométrique. Ainsi, dans un cas comme dans l'autre, la *représentation* ne fait pas seulement référence à l'objet qu'elle représente, à la lame chauffée ou à la corde qui vibre, son expression, en tant que telle, intervient comme dispositif qui a une incidence sur les développements mathématiques considérés, en l'occurrence dans l'introduction d'une nouvelle représentation des fonctions. Les mêmes conditions sémiotiques qui permettent l'introduction de la représentation par des séries trigonométriques se retrouvent dans les deux cas. On peut ainsi rendre compte sémiotiquement de cette introduction. Mais il en résulte aussi que cette analyse sémiotique *ne* peut *pas* rendre compte du fait que l'énoncé inaugural ait été *soutenu* dans un cas et pas dans l'autre : rien sémiotiquement n'empêcherait qu'il le soit à partir du problème de la corde vibrante. Fourier dispose certes de la formule intégrale, et il en a fait l'argument principal de son inauguration, mais cette formule n'est pas nécessaire à l'inauguration, d'autres arguments pourraient être donnés. En revanche, les problèmes sémiotiques posés sont eux incontournables. Or, on a vu Fourier, au contraire de Bernoulli, s'y confronter. Il y a bien là une différence objective.

---

[109] Citons la conclusion de l'analyse sur les cordes vibrantes de Darrigol [2007, 385] : « *To summarize, Bernoulli, d'Alembert, Euler, and Lagrange held different positions regarding the permissible string curves, the curves that trigonometric series could represent over a finite interval, and the status of partial vibrations. Bernoulli admitted any string curve for which both the ordinate and the radius of curvature remained very small compared to the length of the string; he believed that trigonometric series could represent any such curve; he asserted the physical existence of the partial vibrations. D'Alembert required the string curve to be analytic and close to the axis; he believed that any "discontinuous" curve and even some "continuous" curves could not be represented by trigonometric series; he regarded partial vibrations as mathematical fictions. Euler admitted any continuous string curve with piecewise continuous slope and curvature, and with small ordinate and slope; he denied that trigonometric series could represent non-analytic curves, at least those which coincide with segments of the axis (pulses); he ascribed some physical reality to partial vibrations of non-necessarily sine form. Lagrange originally admitted the same string curves as Euler (except for polygonal curves), but came to believe that his passage from the discrete to the continuous required that all derivatives should be finite; he believed that trigonometric series could in some "asymptotic" sense represent any curve, but sometimes denied perfect identity between the series and the curve when the latter was non-analytic or pulse-like; like d'Alembert, he denied any physical reality of th partial vibrations.* »

[110] La confusion de l'énoncé imputé à Bernoulli avec celui de Fourier a déjà été clairement dénoncée par Burkhardt [1908, 20 ; Jouve 2007, II, 29].

Le rapprochement de ces deux analyses a fait ressortir les similitudes sémiotiques de la lame et de la corde vibrante ainsi que le rôle de l'état initial dans l'introduction de la représentation d'une fonction quelconque par des séries trigonométriques. C'est en tant qu'expressions que la lame et la corde vibrante interviennent dans l'introduction de cette représentation. Mais ce ne sont évidemment pas des « expressions pures », des expressions qui « parleraient d'elles-mêmes », ce sont des expressions elles-mêmes interprétées. Or, l'expression de la corde vibrante n'est pas la même pour Bernoulli d'une part et pour D'Alembert et Euler d'autre part : pour D'Alembert et Euler l'état initial en fait pas partie, pas pour Bernoulli (c'est pour cela qu'il n'introduit pas lui-même l'énoncé inaugural). Mais on a aussi vu que le statut de la condition initiale dépendait du statut donné à l'équation différentielle associée au phénomène décrit. Ce qui veut dire que la représentation d'une corde vibrante, ou l'interprétation qui en est faite, dépend du statut que l'on attribue à l'équation différentielle : l'équation différentielle change l'interprétation de l'expression d'une corde vibrante.

# IV - Conclusion

Les textes inauguraux sont nécessairement assez rares. Ils marquent des moments toujours singuliers et remarquables par la constitution des moyens d'expression de totalités, et donc de la généralité, propres aux mathématiques. Ainsi, les séries trigonométriques offrent à Fourier une représentation de toutes les courbes et de toutes les fonctions à une époque où les séries entières, et depuis plus longtemps encore les polynômes, ne pouvaient plus prétendre jouer ce rôle. Des théorèmes, comme celui de Dirichlet [1829], et des théories, en premier lieu la théorie de la chaleur de Fourier, mais aussi rapidement la théorie de l'intégration, vont ainsi pouvoir être développés en tirant parti de ces nouveaux moyens d'expression. Les critères qui ont été proposés pour caractériser ces textes et les énoncés inauguraux aident d'abord à les repérer. Il a été ainsi possible de vérifier que la *Théorie analytique de la chaleur* de Joseph Fourier les satisfaisait. Mais il importe surtout de reconnaître la fonction inaugurale de ces textes, c'est-à-dire à la fois leur rôle dans l'introduction de nouvelles représentations, les problèmes sémiotiques qui déterminent en partie leurs caractéristiques, et enfin les développements qu'ils rendent ensuite possibles. Le fait d'établir que le texte de Fourier est, comme *La Géométrie* de Descartes, un texte inaugural établit à la fois l'importance et la récurrence de la *conformité* dans l'histoire des mathématiques. Non seulement la *généralisation* n'est peut-être pas le vecteur principal du développement des mathématiques mais, à l'inverse, la *conformité* apparaît comme un moyen historique privilégié d'introduction des conditions d'expression de la *généralité* en mathématiques. Certains des changements les plus importants en mathématiques ont été faits par *conformité.*

La notion de texte inaugural contribue à faire reconnaître le rôle des *problèmes sémiotiques* en mathématiques et *a fortiori* en histoire des mathématiques. Ces problèmes découlent du fait que les mathématiques sont tributaires de leurs moyens d'expression : le mathématicien veut exprimer quelque chose, parfois en partie suggérée par les moyens d'expression qu'il utilise, sans toujours le pouvoir avec ces moyens. Il a alors un problème d'expression. Les mathématiques sont ainsi confrontées à des problèmes qui ne sont pas que les problèmes de mathématiques reconnus comme tels et la constitution de ces moyens d'expression à partir de ceux disponibles est un élément de l'histoire des mathématiques.

Les énoncés et en particulier les raisonnements mathématiques tirent parti de leurs moyens d'expression : ils s'inscrivent dans des conditions sémiotiques qui font partie de leurs conditions de possibilité. Chaque énoncé devient ainsi l'indice de ses conditions de possibilité. La généralité, en particulier, suppose des moyens d'expression qui la rendent possible. Par une sorte de principe de conservation sémiotique les moyens d'expression ne varient ensuite que s'ils sont modifiés. Autrement dit : les changements dans les caractéristiques sémiotiques des textes mathématiques doivent être rapportés à des causes assignables.

Nous avons pu ainsi déterminer à la fois le rôle de la représentation d'une fonction arbitraire par des séries trigonométriques dans la *généralité de la théorie de la chaleur* et le rôle de cette théorie, et plus précisément de l'étude de la lame, dans l'introduction de cette représentation. L'examen des conditions d'expression de la généralité de la *théorie* a d'abord permis d'établir qu'elle *ne* tirait *pas* sa généralité de l'*équation générale*. La contribution des séries trigonométriques à la généralité de cette théorie *ne* vient donc *pas* de leur rôle dans la résolution de l'*équation générale* mais de la possibilité, introduite par leur expression, de considérer des solides avec des *conditions initiales* qui soient à la fois *quelconques* mais dont l'expression permet aussi dans chaque cas, quoi que de manière *ad hoc*, de *déterminer les coefficients de la solution générale* et ainsi de trouver la distribution de chaleur cherchée.

La représentation des fonctions par des séries trigonométriques a ainsi été à la fois *inaugurée* au cours de l'étude de la distribution de la chaleur dans un solide particulier (en l'occurrence, une lame) et *exploitée* dans l'étude de la distribution de la chaleur des autres solides. On a ainsi un exemple de texte inaugural qui dispose de la représentation qu'il inaugure. Mais il a bien d'abord fallu que Fourier l'inaugure pour pouvoir en disposer ainsi. Il a dû pour cela résoudre les problèmes sémiotiques attendus et qui confèrent aux textes inauguraux leurs caractéristiques. Ne dénoncer qu'un manque de rigueur dans les arguments donnés pour soutenir un énoncé inaugural revient à ignorer les conditions sémiotiques de cette rigueur et en particulier les conditions de la généralité des énoncés et des démonstrations mathématiques. A l'inverse, l'existence de ces conditions reconnue, les commentaires sur les textes qui les modifient sont souvent des révélateurs privilégiés aussi bien de ces conditions que du rôle joué ensuite par le fait de ne pas les prendre en compte.

L'énoncé inaugural des séries trigonométriques avait déjà été introduit au milieu du siècle précédent lors de la controverse sur les cordes vibrantes. C'est un fait bien connu, notamment de Fourier et de ses interlocuteurs. L'énoncé n'avait cependant pas été soutenu (suivant notre acception) et *de fait* Fourier n'a pu disposer de cette représentation sans l'inaugurer lui-même, ce dont les mathématiciens (et les historiens) purent en revanche se dispenser après lui. La controverse des cordes vibrantes réunit donc, comme le texte de Fourier, les conditions de l'introduction de l'énoncé inaugural, mais sans présenter cependant, contrairement au texte de Fourier, les caractéristiques d'un texte inaugural. Il a à nouveau été possible de rendre compte de l'introduction de cet énoncé en déterminant les conditions de la rencontre de l'expression d'une fonction générale et de celle des séries trigonométriques. On a ainsi vu que l'expression des séries trigonométriques avait été introduite par Bernoulli à partir du principe physique des systèmes oscillants qu'il défendait, les deux expressions générales d'une fonction l'ayant en revanche été ultérieurement par D'Alembert et Euler. L'une d'elles est l'expression $\Psi$ d'une fonction quelconque (expression analytique) qui figure dans la solution générale de l'équation différentielle introduite à cette occasion par D'Alembert et reprise par Euler. L'autre est la figure initiale arbitraire de la corde (expression géométrique), mais dont le statut, et en particulier le fait de la considérer elle-même comme faisant partie des solutions de l'équation différentielle, dépend de la conformité accordée à l'équation différentielle par rapport au problème d'une corde vibrante, qui diffère alors selon les protagonistes. Cette courbe initiale qui, pour D'Alembert et Euler, fait *nécessairement* partie des solutions (elle est même *nécessaire* à la résolution de l'équation différentielle) n'en faisait pas partie pour Bernoulli. C'est ainsi seulement avec l'introduction de l'une ou l'autre de ces deux expressions que les séries trigonométriques données par Bernoulli pour solution du problème de la corde vibrante apparaissent, pour D'Alembert et Euler, devoir représenter une fonction quelconque. Le problème de la corde vibrante est bien ainsi le lieu de rencontre des séries trigonométriques et de deux expressions d'une fonction arbitraire mais sans que les conditions de cette rencontre ne soient partagées par ces trois protagonistes. Ainsi, les expressions de la généralité qui complètent les conditions pour que l'énoncé inaugural soit formulé sont introduites par D'Alembert et Euler qui, de fait, et pour des raisons différentes, ne le soutiennent pas.

On a pu ainsi suivre, voir se mêler et se propager les expressions de la généralité aussi bien dans la théorie de la chaleur que dans l'étude d'une corde vibrante. On a vu en particulier le rôle de l'expression de l'état initial par laquelle, dans les deux cas, peut s'introduire l'expression d'une fonction arbitraire qui sera ensuite mise en relation avec l'expression d'une série trigonométrique. Mais c'est en revanche seulement dans la théorie de la chaleur de Fourier que les séries trigonométriques deviennent elles-mêmes inversement le moyen d'expression de la généralité des conditions initiales et qu'elles deviennent ainsi, et seules, constitutives de la généralité de sa théorie. La comparaison de la controverse des cordes vibrantes et de la *Théorie analytique de la chaleur* fait ainsi ressortir la différence entre les conditions d'introduction d'une nouvelle forme d'expression de la généralité, présentes dans les deux cas, et son usage comme expression de la généralité après que Fourier l'ait inaugurée.

   Ces analyses se fondent sur le caractère objectif des moyens par lesquels la généralité est exprimée et sur la nécessité de disposer de tels moyens : tout énoncé général suppose et dépend des moyens d'expression qui rendent possible l'expression de cette généralité. A la manifestation d'une forme de généralité dans un texte doit répondre les moyens d'expression correspondants. Les énoncés et les raisonnements ont leurs conditions de possibilité ; celles-ci ne sont que très partiellement logiques. Ces moyens doivent aussi avoir été introduits d'une manière ou d'une autre. Ainsi, à l'analyse des moyens d'expression utilisés dans un texte répond l'étude historique de leur introduction.

Cette double nécessité (existence et introduction) est le point de départ de l'étude des moyens d'expression utilisés. Ceux-ci peuvent ensuite être repérés, séparés les uns des autres, identifiés et en quelque sorte tracés dans un texte ou un corpus. Il est possible de les prendre comme objets de raisonnement, quasiment de calcul. Il devient possible de se demander comment ils apparaissent, se propagent et disparaissent. L'analyse présentée a simplement consisté à déterminer dans les textes considérés les moyens d'expression de la généralité utilisés et leurs relations aux séries trigonométriques. L'inauguration de ces séries est ainsi apparue comme la condition pour que leur expression puisse devenir elle-même un moyen d'expression général.

# Bibliographie

différentielles partielles et à coefficients constants". *Journal de l'École Polytechnique. Paris.* 12, 511-592, 1823 repris in Cauchy 1905.

Cauchy, Augustin-Louis. "Mémoire sur les intégrales définies". *Mémoires de l'Académie Royale des Sciences*, 1827 repris in Cauchy 1882.

Cauchy, Augustin-Louis. *Œuvres complètes*, Série 1, tome 1. Paris : Gauthier-Villars, 1882.

Cauchy, Augustin-Louis. *Œuvres complètes*, Série 2, tome 1. Paris : Gauthier-Villars, 1905.

Cauchy, Augustin-Louis. *Œuvres complètes*, Série 2, tome 4. Paris : Gauthier-Villars, 1899.

Charbonneau, Louis. *L'Œuvre mathématique de Joseph Fourier, 2 tomes*: EHESS, 1976.

D'Alembert, Jean Le Rond (1747a). *Réflexions sur la cause générale des vents*. Paris : David l'aîné (Jean-Baptiste Coignard, imprimeur).

D'Alembert, Jean le Rond. "Recherches sur la courbe que forme une corde tendue mise en vibration", *Histoire de l'Académie Royale des Sciences et des Belles-Lettres de Berlin*, 3, pp. 214-219, (1747b) 1749.

D'Alembert, Jean le Rond. "Recherches sur la courbe que forme une corde tenduë mise en vibration : Suite", *Histoire de l'Académie Royale des Sciences et des Belles-Lettres de Berlin*, pp. 220-249, (1747c) 1749.

D'Alembert, Jean le Rond. "Recherches sur les vibrations des cordes sonores" *in Opuscules*, t. I, pp. 1-73, 1761.

D'Alembert, Jean le Rond. *Opuscules mathématiques ou mémoires sur différents sujets de Géométrie, de méchanique, d'optique, d'astronomie etc.*, 1761.

D'Alembert, Jean le Rond. *Opuscules Mathématiques*, tome IV, 1768.

D'Alembert, Jean Le Rond & Crépel, Pierre (éd) & Guilbaud, Alexandre (éd) & Jouve, Guillaume (éd). *OEuvres complètes de D'Alembert : "Opuscules mathématiques" tome I (1761). volume III/1*. Paris : CNRS Editions, 2008.

Darrigol, Olivier. "The acoustic origins of harmonic analysis", *Archive for History of Exact Sciences*, 61 (4), pp. 343-424, 2007.

Demidov, Sergei S. "D'Alembert et la naissance de la théorie des équations différentielles aux dérivées partielles" in (Emery & Monzani 1989) , 333-350,